\documentclass[11pt,letterpaper,reqno]{amsart}
\usepackage{amsfonts}
\usepackage{amsmath}
\usepackage{amsthm, amscd}
\usepackage{mathtools}
\usepackage{amssymb}
\usepackage{mathrsfs}
\usepackage{graphicx}
\usepackage{caption}
\usepackage{subcaption}
\usepackage{float}
\usepackage{xypic}
\usepackage[abs]{overpic}
\usepackage[alphabetic]{amsrefs}
\usepackage{etex}
\usepackage{tikz-cd}
\usepackage{enumerate}
\usepackage{xcolor}

\usepackage{hyperref}
\hypersetup{
	colorlinks,
	linkcolor=[rgb]{0,0,0.7},
	urlcolor=[rgb]{0,0,0.4},
	citecolor=[rgb]{0.4,0.1,0}
}

\allowdisplaybreaks

\theoremstyle{plain}\newtheorem{Theorem}{Theorem}[section]
\theoremstyle{plain}\newtheorem{Corollary}[Theorem]{Corollary}
\theoremstyle{plain}\newtheorem{Lemma}[Theorem]{Lemma}
\theoremstyle{plain}\newtheorem{Definition}[Theorem]{Definition}
\theoremstyle{plain}\newtheorem{Proposition}[Theorem]{Proposition}
\theoremstyle{plain}
\theoremstyle{plain}
\theoremstyle{plain}\newtheorem*{Claim*}{Claim}
\theoremstyle{plain}\newtheorem*{Theorem*}{Theorem}
\theoremstyle{plain}\newtheorem{Question}[Theorem]{Question}

\theoremstyle{remark}\newtheorem{remark}[Theorem]{Remark}
\theoremstyle{remark}\newtheorem{Example}[Theorem]{Example}
\theoremstyle{remark}\newtheorem*{Notation*}{Notation}

\theoremstyle{plain}
\makeatletter
\newtheorem*{rep@theorem}{\rep@title}
\newcommand{\newreptheorem}[2]{
\newenvironment{rep#1}[1]{
 \def\rep@title{#2 \ref{##1}}
 \begin{rep@theorem}}
 {\end{rep@theorem}}}
\makeatother
\newreptheorem{Theorem}{Theorem}
\newreptheorem{Proposition}{Proposition}
\newreptheorem{Corollary}{Corollary}

\numberwithin{equation}{section}

\usepackage{lipsum}

\DeclareMathOperator{\Imm}{Im}

\DeclareMathOperator{\SL}{SL}

\DeclareMathOperator{\pt}{pt}
\DeclareMathOperator{\id}{id}

\DeclareMathOperator{\Diff}{Diff}

\DeclareMathOperator{\Emb}{Emb}

\DeclareMathOperator{\Mcg}{MCG}
\DeclareMathOperator{\Aut}{Aut}
\DeclareMathOperator{\inte}{int}

\DeclareMathOperator{\Fr}{Fr}

\DeclareMathOperator{\Dax}{Dax}

\newcommand{\bZ}{\mathbb{Z}}

\newcommand{\cM}{\mathcal{M}}

\newcommand{\cS}{\mathcal{S}}

\newcommand{\PSL}{\mbox{PSL}}

\author{Jianfeng Lin}
\address{Yau Mathematical Sciences Center, Tsinghua University, Beijing 100084, China}
\email{linjian5477@mail.tsinghua.edu.cn}
\author{Weiwei Wu}
\address{Zhejiang University, Hangzhou, Zhejiang 310058, China}
\email{wwwu9693@zju.edu.cn}
\author{Yi Xie}
\address{Beijing International Center for Mathematical Research, Peking University, Beijing 100871, China}
\email{yixie@pku.edu.cn}
\author{Boyu Zhang}
\address{Department of Mathematics, The University of Maryland at College Park, Maryland 20742, USA}
\email{bzh@umd.edu}

\title[Dax invariants, light bulbs, and symplectic structures]{Dax invariants, light bulbs, and isotopies of symplectic structures}

\begin{document}

\maketitle

\begin{abstract}
This paper addresses several isotopy problems on $4$-manifolds.  First, we classify the isotopy classes of embeddings of $\Sigma$ in $\Sigma\times S^2$ that are geometrically dual to $\{\mbox{pt}\}\times S^2$, where $\Sigma$ is a closed oriented surface with a positive genus, and show that there exist infinitely many such embeddings that are homotopic to each other but mutually non-isotopic, thereby answering a question of Gabai \cite[Question 10.15]{Gabai20}.

By combining this construction with techniques from symplectic topology, we also answer Problem 2(a) in McDuff--Salamon's problem list \cite{MSIntro} and a question of Cieliebak--Eliashberg--Mishachev \cite{CieliebakHprinciple}, which concern the uniqueness and $h$--principle of symplectic structures on closed $4$--manifolds. We answer these questions by establishing the following results:
(1) The space of symplectic forms on every irrational ruled surface homologous to a fixed symplectic form has infinitely many connected components;
(2) There exist infinitely many symplectic forms on every irrational ruled surface that are formally homotopic, cohomologous, but not homotopic to each other.
Both are the first such examples for closed $4$--manifolds. 

The proofs are based on a generalization of the Dax invariant to embedded closed surfaces. In the course of the proof, we obtain several properties of the smooth mapping class group of $\Sigma\times S^2$, which may be of independent interest. For example, we show that there exists a surjective homomorphism from $\pi_0\operatorname{Diff}(\Sigma\times S^2)$ to $\mathbb{Z}^\infty$, such that its restriction to the subgroup of elements pseudo-isotopic to the identity is of infinite rank. 
\end{abstract}

\section{Introduction}
\label{sec_intro}
\subsection{Lightbulbs and the mapping class group of $\Sigma\times S^2$}
Let $\Sigma$ be a closed oriented surface with positive genus, let $M=\Sigma\times S^2$. The first main result of the paper is the following theorem.

\begin{Theorem}
\label{thm_homot_non_isotopic_intro}
    There exists an infinite collection of embeddings of $\Sigma$ in $M$ that are all geometrically dual to $G:=\{\pt\}\times S^2$, are homotopic to each other relative to a neighborhood of $G$, and are mutually non-isotopic as embedded surfaces in $M$. 
\end{Theorem}

In fact, we will also give a complete classification of isotopy classes of embeddings of $\Sigma$ in $\Sigma\times S^2$ that are geometrically dual to $\{\pt\}\times S^2$, which will be stated in Theorem \ref{intro thm: classification} below. We remark that in Theorem \ref{thm_homot_non_isotopic_intro}, the surfaces are still non-isotopic after reparametrizations, and the isotopies are not required to fix the dual sphere.

Theorem \ref{thm_homot_non_isotopic_intro} is closely related to Gabai's 4-dimensional lightbulb theorem.
Given an embedded surface $R$ with a geometrically dual embedded sphere $G$ in a 4-manifold $X$, the surface
$R$ is called $G$-\emph{inessential} if the induced map $\pi_1(R\setminus G)\to \pi_1(X\setminus G)$  is trivial \cite[Definition 9.6]{Gabai20}. The following is a version
of Gabai's $4$-dimensional light bulb theorem for embedded surfaces.
\begin{Theorem*}[{\cite[Theorem 10.1]{Gabai20}}]
Suppose $X$ is an oriented $4$-manifold and
$\pi_1(X)$ has no $2$-torsion.
Let $R$ and $S$ be two homotopic closed embedded surfaces in $M$ with a common geometrically dual sphere $G$.
If both $R$ and $S$ are $G$-inessential, then $R$ and $S$
are isotopic.
\end{Theorem*}
In \cite[Question 10.15]{Gabai20}, Gabai asked whether
the ``$G$-inessential'' condition in the above theorem can be removed. Since $\pi_1(\Sigma\times S^2)$ has no 2-torsion, Theorem \ref{thm_homot_non_isotopic_intro} gives a negative answer to Gabai's question.

Example \ref{example_self-referential} below describes the surfaces constructed in Theorem \ref{thm_homot_non_isotopic_intro}.

\begin{figure}
	\begin{overpic}[width=0.27\textwidth]{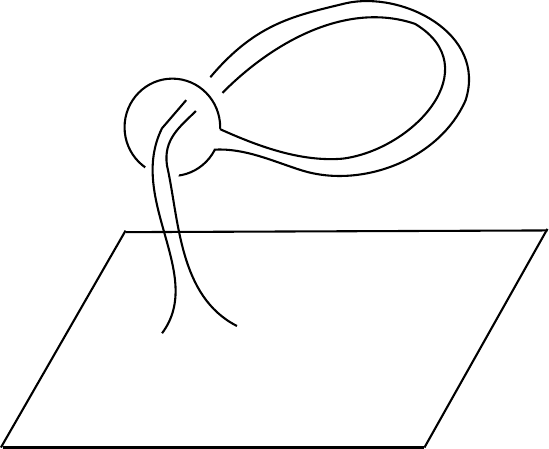}
		\put(16,67){\small{$\partial D$}}
		\put(60,15){$\Sigma_0$}
        \put(85,65){\small{$\gamma_i$}}
	\end{overpic}
	\caption{Example \ref{example_self-referential}}
    \label{fig_self-referential_surface}
\end{figure}

\begin{Example}
\label{example_self-referential}
Take $b_0\in \Sigma, b_1\in S^2$. Consider the embedded surfaces $\Sigma_0=\Sigma\times \{b_1\}$ and $S_0=\{b_0\}\times S^2$. Take an embedded disk $D:D^{3}\hookrightarrow M$ disjoint from $\Sigma_0$ and $S_0$. Take a sequence of embedded arcs $\{\gamma_{i}\}_{i\geq 1}$ in $M$ that satisfy the following conditions:
\begin{enumerate}
    \item $\gamma_{i}(0)=(b_0,b_1), \gamma_{i}(1)\in \partial D$;
    \item The interior of $\gamma_i$ is disjoint from $\Sigma_0, S_0$ and intersects $D$ transversely at a single point $\gamma_{i}(\frac{1}{2})\in \operatorname{int}(D)$.
    \item Take any path $\eta_{i}:I\to D$ from $\gamma(1)$ to $\gamma(\frac{1}{2})$. Set $\gamma'_{i}=\gamma_{i}|_{[\frac{1}{2},1]}* \eta_{i}$. Then the loops $\{\gamma'_{i}\}_{i\in \mathbb{N}}$ are all non-contractible and not freely homotopic to each other and their inverses.  
\end{enumerate}
We take a tube along $\gamma_i$ and form the connected sum $\Sigma_{i}:=\Sigma_0\sharp_{\gamma_i} \partial D\hookrightarrow M$. 
See Figure \ref{fig_self-referential_surface}. 
Then the surfaces $\{\Sigma_{i}\}_{i\geq 0}$ all have a common geometric dual and are homotopic to each other. We will prove that $\{\Sigma_{i}\}_{i\geq 0}$ are not smoothly isotopic to each other. We will also prove that $\Sigma_i$ is isotopic to $f_{i}(\Sigma_0)$, where $f_{i}:M\to M$ is the so-called self-referential barbell diffeomorphism (see Section \ref{sec_prelim}), supported in the regular neighborhood 
\[\nu(S_0\cup \gamma_{i}\cup \partial D)\cong (S^2\times D^2)\natural (S^2\times D^2).\] These diffeomorphisms are all pseudo-isotopic to the identity. It will be proved that they are linearly independent in $\Mcg(M)$ and $\Mcg_{PI}(M)$.
\end{Example}

During the proof, we also establish a theorem on the mapping class group of $M = \Sigma\times S^2$. 
Let $\Diff(M)$ denote the group of diffeomorphisms on $M$, equipped with the $C^{\infty}$ topology. 
Define
\begin{align*}
    \Mcg(M) &= \pi_0\Diff(M),\\
    \Mcg_{PI}(M) &= \{[f]\in \Mcg(M)\mid f \text{ is pseudo-isotopic to }\id\}.
\end{align*}
Let $\mathbb{Z}^{\infty}$ be the direct sum of countably infinite copies of $\mathbb{Z}$.  Then we have the following result.
\begin{Theorem}\label{thm: MCG infinitely generated}
There exists a surjective homomorphism from $\Mcg(M)$ to $\mathbb{Z}^{\infty}$ such that its restriction to $\Mcg_{PI}(M)$ is also of infinite rank.
\end{Theorem}

As an immediate corollary, every subgroup of $\Mcg(M)$ that contains $\Mcg_{PI}(M)$ is infinitely generated. For example, the Torelli group of $M$ is infinitely generated. 

Results on $\Mcg_{PI}(X)$ for various 4-manifolds $X$ were obtained in recent years by 
Budney--Gabai \cite{BG2019,BG2023}, Watanabe \cite{Watanabe2020}, Igusa \cite{Igusa2021}, Singh \cite{Singh2021} and  Fern\'{a}ndez--Gay--Hartman-Kosanovi\'{c} \cite{FGHK2024}. Theorem \ref{thm: MCG infinitely generated} adds to the list of known examples of closed $4$--manifolds whose $\Mcg_{PI}$ is non-trivial. 

\subsection{Nonisotopic symplectic structures on ruled surfaces}
A fundamental question in symplectic topology is to study the existence and uniqueness of symplectic structures on a given smooth manifold.  
Let $X$ be a closed 4-manifold and let $a\in H^{2}(X;\mathbb{R})$ be a class. Consider the space of symplectic forms
\[
\mathcal{S}_{a}:=\{\text{symplectic form $\omega$ on $X$ with $[\omega]=a$}\}
\]
It is clear that $\mathcal{S}_{a}$ is an open subset of the space of closed $2$--forms homologous to $a$, so $\mathcal{S}_{a}$ is locally path connected. 
By Moser's argument, every smooth path between two points in $\mathcal{S}_a$ can be realized by the push-forward of a fixed symplectic form under a 1-parameter family of diffeomorphisms on $X$, so such a path is called an \textit{isotopy} of symplectic forms. 
We study the following well-known open question  (see \cite{SalamonSurvey}, \cite[Problem 2(a)]{MSIntro}). 

\begin{Question}\label{question isotopy} Does there exist a closed 4-manifold $X$ such that $\mathcal{S}_{a}$ is disconnected? 
\end{Question}

 Question \ref{question isotopy} is equivalent to asking whether cohomologous symplectic forms are unique up to isotopy.  It is known to McDuff \cite{McEx} that $\cS_a$ can be disconnected for some six-dimensional manifolds.  For the simplest cases of symplectic 4-manifolds such as complex projective planes and quadric surfaces, the question boils down to the connectedness of diffeomorphism groups of the underlying smooth manifolds \cite{SalamonSurvey}, due to the series of works of Gromov \cite{Gromov85}, Taubes \cite{TaubesGr,TaubesGrSW,TaubesSW}, McDuff \cite{McDuffRuled}, Abreu-McDuff \cite{Abreu}. For non-compact manifolds, Wang \cite{wang2024note} showed that there are infinitely many pairwise distinct non-standard symplectic forms that are standard at infinity on $S^1\times D^3$.  However, the closed case is still open. 

The following theorem, proved in Section \ref{sec_symplectic}, gives an affirmative answer to Question \ref{question isotopy}. Recall that a closed $4$--manifold is called a \emph{ruled surface} if it can be written as an $S^2$ bundle over a closed oriented $2$--manifold. 
\begin{Theorem} 
\label{thm_symplectic_intro}
Let $S^2\to X\to\Sigma$ be a ruled surface with $g(\Sigma)>0$. For every class $a\in H^2(X;\mathbb{R})$ with $a^2([X])>0$, the space $\mathcal{S}_{a}$
 has infinitely many components. 
\end{Theorem}

Theorem \ref{thm_symplectic_intro} indeed provides an answer to an even stronger non-uniqueness problem.  Recall that two symplectic forms are called \emph{formally homotopic} if there exists a path of nondegenerate forms that connects them; and \emph{homotopic} if such a path consists of symplectic forms. In other words, if we denote $\mathcal{S}=\cup_{a}\mathcal{S}_{a},$ 
which consists of \emph{all} symplectic forms on $X$, then two symplectic forms are homotopic if they are in the same connected component of $\mathcal{S}$. An open question in this direction asks whether there are symplectic forms that are formally homotopic but not homotopic in dimension $4$, see \cite[Page 252, Question (b)(c)]{CieliebakHprinciple}.  Again, this is known in higher dimensions: examples were found by Ruan \cite{Ruan94} using the Gromov--Witten invariants.  In contrast, by the work of Taubes, the Gromov--Witten invariants of symplectic 4-manifolds only depend on the smooth structure \cite{TaubesGr,TaubesGrSW,TaubesSW}, which adds to the subtlety in dimension 4. 

The following corollary gives an affirmative answer to this question.

\begin{Corollary}
\label{cor_intro_homotopic}
    For ruled surfaces $X$ as in Theorem \ref{thm_symplectic_intro}, there is an infinite family of symplectic forms $\{\omega_n\}_{n\in \mathbb{Z}^+}$ which are pairwise formally homotopic, cohomologous, but not homotopic.
\end{Corollary}

\begin{proof}
    Recall that symplectic forms on a ruled surface $X$ have been completely classified up to diffeomorphisms by Lalonde--McDuff \cite{LalondeMcDuff}. In particular, $\mathcal{S}_{a}$ is nonempty if and only if $a^2> 0$. And they further proved that every pair of cohomologous symplectic forms on $X$ is diffeomorphic.

By \cite{LalondeMcDuff,LalondeMcDuffJ}, two cohomologous symplectic forms on $X$ are isotopic if and only if they are homotopic. Therefore, Theorem \ref{thm_symplectic_intro} implies that there exists an infinite family of cohomologous symplectic forms $\{\omega_n\}_{n\in\mathbb{Z}^+}$, each belongs to a different connected component of $\mathcal{S}$.
 On the other hand, it is known that $\mathcal{S}$ only contains finitely many formal homotopic classes (see \cite[Section 2.4]{SalamonSurvey}).  This concludes the assertion.
\end{proof}

\subsection{The Dax invariant for closed surfaces} The proofs of the main theorems are based on a construction of the \emph{Dax invariant} for closed embedded surfaces. 
In 1972, Dax \cite{dax1972etude} showed that the relative homotopy group $\pi_{k}(\operatorname{Map}(N,M),\Emb(N,M))$ is isomorphic to a certain bordism group when $2\leq k\leq 2m-3n-3$.  This implies a description of $\pi_{1}(\Emb_{\partial}(I,X))$, where $\Emb_{\partial}(I,X)$ denotes the space of embeddings from the unit interval $I$ to a 4-manifold $X$ that are fixed near the boundary. Based on this work, Gabai \cite{gabai2021self} defined the relative Dax invariant for properly embedded disks $D_1, D_{2}\hookrightarrow M$ that are homotopic relative to the boundary. The relative Dax invariant is an obstruction to the existence of an isotopy between the two disks. The idea of the construction is to decompose the embedded disks $D_{1}, D_2$ into 1-parameter families of embedded arcs, which correspond to two paths $\gamma_1,\gamma_2$ in $\Emb_{\partial}(I,M)$ with the same endpoints, and then apply the Dax isomorphism to the element $[\gamma_{1}*\overline{\gamma}_{2}]\in \pi_{1}(\Emb_{\partial}(I,M))$ to obtain the Dax invariant $\Dax(D_1,D_2)$. 

Using the relative Dax invariant, Gabai \cite{gabai2021self} constructed a pair of embedded disks in $(S^2\times D^2)\natural (S^1\times D^3)$ with a common geometrically dual sphere that are homotopic but not isotopic relative to the boundary. Properties of the Dax invariant for embedded disks in 4-manifolds were also studied in \cite{KT24, schwartz20214}. The higher-dimensional cases were studied in \cite{kosanovic2024space}.

We extend the construction of the Dax invariant to \emph{closed} embedded surfaces in $\Sigma\times S^2$. Fix base points $b_0\in \Sigma$, $b_1\in S^2$, and let $b=(b_0,b_1)\in M$ be the base point of $M$. 
For notational convenience, let 
\[
\mathcal{C}=(\pi_{1}(M,b)\setminus \{1\})/\text{conjugation}.
\]
Let $\mathscr{I}_0: \Sigma\hookrightarrow \Sigma\times S^2$ be the embedding defined by $\mathscr{I}_0(x)=(x,b_1)$. 
More generally, for $d\in \mathbb{Z}$, let $\mathscr{I}_d$ be an embedding of $\Sigma$ in $\Sigma\times S^2$ whose image is given by the union of $\Sigma\times \{b_1\}$ and $|d|$ parallel copies of $\{b_0\}\times S^2$ after resolving singularities, where the orientations of the $S^2$'s depend on the sign of $d$ (see Section \ref{sec_isotopy_one_handles}). 
We define a relative Dax invariant for pairs of embeddings of $\Sigma$ in $M$ that are homotopic to $\mathscr{I}_d$. More precisely, suppose $i_1,i_2:\Sigma \hookrightarrow \Sigma\times S^2$ are two such embeddings, we will define an invariant $\Dax(i_1,i_2)\in \bZ[\mathcal{C}]^\sigma$. 
Here, $\bZ[\mathcal{C}]$ denotes the free abelian group generated by $\mathcal{C}$, and $\bZ[\mathcal{C}]^\sigma$ denotes its subgroup generated by elements of the form $[g]+[g^{-1}]$. 

\begin{remark}
\label{rmk_involution_fixed_point_C}
Since every pair of non-commutative elements in $\pi_1(M,b)$ generate a free group, the equation $xyx^{-1} = y^{-1}$ has no non-trivial solution in $\pi_1(M,b)$, so $\bZ[\mathcal{C}]^\sigma$ is also the subgroup of $\bZ[\mathcal{C}]$ consisting of fixed elements under the involution $\sigma:\bZ[\mathcal{C}]\to \bZ[\mathcal{C}]$ that takes $[g]$ to $[g^{-1}]$. 
\end{remark}

The relative Dax invariant satisfies the following properties.

\begin{enumerate}
\item (Isotopy invariance) $\Dax(i_1, i_2)=0$ if $i_1$ is smoothly isotopic to $i_2$.
\item (Additivity) $\Dax(i_1, i_2)+\Dax(i_2, i_3)=\Dax(i_1, i_3)$.
\item If $f$ is a diffeomorphism on $M$ that is homotopic to the identity, then $\Dax(i_1,i_2)=\Dax(f\circ i_1,f\circ i_2)$.
\item For every $d$ and every $[g]\in\mathcal{C}$, there exists a diffeomorphism $f:M\to M$ that is pseudo-isotopic to the identity, such that \[\Dax(\mathscr{I}_d, f\circ \mathscr{I}_d)=[g]+[g^{-1}],\] and $\{b_0\}\times S^2$ is a common geometric dual of $\mathscr{I}_d$ and $f\circ \mathscr{I}_d$. 
\end{enumerate}  

Properties (1), (2) above will follow immediately from the definition, 
Property (3) is proved by Lemma \ref{lem_Dax_change_by_MCG0_*}, and Property (4) follows from Lemma \ref{lem_self_referential_barbell_Dax}. 
Theorem \ref{thm_homot_non_isotopic_intro} then follows immediately from the isotopy invariance and Property (4) above. See also Lemma \ref{lem_cS_cS'_bijection} below for the discussion on reparametrizations.

In Section \ref{sec_classification}, we will show that the Dax invariant and the homology class of the embedding completely classify all isotopy classes of embeddings of $\Sigma$ in $\Sigma\times S^2$ that are geometrically dual to $\{b_0\}\times S^2$. 

\begin{Theorem}\label{intro thm: classification}
For $d\in \mathbb{Z}$, consider the set of embeddings  $i:\Sigma\hookrightarrow  M$ whose projection to $S^2$ has degree $d$, such that $i$ intersects $\{b_0\}\times S^2$ transversely and positively at one point. Let $\mathcal{F}_d$ be the set of equivalence classes of such embeddings up to reparametrizations and isotopies in $M$. Then each element in $\mathcal{F}_d$
can be represented by an embedding homotopic to $\mathscr{I}_d$, and
the map 
\begin{align*}
\mathcal{F}_{d} &\to \mathbb{Z}\\
[i]&\mapsto \Dax(\mathscr{I}_{d},i)
\end{align*}
is a bijection, where $i:\Sigma\hookrightarrow M$ denotes a representative homotopic to $\mathscr{I}_{d}$.  
\end{Theorem}

To define the relative Dax invariant for closed surfaces, we fix a handle decomposition of $\Sigma=H_0\cup H_1\cup H_2$, where $H_{j}$ denotes the union of $j$-handles. The embeddings $\mathscr{I}_d$ are chosen so that $\mathscr{I}_d|_{H_0\cup H_1} = \mathscr{I}_0|_{H_0\cup H_1}$ for all $d$.  Given  $i_1, i_2$, first apply isotopies $\mathcal{I}_{1}, \mathcal{I}_2$ on $i_1, i_2$ so that 
\[
i_{1}|_{H_0\cup H_1} = i_{2}|_{H_0\cup H_1} = \mathscr{I}_d|_{H_0\cup H_1}.
\]
Let $M_2$ be the complement of a tubular neighborhood of $\mathscr{I}_d(H_0\cup H_1)$ in $M$. By restricting $i_{1}, i_{2}$ to the preimages of $M_2$, we obtain embeddings $i'_{1}, i'_{2}:(D^2,\partial D^2) \hookrightarrow (M_2,\partial M_2)$. We prove that one can always choose $\mathcal{I}_{1}, \mathcal{I}_2$ such that $i_1'$, $i_2'$ are homotopic to each other relative to the boundary. So a relative Dax invariant can be defined for $i_1'$ and $i_2'$, which takes value in $\mathbb{Z}[\pi_1(M_2)\setminus\{1\}]$. Then we define $\Dax(i_1, i_2)$ to be the quotient image of $\Dax(i'_1, i'_2)$ in $\mathbb{Z}[\mathcal{C}]$. The main step is to show that $\Dax(i_1, i_2)$ is independent of the choice of $\mathcal{I}_{1}, \mathcal{I}_2$. To prove this, we note that choices of $\mathcal{I}_{1}, \mathcal{I}_{2}$ differ from each other by loops of embeddings from $H_0\cup H_1$ to $M$. By applying the isotopy extension theorem to such loops, we get compositions of barbell diffeomorphisms on $M_2$. A careful analysis of the effect of these diffeomorphisms shows that they do not change the quotient image of $\Dax(i'_1, i'_2)$.
 
Our construction is the first generalization of the Dax invariant to closed embedded surfaces. We remark that, in principle, the same idea can be applied to the embedding of closed surfaces in arbitrary 4-manifolds. In this case, the range of the relative Dax invariant is given by the quotient of $\bZ[\pi_{1}(M)\setminus\{1\}]$ by the subgroup generated by the image of the Dax homomorphism $d_3$ (see \eqref{eq: Dax homomorphism} below) and the ambiguity from changing the initial isotopy on $H_0\cup H_1$. The difficulty in establishing the general theory is in the computation of the contributions from changing the initial isotopy.

The paper is organized as follows. In Section \ref{sec_prelim}, we review some backgrounds about the barbell diffeomorphism and the relative Dax invariant. In Section \ref{sec_isotopy_one_handles}, we show for any two homotopic embeddings of $\Sigma$ in $\Sigma\times S^2$, one can always isotope them so that they are equal on the 0- and 1-handles and their restrictions to 2-handle are homotopic in the complement of the image of the 0- and 1-handles. This is used in Section \ref{sec_dax_invariant} to prove that the Dax invariant is well-defined and satisfies various properties mentioned above. In Section \ref{sec_Mcg_Sigma_S2}, we use the relative Dax invariant to prove Theorem \ref{thm: MCG infinitely generated}. In Section \ref{sec_classification}, we show that if two embeddings of $\Sigma$ in $\Sigma\times S^2$ are homotopic to each other, are both dual to $\{b_0\}\times S^2$, and have relative Dax invariant zero, then they are isotopic to each other. This implies the classification result for isotopy classes of embeddings of $\Sigma$ in $\Sigma\times S^2$ that are dual to $\{b_0\}\times S^2$ as stated in Theorem \ref{intro thm: classification}. Section \ref{sec_symplectic} proves Theorem \ref{thm_symplectic_intro}.

\subsection*{Acknowledgments} 
The authors  thank Jae Choon Cha for inspiring discussions on the proof of Theorem \ref{intro thm: classification}. J. Lin is partially supported by NSFC grant 12271281. 
W. Wu is partially supported by
NSFC 12471063.
Y. Xie is partially supported by National Key R\&D Program of China 2020YFA0712800 and NSFC 12341105.
B. Zhang is partially supported by NSF grant DMS-2405271 and a travel grant from the Simons Foundation.  

\section{Preliminaries}
\label{sec_prelim}

In this section, we review some key results about the barbell diffeomorphisms introduced by Budney-Gabai \cite{Gabai20}. This diffeomorphism naturally arises when applying the isotopy extension theorem to the so-called ``spinning family'' of embedded arcs.

\subsection{The spinning family of embedded arcs}
We start by introducing some notations. If $A,B$ are subsets of a 4-manifold $X$, we use $\pi_{1}(X;A,B)$ to denote the set of paths in $X$ from some point in $A$ to some point in $B$, up to homotopies keeping the starting point in $A$ and the ending point in $B$. Note that by definition, $\pi_1(X;A,B)$ may not be a group. When $A=B$, we also write $\pi_{1}(X;A,A)$ as $\pi_1(X;A)$. When $A$ is simply-connected,  $\pi_{1}(X;A)$ can be canonically identified with $\pi_{1}(X;a)$ for every $a\in A$.

We also consider  $\hat{\pi}_{1}(X;A,B)$, the set of connected components of the following space 
\[
\{\text{continuous map } \gamma: [0,1]\to X\mid \gamma(0)\in A, \gamma(1)\in B, \gamma(0,1)\cap (A\cup B)=\emptyset\}.
\]
We have a natural map 
\begin{equation}\label{eq: pi1hat to pi1}
\hat{\pi}_{1}(X;A,B)\to  \pi_{1}(X;A,B).   
\end{equation}
If $A,B$ are both embedded arcs, the map (\ref{eq: pi1hat to pi1}) is an isomorphism because of standard transversality results. On the other hand, if $A$ is an embedded arc and $B$ is an embedded surface, the map (\ref{eq: pi1hat to pi1}) is surjective but may not be injective. 

We say that a smooth embedding $k:X\to Y$ between manifolds with boundaries is \emph{neat}, if $k^{-1}(\partial Y) = \partial X$ and $k$ intersects $\partial Y$ transversely. 

\begin{Definition}\label{defn_Emb_I}
Let $X$ be a smooth manifold with nonempty boundary.  Let  $I_0\hookrightarrow X$ be a neatly embedded arc. We use $\Emb_{\partial}(I,X)$ to denote the space of embedded arcs that equal $I_0$ near the boundary. 
Let $s_0$ be a framing on $I_0$ (i.e., a trivialization of its normal bundle). We use $\Emb^{\Fr}_{\partial}(I,X)$ to denote the space of framed embedded arcs that equal $(I_0,s_0)$ near the boundary points. 
\end{Definition}

Now we introduce the model \emph{spinning family}. Consider the sphere 
\[
S_0=\{(x,y,z,0)\in \mathbb{R}^4\mid x^2+y^2+z^2=4\},
\]
and its closed neighborhood 
$X_0=\{\vec{v}\in \mathbb{R}^4\mid \operatorname{distance}(\vec{v},S_0)\leq 1\}\subset \mathbb{R}^4.$
Pick an even smooth function $\iota: [-3, 3]\to [0,2.5]$ that equals $0$ near $\pm 3$, equals $2.5$ near $0$, and that its graph is contained in the set $\{(x,y)\mid 1\le x^2+y^2\le 9\}$. For $\theta\in S^1=\mathbb{R}/2\pi \mathbb{Z}$, we consider the embedded arc 
\[
J_{\theta}=\{(x,\cos \theta \cdot \iota(x), \sin\theta \cdot \iota(x),0)\mid x\in [-3,3]\}\hookrightarrow X_0.
\]
Then $\{J_{\theta}\}_{\theta\in S^1}$ is a loop in $\Emb_{\partial}(I,X_0)$ based at $J_0$, called the ``model spinning family''. The projection of the normal bundle of $J_\theta$ to the plane $x=0$ gives a framing on $J_{\theta}$, denoted by $s'_{\theta}$. Hence $\{(J_{\theta},s'_{\theta})\}$ is a loop in $\Emb^{\Fr}_{\partial}(I,X_0)$. We use $\gamma_0$ to denote the path $\{(0,y,0,0)\mid y\in [2,2.5]\}$ from $J_0$ to $S_0$.

Now we consider the general case. We fix a Riemannian metric on a smooth $4$--manifold $X$.

\begin{figure}
	\begin{overpic}[width=0.6\textwidth]{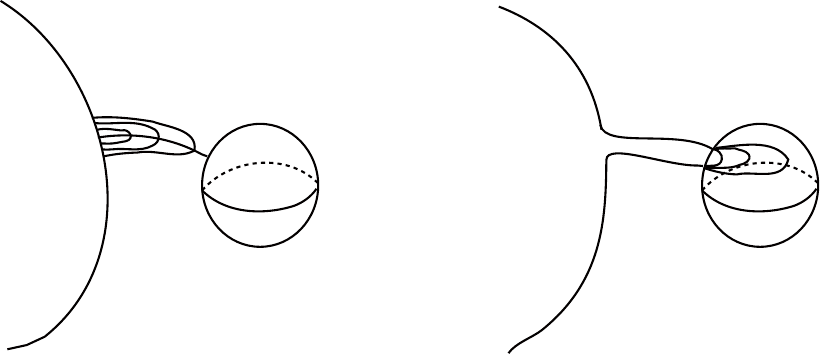}
		\put(10,18){$I_0$}
		\put(60,15){$S_0$}
            \put(40,65){$\gamma$}
	\end{overpic}
	\caption{Definitions \ref{defi: general spining family}}
    \label{fig:def_rho_eta}
\end{figure}

\begin{Definition}\label{defi: general spining family} Let $I_0\hookrightarrow X$ be a neatly embedded arc and let $S\hookrightarrow X$ be an embedded sphere that is disjoint from $I_0$ and has self-intersection $0$. Given any element $g\in \hat{\pi}_{1}(X;I_0,S)$, we represent it by a smoothly embedded arc $\gamma: I\hookrightarrow X$ such that $\gamma$ is perpendicular to $I_0$ at $\gamma(0)\in \operatorname{int}(I_0)$ and perpendicular to $S$ at $\gamma(1)$. 
Pick a framing on $S$ and a framing on $\gamma$. Using the exponential maps, we can find a closed neighborhood $U$ of $S\cup \gamma$ and a diffeomorphism $\varphi:U\xrightarrow{\cong }X_0 $ such that $\varphi(S)=S_0, \varphi(\gamma)=\gamma_0$ and $\varphi(I_0\cap U)=J_0$. For $\theta\in S^1$, we let $I_\theta\hookrightarrow X$ be the embedded arc that equals $I_0$ outside $U$ and equals $\varphi^{-1}(J_\theta)$ inside $U$. Then the loop $\{I_{\theta}\}_{\theta\in S^1}$ in $\Emb_\partial (I,X)$ is called the \emph{spinning family} of $I_0$ along $(\gamma,S)$. This loop represents an element $\rho_{\gamma}\in \pi_{1}(\Emb_{\partial}(I,X),I_0)$. This element $\rho_{\gamma}$ only depends on $g=[\gamma]$ so we also denote it by $\rho_g$. See Figure \ref{fig:def_rho_eta} for an illustration.

Given a framing  $s_0$ on $I_0$, we homotope $s_0$ such that $s_0|_{I_{0}\cap U}$ equals the pull back of $s'_{0}$. Then we define the framing $s_{\theta}$ on $I_{\theta}$ by pulling back $s'_{\theta}$. This gives a lifted loop $\{(I_{\theta},s_{\theta})\}_{\theta\in S^1}$ in $\Emb^{\Fr}_{\partial}(I,X)$. Since $\pi_{2}(SO(3))=0$, such a lift is unique up to homotopy relative to $(I_0,s_0)$. We use $\widetilde{\rho}_{g}\in \pi_{1}(\Emb^{\Fr}_{\partial}(I,X))$ to denote the element represented by this loop. 
\end{Definition}

The following definitions are special cases of Definition \ref{defi: general spining family}. 

\begin{Definition}\label{defi: spinning I0}
Let $I_0$ be a framed embedded arc in $X$. Given 
\[
g\in \hat \pi_{1}(X;I_0,I_0) \cong \pi_{1}(X;I_0,I_0),
\] we represent it by an embedded arc $\gamma: I\to X$ such that $\gamma(0,1)\cap I_0=\emptyset$ and that
$\gamma(0)<\gamma(1)$ under the parametrization of $I_0$. Let $m$ be a meridian sphere of $I_0$ around $\gamma(1)$. We perturb $\gamma$ to an embedded arc $\gamma'$ from $I_0$ to $m$. 
Then we denote the spinning family $\widetilde{\rho}_{\gamma'}$ by 
\[
\widetilde{\tau}_{g}\in \pi_{1}(\Emb^{\Fr}_{\partial}(I,X);I_0).\] We use $\tau_{g}$ to denote the image of $\widetilde{\tau}_{g}$ in $\pi_{1}(\Emb_{\partial}(I,X);I_0)$. 
\end{Definition}

\begin{Definition}\label{defi: spinning I0 to I1}
Let $I_0, I_1$ be embedded arcs in $X$ and let 
\[
g\in \hat \pi_{1}(X;I_0,I_1)\cong \pi_{1}(X;I_0,I_1).
\]
 Let $X_1=X\setminus \nu(I_1)$. Let $m_{1}:S^2\hookrightarrow X_1$ be a meridian sphere of $I_1$. We represent $g$ by an embedded path $\gamma$ from $I_0$ to $I_1$ and slightly perturb it to a path $\gamma'$ in $X_1$ from $I_0$ to $m_1$. Then we denote the spinning family $\widetilde{\rho}_{\gamma'}$ by 
\[
\widetilde{\tau}_{g}\in \pi_{1}(\Emb^{\Fr}_{\partial}(I,X_1);I_0).\] 
We use $\tau_{g}$ to denote the image of $\widetilde{\tau}_{g}$ in $\pi_{1}(\Emb_{\partial}(I,X_1);I_0)$. 
\end{Definition}

\subsection{The barbell diffeomorphism} Given any $I_0\hookrightarrow X$, there is a well-defined map 
\begin{equation}\label{eq: arc pushing}
\partial: \pi_1(\Emb^{\Fr}_{\partial}(I,X);I_0)\to \pi_{0}(\Diff_\partial(X\setminus\nu(I_0))).    
\end{equation}
obtained by the isotopy extension theorem. We call (\ref{eq: arc pushing}) \emph{the arc pushing map}. Here, $\Diff_\partial$ denotes the group of diffeomorphisms that equal the identity near the boundary. 

Now we review the definition of the barbell diffeomorphisms. Let $D^4$ be the unit disk in $\mathbb{R}^4$. Consider the embedded arcs $I_0\in \{(-\frac{1}{2},t,0,0)\in D^4\}$ and $I_1\in \{(\frac{1}{2},t,0,0)\in D^4\}$, equipped with the standard framing. Take the model barbell 
\[\mathcal{B}:=D^4\setminus \nu(I_0\cup I_1)\cong (S^2\times D^2)\natural (S^2\times D^2).\] 
Let $g$ be the unique element in $\pi_{1}(D^4;I_0,I_1)$.  Consider the spinning family 
$\widetilde{\tau}_{g}\in \pi_{1}(\Emb^{\Fr}_{\partial}(I,D^4\setminus\nu(I_1)),I_0)$ given in Definition \ref{defi: spinning I0 to I1} and its image $\partial(\widetilde{\tau}_{g})\in \pi_0(\Diff_\partial (\mathcal{B}))$ under the arc pushing map (\ref{eq: arc pushing}). We fix a representative 
$\beta: \mathcal{B}\to \mathcal{B}$ of $\partial (\widetilde{\tau}_{g})$ and call it \emph{the model barbell diffeomorphism}.

Next, we define the implemented barbell diffeomorphism for any ordered ordered pair of disjoint, embedded spheres with self-intersection $0$ together with an arc connecting them. Let $S_0,S_1\hookrightarrow X$ be disjoint embedded spheres with self-intersection number $0$. Consider an element $g\in \hat{\pi}_{1}(X;S_0,S_1)$, represented by an embedded arc $\gamma:I\to X$. 

Let $\bar{\nu}(S_0\cup \gamma\cup S_1)$ be a regular neighborhood of $S_0\cup \gamma\cup S_1$. We define a diffeomorphism supported on $\overline{\nu}(S_0\cup \gamma\cup S_1)$ as follows: Fix trivializations of the normal bundles $NS_0$ and $NS_1$, which are unique up to isotopy. This gives diffeomorphisms  
\[
\varphi_0: \nu(S_0)\xrightarrow{\cong} S^2\times D^2,\quad \varphi_1: \nu(S_1)\xrightarrow{\cong }S^2\times D^2.
\]
After isotopy of $\gamma$, we may assume 
\[
\gamma(t)=\begin{cases} \varphi^{-1}_{0}((1,0,0)\times (t,0))&\text{ for } t\in [0,\epsilon],\\
  \varphi^{-1}_{1}((-1,0,0)\times (t-1,0))&\text{ for }t\in [1-\epsilon,1].  
\end{cases}
\]
Then $\varphi_{0},\varphi_1$ induces a trivialization of the normal bundle $N(\gamma|_{[0,\epsilon]\cup [1-\epsilon,1]})$. Since $SO(3)$ is connected, we may extend this trivialization to a trivialization $\varphi_{2}:N\gamma\cong \gamma\times \mathbb{R}^{3}$. Putting all these trivializations together, we obtain a diffeomorphism 
\[
\varphi: \bar{\nu}(S_0\cup \gamma\cup S_1)\xrightarrow{\cong }\mathcal{B}.
\]
Then we extend the diffeomorphism 
\[\varphi^{-1}\circ \beta \circ \varphi: \bar{\nu}(S_0\cup \gamma\cup S_1)\to \bar{\nu}(S_0\cup \gamma\cup S_1)\] by the identity map on $X\setminus \bar{\nu}(S_0\cup \gamma\cup S_1)$. This gives the \emph{implemented barbell diffeomorphism} \[\beta_{(S_0,S_1,\gamma)}: X\to X.\]  
The isotopy class of the implemented barbell diffeomorphism only depends on the isotopy class of $S_0\sqcup S_1$ and the homotopy class $g=[\gamma]$. So we also write  $\beta_{(S_0,S_1,\gamma)}$ as $\beta_{(S_0,S_1,g)}$.

\begin{remark}
Up to homotopy, there are two possible choices of $\varphi_{2}$ because $\pi_{1}(SO(3))=\mathbb{Z}/2$. This choice won't affect the isotopy class of $\beta_{(S_0,S_1,\gamma)}$. The resulting diffeomorphisms for different choices of $\varphi_{2}$ are conjugate by a diffeomorphism supported near $\partial \bar{\nu}(S_0\cup\gamma\cup S_1)$. Since $\beta_{(S_0,S_1,\gamma)}$ is supported away from $\partial \bar{\nu}(S_0\cup\gamma\cup S_1)$, its isotopy class is unchanged. 
\end{remark}

The following result follows immediately from the above description of the barbell diffeomorphisms, and it is also a direct consequence of \cite[Construction 5.3]{BG2019}. 

\begin{Proposition}\label{prop: barbell diff as arc pushing map} Let $I_0$ be a framed embedded arc and let $S$ be an embedded sphere with self-intersection $0$ and disjoint from $I_0$. Let $m_0$ be a meridian of $I_0$. Given each $g\in \hat{\pi}_{1}(X;I_0,S)$, one obtains a corresponding element $g'\in \hat{\pi}_{1}(X;m_0,S)$ by a standard local perturbation. Let $\widetilde{\rho}_{g}$ be the spinning family given in Definition \ref{defi: general spining family} and let $\partial$ be the arc pushing map given in (\ref{eq: arc pushing}). Then we have 
\[
\partial(\widetilde{\rho}_{g})=\beta_{(m_0,S,g')}\in \pi_0(\Diff(X\setminus \nu(I_0))).
\makeatletter\displaymath@qed\makeatother
\]
\end{Proposition}

\begin{Corollary}
\label{cor_parallel_barbell_compose_to_id}
Let $X=S^2\times D^2$. Take points $b_0\neq b_1\in D^2$ and $x\in S^2$. Let $S_0=S^2\times \{b_0\}$ and $S_1=S^2\times \{b_1\}$. Let $\gamma:I\to \{x\}\times D^2 $ be an embedded arc from $(x,b_0)$ to $(x,b_1)$. Then $\beta_{(S_0,S_1,\gamma)}$ is smoothly isotopic to the identity relative to $\partial X$.  
\end{Corollary}
\begin{proof}
We attach a 1-handle along the knot $\{x\}\times \partial D^2\subset \partial X$ and obtain a manifold $X'\cong D^4$ with an embedded arc $I_0$ such that $X=X'\setminus \nu(I_0)$. Then $S_0\cup S_1\hookrightarrow X$ is isotopic to a union of two parallel meridians of $I_0$.   By Proposition \ref{prop: barbell diff as arc pushing map}, we have 
\[
\beta_{(S_0,S_1,\gamma)}=\partial(\widetilde{\rho}_{\gamma'})\in \pi_{0}(\Diff(S^2\times D^2)),
\]
where $\gamma'$ is the standard arc from $I_0$ to its meridian $S_1$. It is straightforward to check that the spinning family 
$\widetilde{\rho}_{\gamma'}\in \pi_{1}(\Emb^{\Fr}_{\partial}(I,X');I_0)$
is trivial, so the desired result follows. 
\end{proof}

Given two implemented barbell diffeomorphisms, if one can isotope them such that their supports are disjoint, then they commute up to isotopy.

The following connected sum formula will be useful later.

\begin{Proposition}\label{prop: connected sum}
Let $S,S_1, S_2$ be disjoint embedded spheres in $X$ with self-intersection number $0$. 
Take an embedded arc $\gamma_{01}: I\to X$ from $S$ to $S_1$, and an embedded arc $\gamma_{12}: I\to X$ from $S_1$ to $S_2$. We assume that $\gamma_{01}$ and $\gamma_{12}$ are disjoint from each other and that
\[
(\gamma_{01}(0,1)\cup \gamma_{12}(0,1))\cap (S\cup S_1\cup S_2)=\emptyset.
\]
Take a path $\eta:I\to S_1$ from $\gamma_{01}(1)$ to $\gamma_{12}(0)$. Let $\gamma_{02}:I\to X$ be an embedded arc from $S$ to $S_2$ obtained by a smoothing of $\gamma_{01}\ast \eta\ast \gamma_{12}$. We take a tube along $\gamma_{12}$ and let $S_{12}$ be the connected sum of $S_1$ and $S_2$ along this tube. Then we have 
\begin{equation}\label{eq: barbell connected sum}
\beta_{(S,S_2,\gamma_{02})}\circ \beta_{(S,S_1,\gamma_{01})}=\beta_{(S,S_{12},\gamma_{01})}\in \pi_{0}(\Diff_\partial(X)).
\end{equation}

\end{Proposition}
\begin{proof} Without loss of generality we may assume that \[X=\bar{\nu}(S\cup S_1\cup S_2\cup \gamma_{01}\cup \gamma_{02})\cong (S^2\times D^2)\natural (S^2\times D^2)\natural (S^2\times D^2).\]
Then \[\partial X\cong (S^1\times S^2)\sharp (S^1\times S^2)\sharp (S^1\times S^2).\]
By attaching a 2-handle along the curve $S^1\times \{\pt\}$ in the first summand, we obtain a manifold $X'$ and a framed embedded arc $I_0\hookrightarrow X'$ such that $X=X'\setminus \nu(I_0)$.  
Extend $\gamma_{02},\gamma_{01}$ to $I_0$. By Definition \ref{defi: general spining family}, it is straightforward to check the relation 
$
\widetilde{\rho}_{(I_0,S_{2},\gamma_{02})}\cdot \widetilde{\rho}_{(I_0,S_{1},\gamma_{01})}=\widetilde{\rho}_{(I_0,S_{12},\gamma_{01})}.
$
Note that $S\hookrightarrow X'$ is a meridian of $I_0$. So (\ref{eq: barbell connected sum}) follows from Proposition \ref{prop: barbell diff as arc pushing map}. 
\end{proof}

\subsection{The relative Dax invariant for embedded disks} In this subsection, we briefly review the definition of the Dax invariant for embedded disks. As before, we let $I_0$ be a neatly embedded arc in a 4-manifold $X$. Every loop in $\Emb_{\partial}(I,X)$ with base point $I_0$ defines a map $I\times I\to X$, where $\partial(I\times I)$ is mapped to $I_0$. This gives rise to a homomorphism 
\begin{equation}
\label{eqn_pi1Emb_to_pi2}
\mathcal{F}: \pi_{1}(\Emb_{\partial}(I,X);I_0)\to \pi_2(X;I_0).
\end{equation}
Here, $\pi_2(X;I_0)$ denotes the second homotopy group of $X$ with base point in $I_0$. Since $I_0$ is contractible, the homotopy group is well-defined (see also \cite[Section 3.3]{lin2025mapping}). 
 Following \cite{gabai2021self}, we denote the kernel of $\mathcal{F}$ by 
\[\pi^{D}_{1}(\Emb_{\partial}(I,X);I_0).\] 
Here, $D$ stands for ``Dax''. 

Let $\pi=\pi_{1}(X;I_0)$. By the Dax isomorphism theorem (see \cite[Theorem 0.3]{gabai2021self}), we have an isomorphism
\begin{equation}\label{eq: Dax homomorphism}
  \Dax: \pi^D_{1}(\Emb_{\partial}(I,X);I_0)\to \mathbb{Z}[\pi\setminus \{1\}]/\Imm(d_{3}),
\end{equation}
where $d_{3}$ is a homomorphism of the form
$
d_3: \pi_{3}(X;I_0)\to \mathbb{Z}[\pi\setminus \{1\}].
$
We refer to \cite{gabai2021self} for a description of the map $d_3$.

For later reference, we briefly recall the definition of the Dax isomorphism (\ref{eq: Dax homomorphism}) as follows. 
Given any $[\gamma]\in \pi^D_{1}(\Emb_{\partial}(I,X);I_0)$, one can extend $\gamma$ to a map $h:D^2\to \operatorname{Imm}_{\partial}(I,X)$. By a general position argument, we can choose $h$ such that $h(x)\in \Emb_{\partial}(I,X)$ for all 
$
x\in D^2\setminus \{x_{1},x_{2},\cdots x_{m}\}.
$
Moreover, we may assume that each $h(x_{j})$ is an immersion with a single double point, and that the image of $h$ on a small circle around $x_{j}$ is a resolution of the double point. 
Suppose $h(x_j): I\to X$ is an immersion with a double point $h(x_j)(t_0)=h(x_j)(t_{1})$ for $t_0<t_{1}$. One may concatenate the paths $h(x_j)|_{[0,t_1]}$ and the reverse of $h(x_j)|_{[0,t_0]}$ to obtain a loop in $X_2$, which represents an element in $\pi$. The image $\Dax([\gamma])$ is defined to be the formal sum of the corresponding elements in $\pi$, with signs given by local orientations at the double point, for all $x_j$.

We also introduce the following definitions. 
\begin{Definition}\label{defn_Emb_D_X}
Fix a neat embedding $i_0: D^2\to X$ and a collar neighborhood $N(\partial D^2)$ of $D^2$. We let  $\Emb_{\partial}(D^2,X)$ be the space of neatly embedded disks that equals $i_0$ on $N(\partial D^2)$.
\end{Definition}

\begin{Definition}
\label{def_scan_into_pi1Emb}
Given $i_1,i_2\in \Emb_{\partial}(D^2,X)$, we obtain two paths in the space $\Emb_{\partial}(I, X)$ by scanning $D^2$ with a 1-parameter family of arcs. Since $i_1=i_2$ near $\partial D^2$, we may choose the scanning of $D^2$ so that the two paths have the same endpoints and that their concatenation is a loop in $\Emb_{\partial}(I, X)$ with base point $I_0$, where $I_0$ is a fixed arc in the neighborhood $N(\partial D^2)$ in Definition \ref{defn_Emb_D_X}. 
This gives a well-defined element in $\pi_{1}(\Emb_{\partial}(I, X);I_0),$ which we denote by $\langle i_1, i_2\rangle$.
\end{Definition}

By Smale's theorem, the diffeomorphism group of $D^2$ relative to the boundary is contractible. So the value of $\langle i_1, i_2\rangle$ only depends on the images of $i_1$ and $i_2$. As a result, if $D_1$ and $D_2$ are neatly embedded disks in $X$ such that there exist $i_1,i_2\in \Emb_{\partial}(D^2,X)$ with $D_1 = \Imm(i_1)$, $D_2 = \Imm(i_2)$, then we will use $\langle D_1,D_2\rangle$ to denote the value of $\langle i_1,i_2\rangle$. Similarly, since the orientation-preserving diffeomorphism group of $I$ is contractible, we will often identify $I_0$ with its image in $X$ when there is no risk of confusion.

Now assume $i_1, i_2$ are homotopic relative to their boundaries. Then $\langle i_1,i_2\rangle$ belongs to the Dax group $\pi^{D}_{1}(\Emb_{\partial}(I,X))$. Following \cite{gabai2021self}, we define the relative Dax invariant 
\[
\Dax(i_1,i_2):=\Dax(\langle i_1,i_2\rangle )\in \mathbb{Z}[ \pi\setminus 1]/\operatorname{Im}(d_3).
\]

\section{Isotopy of $0$- and $1$-handles}
\label{sec_isotopy_one_handles}

Let $\Sigma$ be a closed surface with genus $\ell\ge 1$, and let $M=\Sigma\times S^2$.
Pick a base point $b=(b_0,b_{1})\in \Sigma\times S^2 = M$. Let $ \mathscr{I}_0: \Sigma\to M$ be the embedding defined by
 $\mathscr{I}_0(x)=(x,b_1)$. Set $\Sigma_{(0)}=\mathscr{I}_0(\Sigma)$. 
 Given integers $d>0$ and $1\leq s\leq d$, let $p_s$ be distinct points in $\Sigma-\{b_0\}$.
 Let $\Sigma_{(d)}\subset M$ be the surface obtained
by resolving the intersection points in the union of $\Sigma_{(0)}$ and $\{p_1,\cdots p_{d}\}\times S^2$. Given 
$d<0$, $\Sigma_{(d)}$ is defined  in a similar manner, but using $\{p_1,\cdots ,p_{|d|}\}\times S^2$ with the opposite orientation. Let $\mathscr{I}_d:\Sigma\to \Sigma_{(d)}$ be a diffeomorphism obtained by modifying $\mathscr{I}_0$ near $p_s\in \Sigma$
($1\le s\le d$).
By definition,
 $\mathscr{I}_d$ is assumed to coincide with $\mathscr{I}_0$ near $b_0$.

Let $\mathcal{E}_d$ be the set of smooth embeddings $i:(\Sigma,b_0)\to (M,b)$ such that $i$ is homotopic to $\mathscr{I}_d$ relative to $b_0$. The next lemma shows that based isotopy is equivalent to free isotopy for embedded surfaces homotopic to $\Sigma_{(d)}$. 

\begin{Lemma}
\label{lem_cS_cS'_bijection}
    Let $\sim$ be the equivalence relation on $\mathcal{E}_d$ given by isotopy relative to $b_0$. 
    Let $\mathcal{E}'_d$ be the set of (non-parametrized) embedded surfaces in $M$ that are homotopic to $\Sigma_{(d)}$. Let $\sim'$ be the equivalence relation on $\mathcal{E}'_d$ such that $S_1\sim' S_2$ if and only if there exist diffeomorphisms $i_1:\Sigma\to S_1$ and $i_2:\Sigma\to S_2$ such that $i_1$ and $i_2$ are isotopic in $M$. Then the map from $\mathcal{E}_d/\sim$ to $\mathcal{E}'_d/\sim'$ , defined by taking each embedding to its image, is a bijection. 
\end{Lemma}

\begin{proof}
    Let $\Phi: \mathcal{E}_d\to \mathcal{E}'_d$ be the map such that $\Phi(i)$ is the image of $i$. It is clear that if $i\sim i'$, then $\Phi(i)\sim'\Phi(i')$, so the map $\Phi$ induces a map from $\mathcal{E}_d/\sim $ to $\mathcal{E}_d/\sim'$, which we denote by $\overline{\Phi}$. We show that $\overline{\Phi}$ is a bijection.

    We first show that $\overline{\Phi}$ is surjective. For each $S\in \mathcal{E}'_d$, let $i:\Sigma\to S$ be a diffeomorphism that is homotopic to $\mathscr{I}_d$ in $M$. After perturbation of $i$, we may assume that $i(b_0)\neq b$. Let $H:\Sigma\times[0,1]\to M$ be a homotopy such that $H(-,0) = \mathscr{I}_0$ and $H(-,1) = i$. Without loss of generality, we may assume that $H$ is smooth. Let $\gamma:[0,1]\to M$ be defined by $\gamma(t) = H(b_0,t)$. After a generic perturbation of $H$ relative to $\Sigma\times \{0,1\}$, we may assume that $\gamma$ is an embedded arc and the image of $[0,1)$ under $\gamma$ is disjoint from $S$. Let $S'$ be obtained from $S$ by attaching a tube along $\gamma$. Then $S' \sim' S$ and $S'$ is in the image of $\overline{\Phi}$. 

    Now we show that $\overline{\Phi}$ is injective. Suppose $i_1,i_2\in\mathcal{E}$ and $\Phi(i_1) \sim' \Phi(i_2)$. Then there exists a diffeomorphism $\varphi:\Sigma\to \Sigma$ such that $i_1$ and $i_2\circ \varphi$ are isotopic in $M$. Let $p:M=\Sigma\times S^2\to \Sigma $ be the projection map. Since $i_1,i_2$ are homotopic to $\mathscr{I}_d$, we know that $p\circ i_1$ is homotopic to $\id_\Sigma$, and $p\circ i_2\circ \varphi$ is homotopic to $\varphi$. As a result, $\varphi$ is homotopic to $\id_\Sigma$, so it is isotopic to $\id_\Sigma$. Hence $i_1$ and $i_2$ are isotopic in $M$.

    We need to show that $i_1$ and $i_2$ are isotopic \emph{relative to $b_0$} in $M$. Let $H:\Sigma\times [0,1]\to M$ be an isotopy such that $H(-,0)=i_1$ and $H(-,1)=i_2$. Let $\gamma:[0,1]\to M$ be defined by $\gamma(t) = H(b_0,t)$, then $\gamma$ is a loop based at $b$. We discuss two cases.

    Case 1: The genus $\ell$ of $\Sigma$ is greater than $1$. Since both $i_1$ and $i_2$ are homotopic to $\mathscr{I}_d$ relative to $b_0$, the loop $\gamma$ must represent an element in the center of $\pi_1(M)$. When $\ell>1$, the center of $\pi$ is trivial, and hence $\gamma$ must be contractible. So there exists a loop of diffeomorphisms $f_t:M\to M$ such that $f_t(b) = \gamma(t)$, $f_0=f_1=\id_M$. Let $H'(x,t) = f_t^{-1}\circ H(x,t)$, then $H'$ is an isotopy from $i_1$ to $i_2$ relative to $b_0$.

  Case 2: The genus $\ell$ of $\Sigma$ is $1$. Recall that $p: M=\Sigma\times S^2\to \Sigma $ denotes the projection map. Since $\Sigma$ is diffeomorphic to a Lie group, there exists a 1-parameter family of diffeomorphisms $f_t:M\to M$ such that $p\circ f_t(b) = p\circ \gamma(t)$, $f_0=f_1=\id_M$. Let $\hat H(x,t) = f_t^{-1}\circ H(x,t)$, let $\hat \gamma(t) = \hat H(b_0,t)$, then $\hat H$ is another isotopy from $i_1$ to $i_2$, and $\hat \gamma$ is contained in $\{b_0\}\times S^2$. Since $\pi_1(S^2)=1$, the loop $\hat \gamma$ is contractible, and hence there exists a 1-parameter family of diffeomorphisms $u_t:M\to M$ such that $u_t(b) = \hat \gamma(t)$, $u_0=u_1=\id_M$. Let $H'(x,t) = u_t^{-1}\circ \hat H(x,t)$, then $H'$ is an isotopy from $i_1$ to $i_2$ relative to $b_0$.
\end{proof}

\begin{remark}
\label{rmk_relative_Dax_surfaces}
Because of Lemma \ref{lem_cS_cS'_bijection}, for $S_1,S_2\in \mathcal{E}'_d$, we will sometimes abuse notation and use $\Dax (S_1,S_2)$ to denote the Dax invariant of their images in $\mathcal{E}'_d/\sim'$ via the canonical identification with $\mathcal{E}_d/\sim$. 
\end{remark}

\subsection{Fibration tower of embedding spaces}
\label{subsec_fib_tower_emb}
Fix a handle decomposition of $\Sigma$ with a single $0$-handle $H_0$ that contains $b_0$, $1$-handles $H_{1}^1,\cdots, H_{1}^{2\ell}$, and a single $2$-handle $H_2$. Let $H_{1}$ be the union of the $1$-handles. By definition, the handles are codimension-0 submanifolds of $\Sigma$ with corners. 
Let $e^{0}, e^{1}_{1},\cdots, e^{1}_{2\ell}, e^{2}$ be the cells of the corresponding CW structure of $\Sigma$, where each $1$-cell deformation retracts to the core of the corresponding $1$-handle, and we require that $e^0 = \{b_0\}$. Let $U\subset S^2$ be a disk containing $b_1$ in its interior. Let $\nu(e^0) = H_0\times U$, let $\nu(e_j^1) = (H_0\cup H_{1}^{j})\times U$. Then $\nu(e^0)$ is a tubular neighborhood of $\mathscr{I}_d(e^0)$, and $\nu(e_j^1)$ is a tubular neighborhood of $\mathscr{I}_d(e^0\cup e_j^1)$. We also assume that the embeddings $\mathscr{I}_d$ satisfy $\mathscr{I}_d|_{H_0\cup H_1} = \mathscr{I}_0|_{H_0\cup H_1}$.

Let  $\pi$ denote the fundamental group of $\Sigma$ with base point in $H_2$.

\begin{Definition}
Let 
$M_1$ be the closure of $M\setminus \nu(e^0)$, 
and let 
$M_2$ be the closure of $M\setminus (\cup^{2\ell}_{j=1}\nu(e^{1}_{j})))$.
We smooth corners and regard both $M_1$ and $M_2$ as codimension-0 submanifolds of $M$. 
\end{Definition}

\begin{Definition}
Suppose $A$ is a codimension-0 submanifold of $\Sigma$ that contains $b_0$. Let $\Emb^{[\mathscr{I}_d]}_\bullet(A, M)$ denote the space of smooth embeddings $i:A\hookrightarrow M$ such that $i(b_0)=b$ and $i$ is homotopic to  $\mathscr{I}_d|_{A}$ relative to $b_0$. Set $\mathscr{I}_d|_A$ to be the base point of $\Emb^{[\mathscr{I}_d]}_\bullet(A, M)$.
\end{Definition}

By definition, we have $\mathcal{E}_d = \Emb^{[\mathscr{I}_d]}_\bullet (\Sigma, M)$. For $i_1,i_2\in \mathcal{E}_d$, we have $i_1\sim i_2$ if and only if they are on the same connected component of $\mathcal{E}_d$. 

Consider the fibration tower 
\begin{equation}\label{eq: fibration tower}
  \Emb^{[\mathscr{I}_d]}_\bullet (\Sigma, M)\xrightarrow{r_2} \Emb^{[\mathscr{I}_d]}_\bullet(H_0\cup H_{1}, M)\xrightarrow{r_1} \Emb^{[\mathscr{I}_d]}_\bullet(H_{0}, M)\xrightarrow{r_0} \ast,
\end{equation}
where the maps $r_1,r_2$ are given by restrictions.  For $j=0,1,2$, let $F_j$ be the preimage of the base point under $r_j$. Then $r_j$ is a fibration and $F_j$ is the fiber of $r_j$.  

\begin{remark}
We will define a Dax invariant on $\pi_0\Emb^{[\mathscr{I}_d]}_\bullet (\Sigma, M)$. Our strategy is to first define a Dax invariant on $\pi_0 (F_2)$, and then study the properties of the Dax invariant when two elements of $\pi_0 (F_2)$ map to the same element in $\pi_0\Emb^{[\mathscr{I}_d]}_\bullet (\Sigma, M)$. The definition will be formulated in Section \ref{subsec_Dax_closed_surface}. 
\end{remark}

To describe the homotopy type of $F_j$, we introduce the following notations.

\begin{Definition}
\label{defn_i_0^(1)}
Let $\mathscr{I}_d^{(1)}$ denote the restriction of $\mathscr{I}_d$ to $(\cup_{j=1}^{2\ell} e_j^1)\cap \mathscr{I}_d^{-1}(M_1)$. Then $\mathscr{I}_d^{(1)}$ is a neat embedding of $2\ell$ disjoint arcs into $M_1$. Let $Dom(\mathscr{I}_d^{(1)})$ denote the domain of $\mathscr{I}_d^{(1)}$. Then $Dom(\mathscr{I}_d^{(1)})$ consists of $2\ell$ smoothly embedded disjoint arcs in $\Sigma$. Let $n(\mathscr{I}_d^{(1)})$ be a fixed non-vanishing section of the normal bundle of $Dom(\mathscr{I}_d^{(1)})$ in $\Sigma$, such that  $n(\mathscr{I}_d^{(1)})$ is tangent to the boundary of $\mathscr{I}_d^{-1}(M_1)$ on $\partial Dom(\mathscr{I}_d^{(1)})$.
\end{Definition}

\begin{Definition}
Identify $Dom(\mathscr{I}_d^{(1)})$ with $\sqcup_{2\ell}I$.
Define
$\Emb'_{\partial}(\sqcup_{2\ell}I, M_1)$
to be the space of neat embeddings of $2\ell$ disjoint arcs into $M_1$, each equipped with a non-vanishing normal section (i.e. a section of the normal bundle), such that the images near the $4\ell$ boundary points are equal to the restriction of $\mathscr{I}_d^{(1)}$, and the normal section near the boundary points are given by the images of $n(\mathscr{I}_d^{(1)})$ under the tangent map of $\mathscr{I}_d$. 
\end{Definition}

\begin{Definition}
\label{defn_Emb_D2_M2}
Identify $H_2$ with the standard disk $D^2$.  
Define $\Emb_{\partial}(D^2,M_2)$ by setting $X=M_2$ and setting $i_0=\mathscr{I}_d|_{H_2}$ in Definition \ref{defn_Emb_D_X}.
\end{Definition}

\begin{Definition}
\label{defn_Emb_D2_M2_hom}
Define $\Emb_{\partial}^{hom}(D^2,M_2)$ to be the subspace of $\Emb_{\partial}(D^2,M_2)$ consisting of elements $i\in \Emb_{\partial}(D^2,M_2)$ such that the fundamental class of $i(D^2)$ equals the fundamental class of $\mathscr{I}_d(H_2)$ in $H_2(M_2,\partial M_2)$.
\end{Definition}

By definition, $\Emb_{\partial}^{hom}(D^2,M_2)$ is the union of a collection of connected components of $\Emb_{\partial}(D^2,M_2)$. 
We have the following result.

\begin{Lemma}
\label{lem_hmty_type_Fj}
The homotopy types of $F_j$ can be described as follows.
    \begin{enumerate}
        \item $F_0$ is homotopy equivalent to the space of pairs of orthogonal unit vectors in $T_{b}M$. In particular, $\pi_{1}(F_0)\cong 1$.
        \item $F_1$ is homotopy equivalent to the component of $\Emb'_{\partial}(\sqcup_{2\ell}I, M_1)$ containing $\mathscr{I}_d^{(1)}$.
        \item $F_2$ is homotopy equivalent to $\Emb_{\partial}^{hom}(D^2,M_2)$. 
    \end{enumerate}
\end{Lemma}

\begin{proof}
    Most of this lemma follows from standard arguments. The only non-standard part is to verify the following: if $s\in \Emb_{\partial}(D^2,M_2)$, identify the domain of $s$ with $H_2\subset \Sigma$ and extend $s$ to a map from $\Sigma$ to $M$ by defining the extension on $H_0\cup H_1$ to be $\mathscr{I}_d|_{H_0\cup H_1}$, then the extension is homotopic to $\mathscr{I}_d$ relative to $b_0$ if and only if $s\in \Emb_{\partial}^{hom}(D^2,M_2)$. 
    For later reference, we prove it in two separate lemmas below. The ``only if'' part is given by Lemma \ref{lem_F2_in_hom}, and the ``if'' part follows immediately from Lemma \ref{lem_hom_in_F2}.
\end{proof}

\begin{Lemma}
\label{lem_F2_in_hom}
Suppose $i\in F_2$, then the fundamental class of $i(H_2)$ equals the fundamental class of $\mathscr{I}_d(H_2)$ in $H_2(M_2,\partial M_2)$.
\end{Lemma}

\begin{proof}
Since $H_2(M_2,\partial M_2)$ has no torsion, by Poincar\'e duality, we only need to show that for each $\xi \in H_2(M_2)$, its intersection number with $\mathscr{I}_d(H_2)$ and $i(H_2)$  are the same. The intersection number of $\xi$ with $i(H_2)$ in $M_2$ equals the intersection number of $\xi$ with $i(\Sigma)$ in $M$, and the same holds when replacing $i$ with $\mathscr{I}_d$. Since $\mathscr{I}_d$ is homotopic to $i$ in $M$, the lemma is proved. 
\end{proof}

\begin{Lemma}
\label{lem_hom_in_F2}
Suppose $i:(\Sigma,b_0)\to (M,b)$ is a map such that $\mathscr{I}_d$ and $i$ induce the same maps on $\pi_1(\Sigma,b_0)\to \pi_1(M,b)$ and $H_2(\Sigma)\to H_2(M)$. Then $i$ is homotopic to $\mathscr{I}_d$ relative to $b_0$.
\end{Lemma}

\begin{proof}
We use $\langle\cdot, \cdot\rangle$ to denote
the set of pointed homotopy classes.
 We have
$$
\langle\Sigma, M \rangle=\langle\Sigma, S^2\rangle
\times \langle \Sigma, \Sigma\rangle.
$$
The component in $\langle\Sigma, S^2\rangle$
is determined by the induced map on the second homology group because $S^2$ is equal to the 3-skeleton of $\mathbb{C}P^\infty$, which is the Eilenberg--MacLane space $K(\bZ,2)$.
The component in  $\langle \Sigma, \Sigma\rangle$
is determined by the induced map
on $\pi_1$
because $\Sigma$ is
$K(\pi_1(\Sigma),1)$.
\end{proof}

Since $F_0$ is simply connected, the inclusion from $F_1$ to $\Emb^{[\mathscr{I}_d]}_\bullet(H_0\cup H_{1}, M)$ induces a surjection on $\pi_1$. 
The fibration tower gives an action of $\pi_{1}(F_1)$ on $\pi_{0}(F_2)$. 
Since $\Emb^{[\mathscr{I}_d]}_\bullet(H_0\cup H_{1}, M)$ is connected, we have a bijection
\begin{equation}\label{eq: q_0}
\pi_0(\Emb^{[\mathscr{I}_d]}_\bullet(\Sigma, M))\to \pi_{0}(F_2)/\pi_{1}(F_1).\end{equation}

\begin{Definition}
Let $\xi_0:I\to D^2$ be a fixed neat embedding such that its image is contained in $N(\partial D^2)$, where $N(\partial D^2)$ is the neighborhood of $\partial D^2$ as in Definition \ref{defn_Emb_D_X}. As in Definition \ref{defn_Emb_D2_M2}, we identify $H_2$ with $D^2$. 

Let $I_0$ be the composition map $\mathscr{I}_d\circ \xi_0$. We fix a framing $s_0$ on $I_0$.  
Let $\Emb_\partial(I,M_2)$ and $\Emb^{\Fr}_{\partial}(I,M_2)$ be the space of (framed) embeddings of as in Definition \ref{defn_Emb_I}.  
\end{Definition}

\subsection{Properties of $\pi_1(\Emb_{\partial}(I, M_2);I_0)$}
\label{subsec_Dax_pi1_EmbIM}
This subsection studies the properties of $\pi_1(\Emb_{\partial}(I, M_2);I_0)$. The computation will later be used to compute $\pi_1(F_1)$ in \eqref{eq: q_0}.

We consider two collections of elements in $\pi_{1}(\Emb_{\partial}(I, M_2);I_0)$:
\begin{enumerate}
    \item For each $g\in \pi_{1}(M_2;I_0)$,  we let $\tau_{g}$ be the spinning family around $I_0$, as defined in Definition \ref{defi: spinning I0}.
    \item Let $b'_0$ be a fixed point in $H_2\subset \Sigma$, Let $S_0=\{b'_0\}\times S^2$, suppose $b_0'$ is chosen so that $S_0$ is disjoint from $I_0$. For $g\in \hat{\pi}_{1}(M_2;I_0,S_0)$, we let $\rho_{g}$ be the spinning family around $S_0$ as defined in Definition \ref{defi: general spining family}. 
\end{enumerate}

\begin{Lemma}
\label{lem_pi1Emb_commute}
The elements in $\pi_1\Emb_{\partial}(I, M_2)$ given by $\tau_{g}$ and $\rho_{\gamma}$
are commutative with each other. 
\end{Lemma}
\begin{proof}
Given $g_{1}\in \pi_{1}(M_2;I_0)$ and $g_{2}\in \pi_{1}(M_2;I_0)$, we represent them by paths $\gamma_{1}$ and $\gamma_{2}$  respectively such that $\gamma_1\cap \gamma_2=\emptyset$. 
Then there is a map $T^2\to \Emb_{\partial}(I,M_2)$, such that the first $S^1$ factor of $T^2$ parametrizes a spinning along $\gamma_1$, and the second $S^1$ factor parametrizes a spinning along $\gamma_2$.  So $\tau_{g_1}$ and $\tau_{g_2}$ commute. The other cases are similar.\end{proof}

We will show in Corollary \ref{cor_pi1Emb_abel} below that elements of forms $\tau_{g}$ and $\rho_{\gamma}$ generate $\pi_1(\Emb_{\partial}(I,M_2);I_0)$, so $\pi_1(\Emb_{\partial}(I,M_2);I_0)$ is abelian.

In the following, let $D_0$ denote $\mathscr{I}_0(H_2)$, then $D_0$ is a neatly embedded disk in $M_2$. 
Note that the embedding of $M_2$ into $M$ induces an isomorphism on $\pi_1$, so $\pi_1(M_2;I_0)$ is canonically isomorphic to $\pi$. 

\begin{Lemma}\label{lem: pi2} We have an isomorphism $\pi_{2}(M_2;I_0)\cong \mathbb{Z}[\pi]$ defined by sending $[S]\in \pi_{2}(M_2;I_0)$ to its equivariant intersection number with $[D_0]$. A basis of $\pi_{2}(M_2;I_0)$ as a free abelian group is given by the following. For each $g\in \pi_1(M_2;I_0,S_0)$, let $\gamma_g$ be an embedded arc from $I_0$
to $S_0$ such that $[\gamma_g]=g$ and $\inte(\gamma)$ is disjoint from $I_0\cup S_0$. 
Then $\{\mathcal{F}(\rho_{\gamma_g})\}_{g\in \pi_1(M_2;I_0,S_0)}$ is a basis of $\pi_{2}(M_2;I_0)$, where $\mathcal{F}$ is defined by \eqref{eqn_pi1Emb_to_pi2}, and $\rho_{\gamma_g}$ is defined by Definition \ref{defi: general spining family}.
\end{Lemma}

\begin{remark}
    Note that although $\rho_{\gamma_g}$ depends on the choice of $\gamma_g$, the value of $\mathcal{F}(\rho_{\gamma_g})$ only depends on $g$ and does not depend on the choice of $\gamma_g$.
\end{remark}

\begin{proof}
  Let $\widetilde{M}_2$ be the universal cover of $M_2$ and let $\widetilde{\Sigma}$ be the universal cover of $\Sigma$. Let $\operatorname{sk}_{1}\widetilde{\Sigma}$ denote the pre-image of $e^0\cup (\cup_{j=1}^{2\ell} e_j^1)$ in $\widetilde{\Sigma}$.  Then 
  \[\widetilde{M}_2\simeq (\widetilde{\Sigma}\times S^2)\setminus ( \operatorname{sk}_{1}\widetilde{\Sigma}\times \{b_1\}).\] 
  By a straightforward application of the Mayer-Vietoris sequence, we see that the lifts of $S_0$ give a set of free generators for $\pi_{2}(M)\cong H_{2}(\widetilde{M}_2;\mathbb{Z})$, and the result is proved.
\end{proof}

\begin{remark}
    Note that the equivariant intersection number of $\mathcal{F}(\rho_\gamma)$ and $[D_0] = [\mathscr{I}_0(H_2)]$ is equal to the equivariant intersection number of $\mathcal{F}(\rho_\gamma)$ and $[\mathscr{I}_d(H_2)]$ for all $d$. 
\end{remark}

The next lemma is essentially proved in \cite[Theorem 0.3(2)]{gabai2021self}. 
\begin{Lemma}\label{lem: pi1D generators}
$\pi^{D}_{1}(\Emb_{\partial}(I,M_2);I_0)$ is generated by $\{\tau_{g}\mid g\in \pi_{1}(M_2;I_0)\}.$
\end{Lemma}
\begin{proof}
First, we note that each $\tau_g$ is in the center of $\pi_{1}(\Emb_{\partial}(I,M_2);I_0)$. In fact, this follows from \cite[Theorem 1.1]{kosanovic2024homotopy}. To verify the claim directly, observe that for each $[\theta]\in \pi_{1}(\Emb_{\partial}(I,M_2);I_0)$, we may isotope by re-parameterizing $I$ so that for every arc on the representative $\theta$, it coincides with $I_0$ on $[1/3,1]$. Similarly, we choose a representative of $\tau_g$ so that every point on the representative, which is an arc in $M_2$, coincides with $I_0$ on $[0,2/3]$. Since the trajectory of $\tau_g$ is contained in a regular neighborhood of a $1$--dimensional complex, we may further isotope the representatives of $\theta$ and $\tau_g$ so that the trajectory of $\theta$ on $[0,1/3]$ and the trajectory of the representative of $\tau_g$ on $[2/3,1]$ are disjoint. As a result, $[\theta]$ and $\tau_g$ are commutative in $\pi_{1}(\Emb_{\partial}(I,M_2);I_0)$.

Given an element $[\theta]\in \pi^{D}_{1}(\Emb_{\partial}(I,M_2);I_0).$ Since homotopies of arcs in a 4-manifold can be perturbed to \emph{regular} homotopies, we can extend $\theta$ to a map $h: D^2\to \operatorname{Imm}_{\partial}(I, M_2)$, where $\operatorname{Imm}_{\partial}(I, M_2)$ denotes the space of immersions that agree with $I_0$ near the boundary. By a general position argument, we can choose $h$ such that $h(x)\in \Emb_{\partial}(I,M_2)$ for all $x\in D^2\setminus \{x_{1},\cdots x_{m}\}$. Moreover, we may assume that each $h(x_{j})$ is an immersion with a single double point, and that the image of $h$ on a small circle around $x_{j}$ is a resolution of the double point. Then $[\theta]$ is the sum of some conjugations of the spinning families corresponding to these immersions.  Since each spinning family is in the center of $\pi_{1}(\Emb_{\partial}(I,M_2);I_0)$, the lemma is proved. 
\end{proof}

\begin{Corollary}
\label{cor_pi1Emb_abel}
    The group $\pi_{1}(\Emb_{\partial}(I,M_2);I_0)$ is generated by elements of the form $\tau_g$ and
    $\rho_\gamma$, and hence is abelian.
\end{Corollary}

\begin{proof}
    By definition, we have the following exact sequence
\[
\pi_{1}^D(\Emb_{\partial}(I,M_2);I_0) \hookrightarrow \pi_{1}(\Emb_{\partial}(I,M_2);I_0) \xrightarrow{\mathcal{F}} \pi_2(M_2;I_0).
\]
By Lemma \ref{lem: pi2}, elements of the form $\mathcal{F}(\rho_{\gamma_g})$, where $g\in \pi_1(M_2;I_0,S_0)$, generate $\pi_2(M_2;I_0)$. By Lemma \ref{lem: pi1D generators}, elements of the form $\tau_{g}$, where $g\in \pi_{1}(M_2;I_0)$, generate $\pi_{1}^D(\Emb_{\partial}(I,M_2);I_0)$. Hence $\pi_{1}(\Emb_{\partial}(I,M_2);I_0)$ is generated by elements of the form $\rho_\gamma$ and $\tau_g$. By Lemma \ref{lem_pi1Emb_commute}, we conclude that the group $\pi_{1}(\Emb_{\partial}(I,M_2);I_0)$ is abelian.
\end{proof}

Now we consider the Dax isomorphism (\ref{eq: Dax homomorphism}), with $X=M_2$. 
Note that $\pi_3(M_2;I_0)$ is generated by Whitehead products of elements in $\pi_2(M_2)$.
By \cite[Proposition 3.14]{KT24}, we have
$$
d_3([a,b])=\lambda(a,b)+\lambda(b,a)
$$
where $[a,b]$ denotes the Whitehead
product and 
$\lambda$ denotes the equivariant
intersection form omitting
the coefficient of $1$.
Since the equivariant intersection form of $\pi_2(M_2;I_0)$ is trivial,  we have $d_3=0$. Therefore, we can rewrite the Dax isomorphism as 
\[
\Dax: \pi^{D}_{1}(\Emb_{\partial}(I,M_2);I_0)\xrightarrow{\cong} \mathbb{Z}[\pi\setminus 1].
\]

\subsection{Properties of $\pi_1(F_1)$}
\label{subsection_pi1F1}
Recall that $F_2$ is the fiber of the map $r_2$ in the fibration tower \eqref{eq: fibration tower}. By Part (3) of Lemma \ref{lem_hmty_type_Fj}, we know that $F_2$ is homotopy equivalent to $\Emb_{\partial}^{hom}(D^2,M_2)$. 

Now, we introduce  three collections of elements in 
\[\pi_{1}(F_1)=\pi_{1}(\Emb'_{\partial}(\sqcup_{2\ell}I,M_1)).\]
Recall from Definition \ref{defn_i_0^(1)} that $\mathscr{I}_d^{(1)}$ denotes the restriction of $\mathscr{I}_d$ to $(\cup_{j=1}^{2\ell} e_j^1)\cap \mathscr{I}_d^{-1}(M_1)$, and that $n(\mathscr{I}_d^{(1)})$ is a fixed non-vanishing normal vector field of the domain of $\mathscr{I}_d^{(1)}$ in $\Sigma$. Let the base point of $\pi_{1}(\Emb'_{\partial}(\sqcup_{2\ell}I,M_1))$ be given by $\mathscr{I}_d^{(1)}$ and the image of $n(\mathscr{I}_d^{(1)})$ under the tangent map of $\mathscr{I}_d$. 

For notational convenience, let $I_j$ denote the image of $e_j^1\cap \mathscr{I}_d^{-1}(M_1)$ under the map $\mathscr{I}_d$, and let $s_j$ denote the image of $n(\mathscr{I}_d^{(1)})$ restricted to $I_j$. Then $s_j$ is a normal vector field on $I_j$. Let $\Emb_{\partial}(\sqcup_{2\ell}I,M_1)$ be the space of embeddings of $\sqcup_{2\ell}I$ in $M_1$ that agree with $\mathscr{I}_d^{(1)}$ near the boundary points, and let the base point of $\Emb_{\partial}(\sqcup_{2\ell}I,M_1)$ be given by $\mathscr{I}_d^{(1)}$.

Recall that we set  $S_0=\{b'_{0}\}\times S^2$ for $b'_{0}\in H_2$. 

\begin{Definition}\label{defi: spinning around S0}
Fix $j\in\{1,\dots,2\ell\}$. 
Given an embedded arc $\gamma$ from $I_j$ to $S_0$ with $\inte (\gamma)\cap S_0=\emptyset$ and $\inte (\gamma)\cap I_u=\emptyset$ for all $u$,
we construct a loop in $\Emb'_{\partial}(\sqcup_{2\ell}I,M_1)$ similar to Definition \ref{defi: spinning I0 to I1} as follows. For $u\neq j$, we keep $I_{u}$ fixed. For $u=j$, we take the spinning family of $I_{j}$ along $(S_0,\gamma)$. As discussed in Definition \ref{defi: general spining family}, the loop of embeddings of $I_j$ lifts to a unique loop $\widetilde{\rho}_{\gamma}$ of framed embeddings up to homotopy. Hence we obtain an element in $\pi_{1}(\Emb'_{\partial}(\sqcup_{2\ell}I,M_1))$. We denote this element by $\hat \rho_{\gamma}\in \pi_{1}(\Emb'_{\partial}(\sqcup_{2\ell}I,M_1))$. 
   \end{Definition}

\begin{Definition}
For $g\in \pi_{1}(M_1;I_{u},I_{v})$, consider the spinning family as follows. 
Represent $g$ by a path $\gamma$ with end points $\gamma(0)\in I_u,\gamma(1)\in I_v$. If $u=v$, we require that $\gamma(0)<\gamma(1)$ under the parametrization of $I_u$. We also assume $\gamma$ is embedded and its interior doesn't intersect $\sqcup_{j=1}^{2\ell} I_j$. For $j\neq u$, we keep $I_{j}$ fixed.
If $u=v$, we let $I_{u}$ be the spinning family as in Definition \ref{defi: spinning I0}. If $u\neq v$, we let $I_{u}$ be the spinning family as in Definition \ref{defi: spinning I0 to I1}.  Note that this loop is canonically null homotopic in the space of immersions, so it has a unique lifting to 
$\Emb_{\partial}'(\sqcup_{2\ell}I,M_1)$ up to homotopy if the lifting of the endpoints is given by $n(\mathscr{I}_d^{(1)})$. This defines an element $\hat \tau_{g}\in \pi_{1}(\Emb'_{\partial}(\sqcup_{2\ell}I,M_1))$. 

%
\end{Definition}

\begin{figure}
	\begin{overpic}[width=0.6\textwidth]{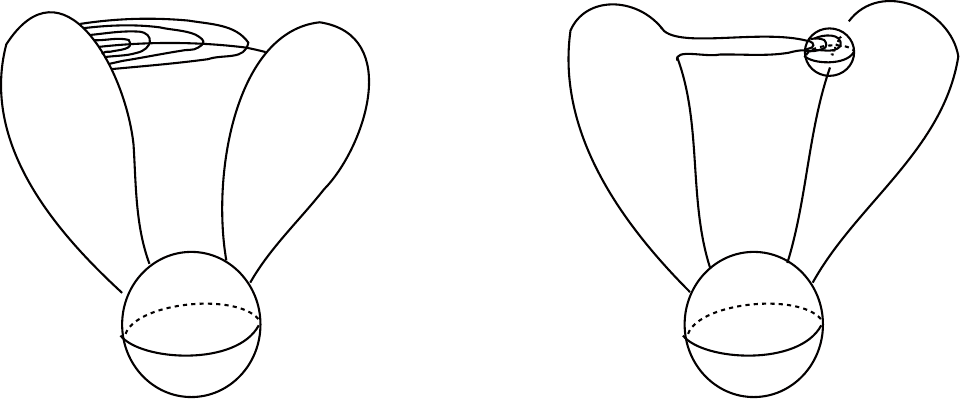}
		\put(62,12){$\partial M_1$}
            \put(190,12){$\partial M_1$}
            \put(-10,70){$I_u$}
            \put(85,70){$I_v$}
            \put(220,70){$I_v$}
	\end{overpic}
	\caption{The Definition of $\hat \tau_g$}
    \label{fig:def_hat_tau_g}
\end{figure}

Figure \ref{fig:def_hat_tau_g} describes the definition of $\hat \tau_g$. 

\begin{Definition}
We will abuse notation and denote the projections of $\hat{\rho}_\gamma,\hat{\tau}_{g}\in \pi_{1}(\Emb_{\partial}'(\sqcup_{2\ell}I,M_1))$ to $\pi_{1}(\Emb_{\partial}(\sqcup_{2\ell}I,M_1))$ by $\rho_\gamma$ and $\tau_{g}$ when there is no risk of confusion with Definitions \ref{defi: spinning I0} and Definition \ref{defi: spinning I0 to I1}. 
\end{Definition}

\begin{Definition}
\label{defn_hat_xi_j}
    For each $1\leq j\le 2\ell$, we fix a trivialization of the normal bundle for $I_{j}$ in $M_1$, so that $s_{j}$ corresponds to a constant section. Then a normal section along $I_{j}$ is just a point in $\Omega S^2$. We take a loop of normal sections $s^{\theta}_{j}$ ($\theta\in [0,1]$) on $I_{j}$ that represents a generator of $\pi_{1}(\Omega S^2)\cong \mathbb{Z}$. Then the loop \[(\sqcup_{u} I_{u}, \sqcup_{u\neq j}s_{u}\sqcup s^{\theta}_{j})_{\theta\in [0,1]}\] defines an element in $\pi_{1}(\Emb'_{\partial}(\sqcup_{2\ell}I,M_1))$, which we denote by $\hat \xi_{j}$. 
\end{Definition}

Now we show that $\pi_{1}(\Emb'_{\partial}(\sqcup_{2\ell}I,M_1))$ is abelian and is generated by elements of the form $\hat{\rho}_\gamma$ and $\hat{\xi}_j$.

\begin{figure}
	\begin{overpic}[width=0.6\textwidth]{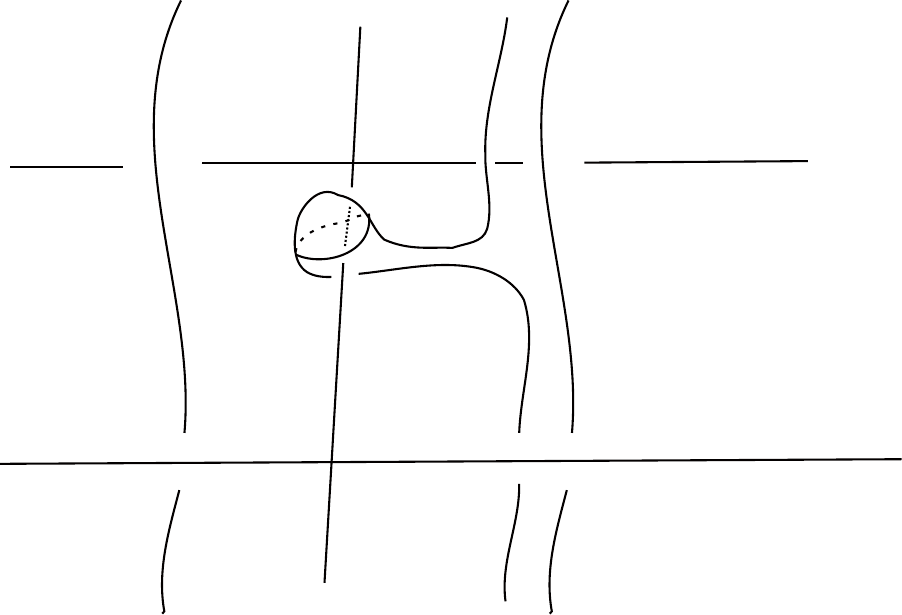}
		\put(23,120){$S_0$}
		\put(135,120){$S_0$}
		\put(65,15){$e_k^1$}
            \put(55,90){$m_k$}
            \put(105,75){$\gamma$}
		\put(95,145){$m_k\# S_0$}
	\end{overpic}
	\caption{The construction of $m_k\# S_0$}
    \label{fig:rho_g_pm}
\end{figure}

\begin{Lemma}
\label{lem_rho_diff_tau}
For any $g\in \pi_{1}(M_1;I_u, I_{v})$, there exist embedded arcs $\gamma_+$ and $\gamma_-$ from $I_u$ to $S_0$ whose interiors are disjoint from $(\sqcup I_u) \cup S_0$, such that
$\rho_{\gamma_{+}}-\rho_{\gamma_-} = \pm \tau_{g}$.  
\end{Lemma}
\begin{proof}
Recall that $e^2$ denotes the 2-cell of $\Sigma$. 
We take an arc $\gamma\subset \mathscr{I}_d(e^2)$ connecting $m_{k}$ to $S_0$ and add a tube around $\gamma$ to form the connected sum $m_k\# S_0$. See Figure \ref{fig:rho_g_pm} for a description of the construction of $m_k\# S_0$ on the universal cover.

The sphere $m_k\# S_0$ is still isotopic to $S_0$ in $M_2$, but the trajectory of a base point on $S_0$ under this isotopy is non-trivial. This shows that the sum of a suitable spinning family around $m_k$ with a suitable spinning family around $S_0$ is another spinning family around $S_0$. A straightforward tracking of the corresponding arcs shows the desired lemma. 
\end{proof}

\begin{Lemma}\label{lem: pi1_emb_H_1_generators}
The group $\pi_{1}(\Emb_{\partial}(\sqcup_{2\ell}I,M_1))$ is generated by elements of the form $\rho_{\gamma}$, and hence is abelian.
\end{Lemma}
\begin{proof}
This proof is similar to Lemma \ref{lem: pi1D generators} and Corollary \ref{cor_pi1Emb_abel}. For each $j$, there is a homomorphism 
\[
\mathcal{F}_j: \pi_1 (\Emb_{\partial}(\sqcup_{2\ell}I,M_1)) \to \pi_2(M_1;I_j)
\]
defined by taking the trajectory of the image of $I_j$. Let 
\[
\hat{\mathcal{F}}: \pi_1 (\Emb_{\partial}(\sqcup_{2\ell}I,M_1)) \to \prod_{j=1}^{2\ell} \pi_2(M_1;I_j)
\]
be the product of the maps $\mathcal{F}_j$, then elements of the form $\hat{\mathcal{F}}( \rho_\gamma)$ generate the codomain of $\hat{\mathcal{F}}$. Let $\pi_1^D (\Emb_{\partial}(\sqcup_{2\ell}I,M_1))$ denote the kernel of $\hat{\mathcal{F}}$, then the same argument as Lemma \ref{lem: pi1D generators} shows that $\pi_1^D (\Emb_{\partial}(\sqcup_{2\ell}I,M_1))$ is generated by elements of the form $\tau_{g}$. By Lemma \ref{lem_rho_diff_tau}, we know that $\pi_{1}(\Emb_{\partial}(\sqcup_{2\ell}I,M_1))$ is generated by elements of the form $\rho_{\gamma}$.
The same argument as Lemma \ref{lem_pi1Emb_commute} shows that $\rho_{\gamma}$ are commutative to each other, so the result is proved.
\end{proof}

\begin{Lemma}\label{lem: pi1F1 generators}
The group $\pi_{1}(\Emb'_{\partial}(\sqcup_{2\ell}I,M_1))$ is generated by elements of the form $\hat \rho_{\gamma}$ and $\hat \xi_{j}$ ($1\leq j\leq 2\ell)$, and hence is abelian. 
\end{Lemma}
\begin{proof}
The projection $\Emb'_{\partial}(\sqcup_{2\ell}I,M_1) \to \Emb_{\partial}(\sqcup_{2\ell}I,M_1)$ is a fibration and the fundamental group of the fiber is generated $\hat \xi_{j}$ $(1\le j\le 2\ell)$. By Lemma \ref{lem: pi1_emb_H_1_generators}, 
$\pi_{1}(\Emb'_{\partial}(\sqcup_{2\ell}I,M_1))$ is generated by $\hat \xi_{j}$'s
and elements of form $\hat \rho_{\gamma}$.

It is clear that $\{\hat \xi_j\}_{1\le j\le 2\ell}$ are commutative to each other, and $\{\hat \rho_\gamma\}_\gamma$ are commutative to each other. 
To show that $\hat \xi_j$ and $\hat \rho_\gamma$ commute, we may assume
the normal sections $s_j^\theta$ ($\theta\in [0,1]$) in the definition of $\hat \xi_j$ only varies with respect to $\theta$ over a small interval 
$\hat{I}_j \subset I_j$ and $\gamma$ is disjoint from $\hat{I}_j$. 
Then there is a map $T^2\to \Emb'_{\partial}(\sqcup_{2\ell}I,M_1)$, such that the first $S^1$ factor of $T^2$ parametrizes a spinning  along $\gamma$, and the second $S^1$ factor parametrizes the variation of $s_j^\theta$. Hence $\hat \xi_j$ and $\hat \rho_\gamma$ commute.
\end{proof}

\begin{Definition}
\label{defn_beta_gamma}
    Let $\hat\rho_\gamma\in \pi_1\Emb'_{\partial}(\sqcup_{2\ell}I,M_1)$ be as above. Then $\hat \rho_\gamma$ lifts to a loop of \emph{framed} arcs, and the lifting is unique up to homotopy.  Define $\beta_\gamma \in \pi_0\Diff_\partial (M_2)$ to be the image of the lifting of $\hat\rho_\gamma$ under the arc pushing map 
    \[
    \pi_1(\Emb^{\Fr}(\sqcup_{2\ell}I,M_1))\to \pi_0\Diff_\partial (M_2).
    \]
\end{Definition}

\subsection{The action of $\pi_1(F_1)$ on $\pi_0(F_2)$}
\label{subsection_pi1F1_action_homotopy}
The action of $\pi_1(F_1)$ on $\pi_0(F_2)$ can be described as follows. 
Given an element 
\[
\alpha\in \pi_1(F_1)\cong \pi_1\Emb'_{\partial}(\sqcup_{2\ell}I,M_1),
\]
we pick an isotopy $f_t:H_1\to M_1$ relative to $H_0$ representing
$\alpha$
and extend it to an isotopy $F_t:M\to M$ with $F_1|_{(\nu(H_0)\cup \mathscr{I}_0(H_1))}=\id_{(\nu(H_0)\cup \mathscr{I}_0(H_1))}$. Then
the action of $\alpha$ on $\pi_0(F_2)=\pi_0(\Emb_{\partial}^{hom}(H_2,M_2))$ is given by $[i]\mapsto [F_1\circ i]$.
In particular, the action of $\hat \rho_\gamma\in \pi_1\Emb'_{\partial}(\sqcup_{2\ell}I,M_1)$ is given by the barbell diffeomorphism $\beta_\gamma$ defined in Definition \ref{defn_beta_gamma}. 

We use $D_0$ to denote the embedded disk $H_2\times \{b_1\}\hookrightarrow M_2$ and use $S_0$ to denote the embedded sphere $\{b'_0\}\times S^2\hookrightarrow M_2$. We use $b'$ to denote their intersection $(b_0',b_1)$.  Pick an arbitrary point $b_{0}''\in \partial H_2$ and set $b''=(b_0'',b_1)\in \partial D_0$ as the base point, then $\pi_{1}(M_2, b'')$ is canonically isomorphic to $\pi$. Since $D_0$ is simply-connected, every path in $M_2$ from $D_0$ to itself represents an element in $\pi$. We use $\pi_{2}(M_2;\partial D_0)$ to denote the set of homotopy classes of maps $D^2\to M_{2}$ that coincide with $\mathscr{I}_0|_{\partial H_2}$ on the boundary. 

Now we recall the definition of the equivariant intersection number 
\[
\lambda: \pi_{2}(M_2;\partial D_0)\times \pi_{2}(M_2;\partial D_0)\to \mathbb{Z}[\pi].
\]
Take two elements in $\pi_{2}(M_2;\partial D_0)$ represented by neatly immersed disks $i_{1},i_{2}: D^2\to M_{2}$. We pick a nonzero tangent vector $v\in T_{b_{1}}S^{2}$, which gives a section of the normal bundle of $D_0$ in $M_{2}$. Using this section, we may perturb $i_{2}$ near its boundary such that $i_{1}(D^2)$ and $i_{2}(D^2)$ only intersect outside a collar neighborhood of $\partial M_{2}$. After a further perturbation, we may assume $i_1$ and $i_2$ intersect transversely. For $j=1,2$ and  $x\in i_{1}(D^2)\cap i_{2}(D^2)$, we let $\gamma_{j}$ be a path in $D^{2}$ from $b_{0}''$ to $i^{-1}_{j}(x)$ . Then $(i_{1}\circ \gamma_{1})* \overline{(i_{2}\circ \gamma_{2})}$ defines a loop based at $b''$ and represents an element $g_{x}\in \pi$. We define 
\[
\lambda([i_{1}],[i_{2}]):=\sum_{x\in i_{1}(\mathring{D}^2)\cap i_{2}(\mathring{D}^2) } \operatorname{sign}(x)\cdot g_{x}\in\mathbb{Z}[\pi].    
\]
Here, $\operatorname{sign}(x)$ is the sign of the intersection point. Since every homotopy between immersed disks can be realized by a sequence of finger/Whitney moves, twist moves, and isotopies, and all these moves don't change the value of $\lambda([i_{1}],[i_{2}])$, the intersection number $\lambda$ is well-defined. 

Similarly, we recall the definition of the self-intersection number \[
\mu: \pi_{2}(M_2;\partial D_0)\to \mathbb{Z}[ \bar\pi].
\]
Here, $ \mathbb{Z}[\bar\pi]$ is the free abelian group generated by the set \[\bar \pi=(\bar \pi\setminus \{1\})/(g\sim g^{-1}).\] 
Consider $[i]\in \pi_{2}(M_2;\partial D_0)$, where $i$ is a neat immersion. For each double point $x$ of $i$ with $i^{-1}(x)=\{x_{1},x_{2}\}$, we pick a path $\gamma_{j}$ in $D^2$ from $b_0''$ to $x_{j}$ ($j=1,2$). Then $(i\circ \gamma_{1})* \overline{(i\circ \gamma_{2})}$ represents an element $g_{x}\in \pi$. Note that there is no canonical ordering of $x_1$ and $x_2$, and if we interchange $x_{1}$ and $x_{2}$, the element $g_{x}$ becomes $g^{-1}_{x}$. We define 
\[
\mu([i]):=\sum_x \operatorname{sign}(x)[g_{x}]\in \mathbb{Z}[ \bar\pi], 
\]
 where the sum is taken over all double points $x$ with $g_{x}\neq e$. Since every homotopy between immersed disks can be realized by a sequence of finger/Whitney moves, twist moves, and isotopies, and all these moves don't change the value of $\mu([i])$, the self-intersection number $\mu$ is well-defined. 

Given an element $\alpha=\sum^{n}_{k=1}\pm g_{k}\in \mathbb{Z}[ \pi]$, we construct an element  $[i_{\alpha}]\in \pi_{2}(M_{2};\partial D_0)$ as follows. We pick embedded spheres $S_{1},\cdots ,S_{n}$ parallel to $S_0$. For each $k$, we take an embedded arc in $M_{2}$ from $D_0$ to $S_{k}$. Since $D_0\cup S_{k}$ is simply-connected, the arc $\gamma_{k}$ represents an element of $\pi$. We pick $\gamma_{k}$ such that $[\gamma_{k}]=g_{k}$. Then we take a tube along $\gamma_{k}$ to connect $D_0$ to $S_{k}$ or its orientation reversal $\overline{S_{k}}$, depending on the sign of $g_{k}$. This gives an immersion $i_{\alpha}$. This construction gives a set-theoretic map 
\begin{align*}
\eta: \mathbb{Z} [\pi] &\to \pi_{2}(M_2;\partial D_0)\\
\alpha &\mapsto [i_{\alpha}].
\end{align*}
\begin{Lemma} 
\label{lem_lambda_0}
The map $\eta$ is a bijection. Its inverse is given by the map
\[\lambda_0:=\lambda(-,\mathscr{I}_0|_{H_2}):\pi_{2}(M_2;\partial D_0)\to \mathbb{Z}[ \pi].\]
\end{Lemma}
\begin{proof} We have a bijection $\pi_{2}(M_2;\partial D_0)\cong \pi_{2}(M_{2})$, which sends $[f]$ to the pointed homotopy class of 
\[
f\cup \mathscr{I}_0|_{H_2}: S^2\cong D^2\cup D^2\to M_{2}.
\]
 Note the isomorphism 
$\pi_{2}(M_{2})\cong H_{2}(\widetilde{M}_2),$
where $\widetilde{M}_2$ denotes the universal cover of $M_{2}$. 
Since $H_{2}(\widetilde{M}_2)$ is freely generated by lifts of $S_0$, we see that $\pi_{2}(M_{2})$ is freely generated by $\{g \cdot S_0\}_{g\in \pi}$. So $\eta$ is bijection. It is straightforward to verify that $\lambda\langle i_{\alpha},\mathscr{I}|_{H_2}\rangle =\alpha$. This finishes the proof.
\end{proof}

\begin{Lemma} 
\label{lem_relation_mu_lambda0}
The map 
$\mu\circ \eta: \mathbb{Z}[ \pi] \to \mathbb{Z}[ \bar\pi]$
is a homomorphism that sends $1\in \mathbb{Z}[ \pi]$ to $0$ and sends $g\in \mathbb{Z}[ \pi]$ $(1\neq g\in \pi)$ to the projection of $g$ in $\bar \pi$. 
\end{Lemma}
\begin{proof} It is straightforward to verify that the self-intersection number of the immersion $i_\alpha$ is equal to the quotient image of $\alpha$ in $\mathbb{Z}[\bar\pi]$. 
\end{proof}

Next, we study the effect of barbell diffeomorphisms on $\pi_{2}(M_2;\partial D_0)$. We are interested in the action of the barbell diffeomorphism $\beta_\gamma$ in Definition \ref{defn_beta_gamma}.
For $[i]\in \pi_{2}(M_{2};\partial D_0)$, we compute the difference 
\[
\lambda_0([\beta_\gamma \circ i])-\lambda_0([i])\in \mathbb{Z}[ \pi],
\]
where the map $\lambda_0$ is given by Lemma \ref{lem_lambda_0}. Recall that $I_k\subset M_1$ denotes the image of $e_k^1\cap \mathscr{I}_d^{-1}(M_1)$ under the map $\mathscr{I}_d$. Suppose $m_k$ is a meridian of $I_k$ and $\gamma$ is an arc connecting $m_k$ and $S_0$. 
We denote the positive and negative intersection point of $m_{k}$ with $D_0$ by $x_{+}$ and $x_{-}$ respectively. We homotope the end points of $\gamma$ along $m_{k}\cup S_0$ to obtain a path $\gamma_{+}$ from $x_{+}$ to $b'$ and a path $\gamma_{-}$ from $x_{-}$ to $b'$. Since $\inte(D_0)$ is simply-connected and $\gamma_\pm$ both have end points in $\inte(D_0)$, they define elements in $\pi$. We denote the corresponding elements by $g_{+}$ and $g_{-}$ respectively. 

\begin{Lemma}
\label{lem_change_of_lambda0_under_barbell}
For every $[i]\in \pi_{2}(M_2;\partial D_0)$, we have   \[\lambda_0([\beta_{\gamma}\circ i])-\lambda_0([i])=g_{+}-g_{-}-g^{-1}_{+}+g^{-1}_{-}.\] 
\end{Lemma}
\begin{proof}
The barbell 
$\mathcal{B}_{\gamma}:=\nu(m_{k})\cup \nu(\gamma)\cup \nu(S_0)$
intersects $D_0$ at three disks $D_{+}$, $D_{-}$ and $D_{S}$, with centers $x_{+}, x_{-}$ and $b'$ respectively. By Lemma \ref{lem_lambda_0}, the map $i$ is homotopic to $i_\alpha$ for some $\alpha\in \bZ[\pi]$. Hence after a suitable homotopy, we may assume \[i(D^2)\cap \mathcal{B}_{g}=D_0\cap \mathcal{B}_{g}=D_{+}\cup D_{-}\cup D_{S}.\]
Then $(\beta_\gamma\circ i)(D^2)$ is obtained by removing $D_{+},D_{-},D_{S}$ from $i(D^2)$ and replacing them by $\beta_\gamma(D_{+}), \beta_\gamma(D_{-})$ and $\beta_\gamma(D_S)$ respectively. Therefore,
$\lambda_0([\beta_\gamma\circ i])-\lambda_0([i])$
is equal to the sum of the equivariant intersection numbers of $\beta_\gamma(D_{+}), \beta_\gamma(D_{-}),\beta_\gamma(D_S)$ with $D_0$.
By the explicit description of the barbell diffeomorphism in \cite[Construction 5.3]{BG2019}, we see that $\beta_\gamma(D_{+})$ is obtained by tubing $D_{+}$ with a parallel copy of $S_0$  along $\gamma_{+}$. Hence, the equivariant intersection number contributed by the intersection of $\beta_\gamma(D_{+})$ and $D_0$ is $g_{+}$. Similarly, $\beta_\gamma(D_{-})$ is obtained by tubing $D_{-}$ with a parallel copy of $\overline{S_0}$, so the equivariant intersection number contributed by the intersection of $\beta_\gamma(D_{-})$ and $D_0$ is $-g_{-}$. On the other hand, the disk $\beta_\gamma(D_{S})$ is obtained by tubing $D_{S}$ with $\overline{m_{k}}$ along $\gamma$, so the corresponding contribution is $-g^{-1}_{+}+g^{-1}_{-}$. This finishes the proof.
\end{proof}

\begin{Proposition}\label{prop_homotopy_after_barbell}
Let $i:D^2\to M_{2}$ be an embedding that equals $\mathscr{I}_{0}$ near $\partial D^2$. Assume the nonequivariant intersection number $i(D^2)\cdot D_0=d\in \mathbb{Z}$, where we perturb the boundary of one of the disks in the direction of the $S^2$-factor of $M$ to define the intersection number. Then there exists a finite sequence of embedded arcs 
$\{[\gamma_{k}]\in \hat\pi_{1}(M;I_{l_{k}},S_0)\}_{1\leq k\leq n}$ such that
\[
[\beta^{\pm 1}_{\gamma_{n}}\circ\cdots \circ \beta^{\pm 1}_{\gamma_{1}}\circ i]=[\mathscr{I}_{d}|_{H_2}] \in \pi_2(M_2;\partial D_0).
\]
\end{Proposition}
\begin{proof}
Let $\alpha=\lambda_{0}([\mathscr{I}_{d}])-\lambda_0([i])\in \mathbb{Z}[\pi]$. By Lemmas \ref{lem_lambda_0} and \ref{lem_change_of_lambda0_under_barbell},
 it suffices to  find $\{g_{k}\}_{1\leq k\leq n}$ such that 
\[
\alpha=\sum^{n}_{k=1}\pm (g_{k,+}-g_{k,-}-g^{-1}_{k,+}+g^{-1}_{k,-}).
\]
Since both $\mathscr{I}_{d}$ and $i$ are embeddings, we have $\mu([\mathscr{I}_{d}])=\mu([i])=0$. Therefore, $\alpha$ belongs to kernel of the map $\mu\circ \eta$. Hence $\alpha$ can be written as a linear combination of $h-h^{-1}$ ($h\in \pi\setminus\{1\}$) and $1\in \pi$. 

The coefficient of $1\in \pi$ in $\alpha$ must be zero because the nonequivariant intersection number of both $\mathscr{I}_{d}(D^2)$ and $i(D^2)$ with $D_0$ are equal to $d$. This implies that $\alpha$ can be written as a linear combination of $h-h^{-1}$ ($h\in \pi\setminus\{1\}$). So it suffices to express $h-h^{-1}$ as a sum of terms of the form $\pm(g_{k,+}-g_{k,-}-g^{-1}_{k,+}+g^{-1}_{k,-})$. 

For this, we take the universal cover  $p: \widetilde{\Sigma}\to \Sigma$. Recall that $e_0=\{b_0\}$ denotes the zero-cell in the cellular decomposition of $\Sigma$. By fixing $\widetilde{b}_0'\in p^{-1}(b_0')$, we can identify elements in $\pi$ as fundamental domains of $\widetilde{\Sigma}$, and we can identify $g\in \pi_{1}(M_2;I_{j},S_0)$ as a lift $e^{1}_{j,g}$ of the $1$-cell $e^{1}_{j}$ in $\mathbb{H}$. Under these correspondences, $g_{\pm}$ are the two domains adjacent to $e^{1}_{j,g}$. Since $\widetilde{\Sigma}$ is connected, we may find a sequence $h=h_{n},\cdots ,h_{0}=1$ in $\pi'$ such that $\{h_{k},h_{k-1}\}=\{g_{k,+},g_{k,-}\}$ for some $g_{k}$. Then $h$ can be expressed as the sum of $\pm(g_{k,+}-g_{k,-})$, and hence $h-h^{-1}$ can be expressed as the sum of $\pm (g_{k,+}-g_{k,-}-g^{-1}_{k,+}+g^{-1}_{k,-})$.
\end{proof}

Let $\Emb_{H_0}(H_1, M_1)$ be the space of embeddings of $H_1$ in $M_1$ that extends to a smooth embedding of $H_0\cup H_1$ to $M$ by defining the extension on $H_0$ to be equal to $\mathscr{I}_0|_{H_0}$. Let the base point of $\Emb_{H_0}(H_1, M_1)$ be $\mathscr{I}_0|_{H_1}$. Then we have a canonical isomorphism
\[
\pi_1(F_1) \cong \pi_1 \Emb_{H_0}(H_1, M_1)\cong \pi_1\Emb_\partial'(\sqcup_{2\ell}I, M_1),
\]
where the fiber $F_1$ is defined in \eqref{eq: fibration tower}. 
The isotopy extension theorem gives an arc pushing map 
\[
\pi_1(F_1)\to \pi_0\Diff(M,\nu(H_0)\cup \mathscr{I}_0(H_1)),
\]
where $\Diff(M,\nu(H_0)\cup \mathscr{I}_0(H_1))$ denotes the group of diffeomorphisms of $M$ that are identity on $\nu(H_0)\cup \mathscr{I}_0(H_1)$. Note that since not all elements of $\pi_1(F_1)$ lift to loops of \emph{framed} arcs, the diffeomorphisms given by the arc pushing map may not be the identity on the normal bundle of $\mathscr{I}_0(H_1)$. 
\begin{Corollary}
\label{Cor_alpha_action}
Suppose $\alpha\in \pi_1\Emb_\partial'(\sqcup_{2\ell}I, M_1)$, and let 
\[
f_\alpha\in \pi_0\Diff(M,\nu(H_0)\cup \mathscr{I}_0(H_1))
\]
be its image under the arc pushing map. Then the difference
\[
\lambda_0([f_\alpha \circ i])-\lambda_0([i])\in \mathbb{Z}[ \pi]	
\] 
does not depend on the choice of $i\in \Emb^{hom}_{\partial}(D^2,M_2)$. Moreover, the map
\begin{align*}
\mathcal{P}: \pi_1\Emb_\partial'(\sqcup_{2\ell}I, M_1) &\to \mathbb{Z}[\pi] \\
\alpha &\mapsto \lambda_0([f_\alpha \circ i])-\lambda_0([i])
\end{align*}
is a homomorphism.
\end{Corollary}
\begin{proof}
By Lemma \ref{lem_lambda_0}, we may homotope $i$ so that its image only differs from $D_0$ in a neighborhood $N$ of $S_0$ and a finite union of arcs in $M_2$. By Lemma \ref{lem: pi1F1 generators}, the support of $f_\alpha$ can be taken to be disjoint from $N$. So $\lambda_0([f_\alpha \circ i])-\lambda_0([i])$ does not depend on $i$. 

By Lemma \ref{lem: pi1F1 generators}, we know that for 
\[
\alpha_1,\alpha_2\in \pi_1\Emb_{H_0}(H_1, M_1)\cong \pi_1\Emb_\partial'(\sqcup_{2\ell}I, M_1),
\]
we may always choose representatives of $f_{\alpha_1}, f_{\alpha_2}$ so that their supports are disjoint, therefore the map $\mathcal{P}$ is a homomorphism. 
\end{proof}

\begin{Definition}
\label{defn_subgroup_K}
Define $K\subset \pi_1\Emb_\partial'(\sqcup_{2\ell}I, M_1)$ to be the kernel of the homomorphism $\mathcal{P}$ in Corollary \ref{Cor_alpha_action}. 
\end{Definition}

\section{The $\Dax$ invariant}
\label{sec_dax_invariant}

Recall that the group $K\subset \pi_1(\Emb_\partial'(H_1, M_1))$ is defined in Definition \ref{defn_subgroup_K}.
Suppose $i:D^2\hookrightarrow M_2$ is a neatly embedded disk with $i(\partial D^2) = \mathscr{I}_0(\partial H_2)$, this section studies the properties of $\Dax(i,f_\alpha\circ i)\in \mathbb{Z}[\pi\setminus\{1\}]$, where $f_\alpha$ is the image of an element $\alpha\in K$ under the arc pushing map. Note that 
by Corollary \ref{Cor_alpha_action} and the definition of $K$, the embeddings $i, f_\alpha\circ i$ are homotopic to each other relative to the boundary if and only if $\alpha\in K$. Therefore, $\Dax( i, f_\alpha\circ i)$ is well-defined if and only if $\alpha\in K$.

We first observe that $\Dax(i,f_\alpha\circ i)$ does not depend on the choice of $i$.
\begin{Lemma}
\label{lem_f_alpha_Dax_equal_D0}
    Suppose $\alpha \in K$, let $f_\alpha \in \Diff(M,\nu(H_0)\cup H_1)$ be the arc pushing map of $\alpha$. Then for every neatly embedded disk $i:D^2\hookrightarrow M_2$ with $i(\partial D^2) = \mathscr{I}_0(\partial H_2)$, we have
\begin{equation}
\label{eqn_f_alpha_Dax_equal_D0}
\Dax( i, f_\alpha\circ i) = \Dax(D_0, f_\alpha(D_0)). 
\end{equation}
\end{Lemma}

\begin{proof}
    First, note that if $i_1\simeq i_2$ rel $\partial D^2$, then 
    \begin{align}
        & \Dax( i_1, f_\alpha\circ i_1) -  \Dax( i_2, f_\alpha\circ i_2) \nonumber \\
        = & \Dax( i_1, i_2) - \Dax(f_\alpha\circ i_1, f_\alpha\circ i_2) = 0,
        \label{eqn_inv_diff_f_alpha_homotopy}
    \end{align}
    where the second equality holds because $f_\alpha$ induces the identity map on $\pi$. Suppose the non-equivariant intersection number of $i$ and $D_0$ is $d$. Then, by Proposition \ref{prop_homotopy_after_barbell}, there exists $\alpha' \in \pi_1(\Emb_\partial'(H_1, M_1))$ such that $f_{\alpha'}\circ i$ is homotopic to $\mathscr{I}_d|_{H_2}$ rel $\partial D^2$. Since $\pi_1(\Emb_\partial'(H_1, M_1))$ is abelian (Lemma \ref{lem: pi1F1 generators}), we know that $f_{\alpha}$ is commutative with $f_{\alpha'}$ up to isotopy relative to $\nu(H_0)\cup H_1$. Hence we have
    \begin{align*}
         \Dax( i, f_\alpha\circ i) & = \Dax( f_{\alpha'}\circ i, f_{\alpha'}\circ f_\alpha\circ i) \\
        & = \Dax( f_{\alpha'}\circ i,  f_\alpha\circ f_{\alpha'}\circ i) 
        = \Dax( \mathscr{I}_d|_{H_2},  f_\alpha\circ \mathscr{I}_d|_{H_2}) ,
    \end{align*}
    where the first equality holds because $f_{\alpha'}$ induces the identity map on $\pi$, the second equality follows from the commutativity of $f_\alpha$ and $f_{\alpha'}$ up to isotopy, and the third equality follows from \eqref{eqn_inv_diff_f_alpha_homotopy}. 

    As a result, we only need to verify \eqref{eqn_f_alpha_Dax_equal_D0} when $i = \mathscr{I}_d|_{H_2}$. In this case, the image of $i$ only differs from $D_0$ in a small neighborhood $N$ of $S_0$. One may choose $f_\alpha$ so that its support is disjoint from $N$. It then follows directly from Definition \ref{def_scan_into_pi1Emb} that $\langle D_0, f_\alpha(D_0)\rangle$ can be decomposed into a finite product in $\pi_{1}\Emb_{\partial}(I, M_2)$ such that after changing the order of the products one obtains $\langle \mathscr{I}_d(H_2), f_\alpha\circ \mathscr{I}_d(H_2)\rangle$. Since $\pi_{1}\Emb_{\partial}(I, M_2)$ is abelian (Corollary \ref{cor_pi1Emb_abel}), we have  
    \[
\langle \mathscr{I}_d(H_2), f_\alpha\circ \mathscr{I}_d(H_2)\rangle =\langle D_0, f_\alpha(D_0)\rangle\in \pi_{1}\Emb_{\partial}(I, M_2), 
    \]
    and hence 
$
\Dax(\mathscr{I}_d|_{H_2}, f_\alpha\circ \mathscr{I}_d|_{H_2}) =\Dax( D_0, f_\alpha(D_0)).
$
\end{proof}

Recall that $\mathbb{Z}[\pi\setminus\{1\}]^\sigma$ denotes the subgroup of $\mathbb{Z}[\pi\setminus\{1\}]$ that is invariant under the involution defined by $g\mapsto g^{-1}$ on generators. This is the codomain of the Dax invariant for embedded disks. 
\begin{Definition}
    Define 
    \begin{align*}
        \mathcal{Q}: K &\to \mathbb{Z}[\pi\setminus\{1\}]^\sigma\\
                    \alpha&\mapsto \Dax( i, f_\alpha\circ i),
    \end{align*}
    where $i:D^2\hookrightarrow M_2$ is a neatly embedded disk with $i(\partial D^2) = \mathscr{I}_0(\partial H_2)$.
\end{Definition}
By Lemma \ref{lem_f_alpha_Dax_equal_D0}, the map $\mathcal{Q}$ does not depend on the choice of $i$. By the additivity property of the Dax invariant, the map $\mathcal{Q}$ is a homomorphism.

Recall that $\mathcal{C}$ denotes the set $(\pi\setminus\{1\})/\text{conjugation}$. There is a canonical projection 
\begin{equation}
\label{eqn_proj_pi_to_C}
\mathfrak{p}:\mathbb{Z}[\pi\setminus\{1\}]\to \mathbb{Z}[\mathcal{C}]
\end{equation}
that takes every generator of $\mathbb{Z}[\pi\setminus\{1\}]$ to its conjugation class. 

The main result of Section \ref{sec_dax_invariant} is the following proposition. 
\begin{Proposition}
\label{prop_arc_pushing_K_Dax}
The following composition map is zero:
\begin{equation}
\label{eqn_arc_pushing_K_Dax_compose}
K \xrightarrow{\mathcal{Q}}\mathbb{Z}[\pi\setminus\{1\}]^\sigma\hookrightarrow\mathbb{Z}[\pi\setminus\{1\}] \xrightarrow{\mathfrak{p}} \mathbb{Z}[\mathcal{C}]. 
\end{equation}
\end{Proposition}

\subsection{The Dax invariant for closed surfaces}
\label{subsec_Dax_closed_surface}
Before proving Proposition \ref{prop_arc_pushing_K_Dax}, we explain how to use this result to define the Dax invariant for pairs of closed surfaces in $M$.

We use the notation in the fibration tower \eqref{eq: fibration tower}. 
The Dax invariant is defined for pairs of embeddings $i_1,i_2\in   \Emb^{[\mathscr{I}_d]}_\bullet (\Sigma, M)$. Since $\Emb^{[\mathscr{I}_d]}_\bullet(H_0\cup H_{1}, M)$ is path connected, the map $\pi_0(F_2)\to \pi_0 \Emb^{[\mathscr{I}_d]}_\bullet (\Sigma, M)$ is surjective, so we may homotope $i_1,i_2$ so that they are both contained in $F_2$. By Proposition \ref{prop_homotopy_after_barbell}, there exist $\alpha_1,\alpha_2\in \pi_1(F_1)$, whose images under the arc pushing map are $f_{\alpha_1}$ and $f_{\alpha_2}$, such that $f_{\alpha_1}\circ i_1|_{H_2}$ and $f_{\alpha_2}\circ {i_2}|_{H_2}$ are homotopic disks relative to boundary in $M_2$. The embeddings $f_{\alpha_1}\circ i_1$ and $f_{\alpha_2}\circ {i_2}$ represent the images of $i_1,i_2$ in $\pi_0(F_2)$ under the action of $\alpha_1,\alpha_2\in \pi_1(F_1)$.

To distinguish notation, in this subsection, we use $\Dax_{M_2}$ to denote the relative Dax invariant for embedded disks in $M_2$, then the relative Dax invariant 
\[
\Dax_{M_2}(f_{\alpha_1}\circ i_1|_{H_2}, f_{\alpha_2}\circ i_2|_{H_2})\in \mathbb{Z}[\pi\setminus\{1\}]^\sigma
\]
is well-defined. 

\begin{Definition}
\label{def_Dax(i1,i2)}
    Define the relative Dax invariant for the pair of embedded closed surfaces $i_1,i_2$ by
    \[
\Dax(i_1,i_2) := \mathfrak{p}\circ \Dax_{M_2}(f_{\alpha_1}\circ i_1|_{M_2}, f_{\alpha_2}\circ i_2|_{M_2})\in \mathbb{Z}[\mathcal{C}]^\sigma.
    \]
\end{Definition}
Recall that $\mathbb{Z}[\mathcal{C}]^\sigma$ denotes the subgroup of $\mathbb{Z}[\mathcal{C}]$ generated by elements of the form $[g]+[g^{-1}]$. 

The elements $\alpha_1,\alpha_2\in \pi_1\Emb'_{\partial}(\sqcup_{2\ell}I,M_1)$ in Definition \ref{def_Dax(i1,i2)} are only determined by $i_1,i_2$ up to compositions with elements in $K$. By Proposition \ref{prop_arc_pushing_K_Dax}, different choices of $\alpha_1,\alpha_2$ do not change the value of $\mathfrak{p}\circ \Dax_{M_2}(f_{\alpha_1}\circ i_1|_{M_2}, f_{\alpha_2}\circ i_2|_{M_2})$. Since \eqref{eq: q_0} is a bijection, we conclude that $\Dax(i_1,i_2)$ in Definition \ref{def_Dax(i1,i2)} is well defined and only depends on the isotopy classes of $i_1,i_2$ relative to $b_0$. 

By Lemma \ref{lem_cS_cS'_bijection}, the Dax invariant extends to pairs of embedded surfaces in $M$ that are freely homotopic to $\mathscr{I}_d$. 

Now we establish the basic properties of the Dax invariant for closed surfaces. First, it follows immediately from the additivity of $\Dax_{M_2}$ that 
\[
\Dax(i_1,i_2) + \Dax(i_2,i_3) = \Dax(i_1,i_3).
\]

A priori, the definition of the Dax invariant may depend on the choice of the base point $b_0$ and the handle body decomposition of $\Sigma$.  We show that the Dax invariant is independent of these choices. First, note that the set $\mathcal{C}$ is canonically identified with the set of free homotopy classes of non-trivial loops on $\Sigma$, so it does not depend on the choice of the base point. 

\begin{Lemma}
\label{lem_Dax_indep_handle_decomp}
    Suppose $\Sigma = H_0'\cup H_1'\cup H_2'$ is another handle decomposition of $\Sigma$ with a single 0-handle, $2\ell$ 1-handles, and a single 2-handle, such that $b_0\in H_0'$. Let $\Dax'(i_1,i_2)$ be the Dax invariant defined with respect to $(H_0',H_1',H_2')$. Then 
$
\Dax'(i_1,i_2) = \Dax(i_1,i_2).
$
\end{Lemma}
\begin{proof}
The definition of $\Dax(i_1,i_2)$ can be alternatively described as follows. We choose isotopies of $i_1,i_2$ relative to $b_0$ such that their restrictions to $H_0\cup H_1$ are both equal to $\mathscr{I}_d|_{H_0\cup H_1}$ and the restrictions to $H_2$ are homotopic to $\mathscr{I}_d|_{H_2}$ relative to the boundary in $M_2$. Then define $\Dax(i_1,i_2)$ to be given by the Dax invariant of $i_1|_{H_1}$, $i_2|_{H_2}$ in $M_2$. Note that in this description, the Dax invariant only depends on the union $H_0\cup H_1$ and not on the full data of the handle decomposition. 
The handle decomposition is only used when we compute the fundamental group of the embedding space of $H_0\cup H_1$ in $M$ using the fibration tower \eqref{eq: fibration tower}. 

Since there is an isotopy of $\Sigma$ relative to $b_0$ that takes $H_0\cup H_1$ to $H_0'\cup H_1'$, the result is proved. 
\end{proof}

Now we consider the dependency of $\Dax(i_1,i_2)$ on the base point. 
\begin{Lemma}
\label{lem_Dax_indep_base_point}
    The value of $\Dax(i_1,i_2)$ does not depend on the choice of $b_0$.
\end{Lemma}
\begin{proof}
    Suppose $b_0^{(1)}$, $b_0^{(2)}$ are two choices of base points. We may choose a handle decomposition such that both points are contained in $H_0$, so the desired result follows from Lemma \ref{lem_Dax_indep_handle_decomp}. 
\end{proof}

Since elements in $\mathcal{C}$ are represented by the free homotopy classes of non-trivial loops on $\Sigma$, every element in $\Diff(\Sigma)$ acts on $\mathbb{Z}[\mathcal{C}]$. For $\varphi\in \Diff(\Sigma)$, let $\varphi_*$ denote its action on $\mathbb{Z}[\mathcal{C}]$. The construction of the Dax invariant and Lemmas \ref{lem_Dax_indep_handle_decomp}, \ref{lem_Dax_indep_base_point} immediately imply the following result. Here, we regard the Dax invariant as defined for pairs of embedded, non-parametrized surfaces as in Remark \ref{rmk_relative_Dax_surfaces}.

\begin{Lemma}
\label{lem_Dax_equiv_diffeo}
    Suppose $\varphi:\Sigma\to \Sigma$ is an orientation-preserving diffeomorphism, and let $\Sigma_1,\Sigma_2\subset M$ be embedded surfaces that are both homotopic to $\Sigma_{(d)}$. Then 
    \[
\Dax((\varphi\times\id)(\Sigma_1), (\varphi\times\id)(\Sigma_2)) = \varphi_*\Dax(\Sigma_1,\Sigma_2).
\makeatletter\displaymath@qed\makeatother
    \]
\end{Lemma}

Suppose $f:(M,b)\to (M,b)$ is a diffeomorphism that induces the identity map on $\pi_1$. Let $i_1,i_2$ be as in Definition \ref{def_Dax(i1,i2)}. Then by Lemma \ref{lem_hom_in_F2}, there exists $d$ such that $f\circ i_1$ and $f\circ i_2$ are both homotopic to $\mathscr{I}_d$ relative to $b_0$.

\begin{Lemma}
\label{lem_Dax_change_by_MCG0_*}
    Suppose $f:(M,b)\to (M,b)$ is a diffeomorphism that induces the identity map on $\pi_1$.  Then
$
\Dax(i_1,i_2) = \Dax(f\circ i_1, f\circ i_2).
$
\end{Lemma}

\begin{proof}
    After isotopy of $f$ relative to $b$, we may assume that $f$ is the identity map on $\mathscr{I}_0(H_0\cup H_1)$ and takes $M_2$ to $M_2$. Since $\pi_1(M_2)\to \pi_1(M)$ is an isomorphism, we know that $f$ induces the identity map on $\pi_1(M_2)$. After isotopies of $i_1,i_2$, we may assume that $i_1|_{H_0\cup H_1} = i_2|_{H_0\cup H_1} = \mathscr{I}_d|_{H_0\cup H_1}$ and $i_1|_{H_2}, i_2|_{H_2}, \mathscr{I}_d|_{H_2}$ are neatly embedded disks in $M_2$ homotopic to each other relative to the boundary. Then 
    \begin{align*}
\Dax (f\circ i_1, f\circ i_2) &= \mathfrak{p}\circ \Dax_{M_2}(f\circ i_1|_{H_2},f\circ i_2|_{H_2}) \\
&= \mathfrak{p}\circ \Dax_{M_2}( i_1|_{H_2}, i_2|_{H_2}) = \Dax(i_1,i_2).
\qedhere
    \end{align*}
\end{proof}

 The rest of the section is devoted to the proof of Proposition \ref{prop_arc_pushing_K_Dax}. We remark that arguments in this proof will not be used later in the paper.

\subsection{A splitting of $\pi_{1}(\Emb'_{\partial}(\sqcup_{2\ell}I,M_1))$}
We introduce a canonical splitting of $\pi_{1}(\Emb'_{\partial}(\sqcup_{2\ell}I,M_1))$ as the direct sum of $\pi_{1}(\Emb_{\partial}(\sqcup_{2\ell}I,M_1))$ and $\mathbb{Z}^{2\ell}$. This splitting will
allow us to reduce the verification of Proposition \ref{prop_arc_pushing_K_Dax} to the subgroup of $K$ generated by elements of the form $\hat \rho_\gamma$.

Recall that $e^1_j$ ($j=1,\dots,2\ell$) denote the 1-cells of $\Sigma$ and $I_j$ denotes $\mathscr{I}_0(e_j^1)\cap M_1$. 
For each $\alpha\in \pi_{1}(\Emb'_{\partial}(\sqcup_{2\ell}I,M_1))$, let $f_\alpha$ be the image of $\alpha$ under the arc pushing map, recall that $f_\alpha$ fixes $\nu(H_0)\cup \mathscr{I}_0(H_1)$. Consider the action of $f_\alpha$ on the normal bundle of each arc $I_j$. It is a bundle map that equals the identity on the boundary and also on the sub-bundle tangent to $H_1$. Therefore, the action of $f_\alpha$ on the normal bundle of $I_j$ defines an element in $\pi_1(\SL(2))\cong \mathbb{Z}$. We call the corresponding integer the \emph{rotation number} of $f_\alpha$ on the normal bundle of $I_j$.

\begin{Definition}
Let 
$\mathcal{R}:\pi_1(\Emb_\partial'(H_1, M_1))\to \bZ^{2\ell}$ be the homomorphism such that the $j^{th}$ coordinate  of $\mathcal{R}(\alpha)$ equals the rotation number of $f_\alpha$ on the normal bundle of $I_j$.
\end{Definition}
It is clear from the definition that $\mathcal{R}$ is a group homomorphism. 

In the following lemma, recall that $\hat \rho_\gamma,\hat\xi_j\in \pi_1\Emb'(\sqcup_{2\ell} I,M_1)$ are defined in Definition \ref{defi: spinning around S0} and Definition \ref{defn_hat_xi_j}. 
\begin{Lemma}
\label{lem_computation_rotation_number}
The map $\mathcal{R}$ satisfies the following properties:
\begin{enumerate}
\item  $\mathcal{R}(\hat \rho_\gamma)=0$.
\item   $\mathcal{R}(\hat \xi_j)=(0,\dots,0,2,0,\dots,0)$, where the coefficient $2$ is on the 
 $j$-th coordinate.
\end{enumerate}
\end{Lemma}
\begin{proof}
By definition, $\hat \rho_\gamma$ lifts to a loop of \emph{framed} embeddings, so $f_{\hat \rho_\gamma}$ acts trivially on the normal bundle of  $\mathscr{I}_0(e^1_j)$, and hence $\mathcal{R}(\hat \rho_\gamma)=0$.

 For Part (2), we need to lift $\hat \xi_j$ to a path of framed embeddings 
 of $\sqcup_{2\ell} I$ and compute the framing on the ending point of the path. For $u\neq j$, we just lift it to the constant path of framings on 
 $I_u$. For $I_j$, consider the following lifting problem
\[\xymatrix{  ([0,1]\times \partial I_j)\cup(\{0\}\times I_j) 
\ar[r]^-{c}\ar@{^{(}->}[d] & \operatorname{SO(3)}\ar[d]^{\operatorname{ev}_{w}}\\
 [0,1]\times I_j\ar@{.>}[ru]^{\widetilde{h}}\ar[r]^{h} & S^2.  
} 
\]
Here $c$ is the contant map with value $\id$, $\operatorname{ev}_{w}$ is defined by the action of $SO(3)$ on 
$w=(0,0,1)\in S^2$, and $h$ is a map that equals $w$ on the boundary and has degree $1$. Since 
$SO(3)\to S^2$ is a principle $SO(2)$-bundle of Euler number $2$, the map
$$
\widetilde{h}|_{\{1\}\times I_j} : \{1\}\times I_j \to \operatorname{ev}_{w}^{-1}(w)=SO(2)
$$
is a degree--$2$ map (the boundary of its domain is mapped to $\id$ by construction). 
Since $\widetilde{h}$ gives the lifting of $\hat\xi_j$ to a path of \emph{framed} embeddings for $I_j$, we conclude that 
 $\mathcal{R}(\hat \xi_j)=(0,\cdots,0,2,0,\cdots,0)$
\end{proof}

Note that we have a fibration
\[ 
 \prod_{j=1}^{2l}\Omega S^2\to \Emb'_{\partial}(\sqcup_{2\ell}I,M_1)\to \Emb_{\partial}(\sqcup_{2\ell}I,M_1).
\]
Hence, there is an exact sequence
\begin{equation}
\label{eqn_exact_seq_pi_1Emb'}
\mathbb{Z}^{2\ell} \to \pi_1\Emb'_{\partial}(\sqcup_{2\ell}I,M_1)\to \pi_1\Emb_{\partial}(\sqcup_{2\ell}I,M_1)\to 0,
\end{equation}
where the image of the $j^{th}$ basis elements of $\mathbb{Z}^{2\ell}$ is $\hat\xi_j$.
The next corollary shows that \eqref{eqn_exact_seq_pi_1Emb'} is a split short exact sequence
\begin{Corollary}
\label{cor_split_pi_1Emb'}
\begin{enumerate}
\item Each element $\bar{\alpha}\in  \pi_1\Emb_{\partial}(\sqcup_{2\ell}I,M_1)$ has a unique lifting  $ \alpha\in  \pi_1\Emb_{\partial}'(\sqcup_{2\ell}I,M_1)$ such that $\mathcal{R}(\alpha)=0$. 
\item The elements $\hat\xi_j$ ($j=1,\dots,2\ell$) are linearly independent.
\item We have
\begin{equation}
\label{eqn_split_pi1_Emb'}
\pi_1\Emb_{\partial}'(\sqcup_{2\ell}I,M_1) = \ker\mathcal{R} \oplus \mathbb{Z}^{2\ell},
\end{equation}
where the second summand is the free abelian group generated by $\hat\xi_j$ ($j=1,\dots,2\ell$). 
\end{enumerate}
\end{Corollary}

\begin{proof}
 By Lemmas \ref{lem_computation_rotation_number} and \ref{lem: pi1F1 generators}, the image of $\mathcal{R}$ is always even, therefore $\frac12 \mathcal{R}$ is a well-defined homomorphism from $\pi_1\Emb_{\partial}'(\sqcup_{2\ell}I,M_1)$ to $\mathbb{Z}^{2\ell}$. 

    Let $t: \mathbb{Z}^{2\ell} \to \pi_1\Emb'_{\partial}(\sqcup_{2\ell}I,M_1)$ be the map in \eqref{eqn_exact_seq_pi_1Emb'}, which takes the $j^{th}$ basis element to $\hat\xi_j$. By Lemma \ref{lem_computation_rotation_number}, the composition $(\frac12 \mathcal{R})\circ t$ is the identity map, so \eqref{eqn_exact_seq_pi_1Emb'} is a split short exact sequence, and there is a unique lifting:
    \[
    l: \pi_1\Emb_{\partial}(\sqcup_{2\ell}I,M_1) \to \pi_1\Emb_{\partial}'(\sqcup_{2\ell}I,M_1)
    \]
    such that $(\frac12 \mathcal{R}) \circ l = 0$. 
\end{proof}

\begin{Corollary}
\label{cor_kerR_gen_by_rho}
    $\ker \mathcal{R}$ is generated by elements of the form $\hat \rho_{\gamma}$.
\end{Corollary}

\begin{proof}
    By Lemma \ref{lem_computation_rotation_number}, we have $\hat \rho_{\gamma}\in \ker \mathcal{R}$ for all $\gamma$. By Lemma \ref{lem: pi1_emb_H_1_generators}, the images of $\hat \rho_{\gamma}$ in $\pi_1\Emb_{\partial}(\sqcup_{2\ell}I,M_1)$ form a generating subset. Therefore, by Corollary \ref{cor_split_pi_1Emb'}, we conclude that $\{\hat \rho_{\gamma}\}$ generates $\ker \mathcal{R}$.
\end{proof}

Recall from Definition \ref{defn_subgroup_K} that $K$ is the subgroup of $\pi_1\Emb_{\partial}'(\sqcup_{2\ell}I,M_1)$ consisting of elements $\alpha$ whose images $f_\alpha$ under the arc pushing map do not change the homotopy class of $D_0$. 
Let $K_0 = K \cap \ker(\mathcal{R})$. The next lemma shows that we only need to verify Proposition \ref{prop_arc_pushing_K_Dax} when $\alpha \in K_0$. 

\begin{Lemma}
For every $\alpha\in K$, there exists an element $\alpha'\in K$ such that 
\begin{enumerate}
	\item  $\alpha-\alpha'\in K_0$,
	\item $\mathcal{Q}(\alpha')=0$. 
\end{enumerate}
\end{Lemma}

\begin{proof}
Suppose $\mathcal{R}(\alpha)=(2a_1,\cdots,2a_{2\ell})$. We have
\[
\langle \Sigma, SO(2)\rangle \cong H^1(\Sigma;\bZ)\cong H^1(H_0\cup H_1;\bZ) \cong \bZ^{2\ell}.
\]
Pick a map $f:\Sigma\to SO(2)$ corresponding to $(2a_1,\cdots,2a_{2\ell})\in \bZ^{2\ell}$ in the above
isomorphism. We may further assume $f(H_0)=\{\id\}$. We denote the composition
\[
 \Sigma \xrightarrow{f}	SO(2) \hookrightarrow SO(3)
\]
by $h_1$, which is null-homotopic. Let $h_t:\Sigma\to SO(3)$ $(t\in [0,1])$ be a hull-homotopy such that
$h_0(\Sigma)=\{\id\}$ and $h_t(H_0)=\{\id\}$ for all $t\in [0,1]$. Define 
the following isotopy
\[
r_t: \Sigma\times S^2 \to \Sigma\times S^2,\quad
      (x,v)\mapsto (x,h_t(x)v),
\]
where we assume the coordinate for $b_1\in S^2$ is $(0,0,1)$ such that it is fixed under the action 
of $SO(2)\hookrightarrow SO(3)$.

Let $\alpha'$ be the loop of embeddings defined by the images of the base embedding under the action of $\{r_t\}_{t\in[0,1]}$. Then by definition, the map $f_{\alpha'}$ is represented by $r_1$. Therefore, we have $\mathcal{R}(\alpha')=(2a_1,\cdots,2a_{2\ell})=\mathcal{R}(\alpha)$. We also have $f_{\alpha'}|_\Sigma = \id$, which implies that $\alpha' \in K$ and $\Dax(D_0,f_{\alpha'}(D_0)) = 0$. So the desired result is proved. 
\end{proof}

\begin{Corollary}
\label{cor_K_image_equal_K0}
    The image of $\mathcal{Q}:K\to \mathbb{Z}[\pi\setminus\{1\}]^\sigma$ is equal to the image of its restriction to $K_0$. \qed
\end{Corollary}

By Corollary \ref{cor_K_image_equal_K0}, in order to prove Proposition \ref{prop_arc_pushing_K_Dax}, we only need to verify that the image of $K_0$ under the composition map \eqref{eqn_arc_pushing_K_Dax_compose} is zero. 

\subsection{Calculus of barbells}
This subsection studies the actions of certain barbell diffeomorphisms on the disk $D_0$. 
From now on, we fix a choice of orientations on $\Sigma$ and $S^2$.

\subsubsection{Self-referential barbell diffeomorphisms}
\label{subsubsec_self-referential_barbell}
We define a class of barbell diffeomorphisms on $M_2$, which will be called the \emph{self-referential barbell diffeomorphisms}.

Recall that $H_2$ denotes the 2-handle of $\Sigma$. Let $p,q$ be distinct points in the interior of $H_2$, let $\gamma$ be a smooth, not necessarily embedded arc in $\Sigma$ from $p$ to $q$ such that 
\begin{enumerate}
    \item there is a unique interior point of $\gamma$ that passes through $q$, and the endpoint $q$ intersects the interior of $\gamma$ transversely, 
    \item the interior of $\gamma$ is disjoint from $p$. 
\end{enumerate}
 Let $\hat \gamma$ be a lifting of $\gamma$ to $M_2$ such that the only self-intersection is the endpoint that lifts $q$. 
 Let $m$ be a meridian of the interior of $\hat \gamma$ at $q$. 
 Orient $m$ so that the linking of $m$ with the interior of $\hat \gamma$ is positive. Modify $\hat \gamma$ near the endpoint $q$ so that it becomes an embedded arc from $\{p\}\times S^2$ to $m$. Consider the barbell diffeomorphism associated with the triple $(\{p\}\times S^2, m, \hat \gamma)$, where the orientation of $\{p\}\times S^2$ is given by the chosen orientation of the $S^2$--factor of $M$. 
We define $\beta_{p,\gamma}$ to be the above diffeomorphism extended to $M_2$ by the identity, and call it a \emph{self-referential barbell diffeomorphism}.

\begin{remark}
\label{rmk_barbell_commutative_sum_notation}
Note that $\beta_{p,\gamma}$ is only determined by $p,\gamma$ up to isotopy. In the following, we will often regard $\beta_{p,\gamma}$ as an element in $\pi_0\Diff_\partial(M_2)$. Up to isotopy, the compositions of self-referential barbell diffeomorphisms are commutative, because each pair of self-referential barbell diffeomorphisms can be isotoped to have disjoint support. As a result, we will use ``$+$'' and ``$-$'' to denote the composition and quotient of self-referential barbell diffeomorphisms. 
\end{remark}

 Suppose $\gamma$ is as above, let $c(\gamma)$ denote the loop that goes from $q$ to $q$ along $\gamma$, where the orientation agrees with $\gamma$. Then $c(\gamma)$ defines a conjugation class in $\pi$. Let $[c(\gamma)]\in \mathbb{Z}[\mathcal{C}]$ be the conjugation class of $c(\gamma)$ if it is non-trivial, and be zero if it is trivial. Then we have the following result.
 \begin{Lemma}
     \label{lem_self_referential_barbell_Dax}
     Let $\beta_{p,\gamma}$ be a self-referential diffeomorphism. Then
     \begin{enumerate}
         \item $\beta_{p,\gamma}$ is pseudo-isotopic to the identity relative to $\partial M_2$.
         \item  We have
         \[
\mathfrak{p}\circ \Dax(\beta_{p,\gamma}(D_0), D_0) = \epsilon([c(\gamma)] + [c(\gamma)^{-1}]),
     \]
     where $\epsilon\in\{1,-1\}$ is a universal constant depending only on the orientation conventions. 
     \end{enumerate} 
   
 \end{Lemma}

 \begin{proof}
     The fact that $\beta_{p,\gamma}$ is pseudo-isotopic to the identity relative to $\partial M_2$ follows from \cite[Proposition 2.6]{BG2023}. The image of $D_0$ under is obtained by attaching a self-referential tube defined by $\hat \gamma$ to $D_0$, so its relative Dax invariant is given by $\pm(c(\gamma)+c(\gamma)^{-1})$ up to conjugation, where the sign is a universal constant. So the desired result is proved. We refer the reader to \cite{gabai2021self,schwartz20214} for the definition of self-referential tubes and the computation of their Dax invariants.
 \end{proof}
 
 \begin{Definition}
 \label{defn_associated_loop_self_referential_barbell}
     We will call the $c(\gamma)$ above the \emph{associated loop} with the self-referential barbell diffeomorphism $\beta_{p,\gamma}$. 
 \end{Definition}

The following lemma will be useful when studying the effect of barbell diffeomorphisms on Dax invariants.
\begin{Lemma}
\label{lem_sum_beta_Dax_rewrite}
    Suppose $\beta_1,\dots,\beta_m$ is a sequence of elements in $\pi_0\Diff_\partial(M_2)$ that are commutative to each other and are homotopic to the identity map on $M_2$ relative to $\partial M_2$. Then
    \[
\Dax((\sum_{j=1}^m\beta_j)(D_0),D_0) = \sum_{j=1}^m \Dax(\beta_j(D_0),D_0). 
    \]
\end{Lemma}

Recall that the summation notation on the left-hand side denotes the composition map. 

\begin{proof}
    Since each $\beta_j$ is homotopic to the identity, it induces the identity map on $\pi_1(M_2)$. As a result,
$\Dax(\beta_j(D_1),\beta_j(D_2)) = \Dax(D_1,D_2)$
    for all embedded disks $D_1,D_2$ that are homotopic relative to the boundary. Therefore,
    \begin{align*}
        \Dax((\sum_{j=1}^m\beta_j)(D_0),D_0) 
       = & \sum_{l=1}^m \Dax  \Big((\sum_{j=1}^l\beta_j)(D_0),(\sum_{j=1}^{l-1}\beta_j)(D_0)\Big)\\
       = & \sum_{l=1}^m \Dax (\beta_l(D_0),D_0). \qedhere
    \end{align*}
\end{proof}

\subsubsection{Vertical barbell diffeomorphisms}
Now we define another class of barbell diffeomorphisms, which will be called \emph{vertical barbell diffeomorphisms}. 

Let $p_0,p_1$ be two distinct points in $\inte(H_2)$, and let $\gamma$ be a smooth, not necessarily embedded, arc in $\Sigma$ connecting $p_1,p_2$ such that $\inte(\gamma)$ is disjoint from $p_0,p_1$. We define a barbell diffeomorphism on $M_2$ associated with $p_0,p_1,\gamma$ as follows. 
Lift $\gamma$ to an embedded arc $\hat{\gamma}$ in $M_2$ connecting $\{p_0\}\times S^2$ and $\{p_1\}\times S^2$. Note that any two such liftings are isotopic to each other. Consider the barbell diffeomorphism associated with the triple $(\{p_0\}\times S^2,\{p_1\}\times S^2,\hat{\gamma})$, and define $\beta_{p_0,p_1,\gamma}$ to be its extension to $M_2$ by the identity. We call $\beta_{p_0,p_1,\gamma}$ a \emph{vertical barbell diffeomorphism}. Note that Remark \ref{rmk_barbell_commutative_sum_notation} also applies to vertical barbells.

We study the behavior of vertical barbell diffeomorphisms under a homotopy of the defining data. Suppose $(p_0(t),p_1(t),\gamma(t))$ ($t\in[0,1]$) is a smooth 1-parameter family, such that for each $t$, the points $p_0(t), p_1(t)$ are distinct and contained in $\inte(H_2)$, and $\gamma(t)$ is a smooth, not necessarily embedded, arc in $\Sigma$ connecting $p_0,p_1$. Assume that when $t=0,1$, the interior of $\gamma(t)$ is disjoint from $\{p_0(t), p_1(t)\}$. We also assume that each $\gamma(t)$ is parametrized by $[0,1]$, and the following transversality properties hold:
\begin{enumerate}
    \item There are only finitely many $t$ such that the interior of $\gamma(t)$ passes through $p_0(t)$ or $p_1(t)$.
    \item  For each $t$ such that the interior of $\gamma(t)$ intersects $\{p_0(t), p_1(t)\}$, there is exactly one point $x\in(0,1)$ such that $\gamma(t)(x) \in \{p_0(t), p_1(t)\}$ and the self-intersection at $\gamma(t)(x)$ is transverse.
    \item The intersection of the interior of $\gamma(t)$ with $\{p_0(t), p_1(t)\}$ are also transverse with respect to the parameter $t$, in the following sense. Suppose $\gamma(t)(x) = p_n(t)$ where $n=0$ or $1$. Then $\frac{\partial}{\partial x} \gamma(t)(x)$ and $\frac{\partial}{\partial t} p_n(t)-\frac{\partial}{\partial t}\gamma(t)(x)$ are linearly independent.
\end{enumerate}
Every 1-parameter family $(p_0(t),p_1(t),\gamma(t))$ can be perturbed to satisfy the above conditions.

Suppose the interior of $\gamma(t)$ passes through $p_0(t)$ or $p_1(t)$ when $t$ is equal to $t_1<\dots<t_m$. For each $t_j$, let $n_j\in\{0,1\}$ be the coordinate such that the interior of $\gamma(t_j)$ passes through the end point $\gamma(t_j)(n_j)=p_{n_j}$, let $x_j\in (0,1)$ be such that $\gamma(t_j)(x_j) = \gamma(t_j)(n_j)$.

For each $j$, define a self-referential barbell $\beta_j$ associated with $\gamma(t_j)$ as follows. If $n_j=1$, let $\beta_j = \beta_{p_0(t_j),\gamma(t_j)}$. If $n_j=0$, let $\bar\gamma(t_j)$ be the same arc as $\gamma(t_j)$ but with the revered orientation, and let $\beta_j = \beta_{p_1(t_j),\bar\gamma(t_j)}$.

Let $c_j\in \pi$ be the loop from $\gamma(t_j)(x_j)$ to itself along $\gamma(t_j)$, endowed with the same orientation as $\gamma(t_j)$. If $n_j=1$, this is the loop associated with the self-referential barbell diffeomorphism $\beta_j$ (see Definition \ref{defn_associated_loop_self_referential_barbell}); if $n_j=0$, this is the reversal of the loop associated with $\beta_j$.  Let $\epsilon_j\in\{1,-1\}$ be the orientation of 
\begin{equation}
\label{eqn_sign_vertical_barbell_crossing}
\frac{\partial}{\partial x} \gamma(t)(x),\quad \frac{\partial}{\partial t} p_{n_j}(t)-\frac{\partial}{\partial t}\gamma(t)(x)
\end{equation}
at $(t_j,x_j)$ as a basis of the tangent space of $\Sigma$, relative to the chosen orientation of $\Sigma$. 

The next result compares the actions of $\beta_{p_0(t),p_1(t),\gamma(t)}$ for $t=0,1$. 
\begin{Proposition}
\label{prop_vertical_barbell_crossing}
 Let $(p_0(t),p_1(t),\gamma(t))$, $\beta_j$, $\epsilon_j$ and $c_j$ be as above. Then 
 \begin{enumerate}
     \item There exists a universal constant $\epsilon^{(1)}\in \{1,-1\}$ such that 
     \[
     \beta_{p_0(1),p_1(1),\gamma(1)} - \beta_{p_0(0),p_1(0),\gamma(0)} = \epsilon^{(1)} \sum_{j=1}^m \epsilon_j\beta_j.
     \]
     
     \item The vertical barbell diffeomorphisms $\beta_{p_1(0),p_2(0),\gamma(0)}$ and $\beta_{p_1(1),p_2(1),\gamma(1)}$ are pseudo-isotopic relative to $\partial M_2$.

     \item There exists a universal constant $\epsilon^{(2)}\in \{1,-1\}$ such that  \begin{align}
& \mathfrak{p}\circ \Dax\big(\beta_{p_0(1),p_1(1),\gamma(1)}(D_0), \beta_{p_0(0),p_1(0),\gamma(0)}(D_0)\big)\nonumber \\
= & \epsilon^{(2)} \sum_{j=1}^m  \epsilon_j([c_j]+[c_j^{-1}])\in \mathbb{Z}[\mathcal{C}]. \label{eqn_vertical_barbell_crossing_RHS}
 \end{align}
 Here, $[c_j]$ denotes the conjugation class of $c_j$ if it is non-trivial, and denotes zero if it is trivial. 
 \end{enumerate}
 Moreover, we have $\epsilon^{(2)} = \epsilon^{(1)}\epsilon$, where $\epsilon$ is the constant in Lemma \ref{lem_sum_beta_Dax_rewrite}. 
\end{Proposition}
The constants $\epsilon^{(1)},\epsilon^{(2)}$ depend only on the orientation conventions in the setup of the notation. 

\begin{proof}
    It is clear that for each $j$, the barbell diffeomorphisms $\beta_{p_0(t),p_1(t),\gamma(t)}$ for $t\in(t_j,t_{j+1})$ are isotopic in $\Diff_\partial(M_2)$. Hence, we only need to compute the change of the barbell diffeomorphism across $t_1,\dots,t_m$. 
    By Proposition \ref{prop: connected sum}, we know that when $t$ goes across $t_j$, the map $\beta_{p_0(t),p_1(t),\gamma(t)}$ changes by a composition with $\pm \beta_j$. A straightforward tracking of the signs proves Part (1) of the desired proposition.

    Part (2) follows from the result of Part (1) and Lemma \ref{lem_self_referential_barbell_Dax} (1). 

    To prove Part (3), note that since the vertical barbell diffeomorphisms induce the identity map on $\pi_1(M_2)$, we have
    \begin{align*}
& \Dax\big(\beta_{p_1(1),p_2(1),\gamma(1)}(D_0), \beta_{p_1(0),p_2(0),\gamma(0)}(D_0)\big) 
\\
= & \Dax\big((\beta_{p_1(1),p_2(1),\gamma(1)}-\beta_{p_1(0),p_2(0),\gamma(0)})(D_0), D_0\big). 
    \end{align*}
    By Lemma \ref{lem_sum_beta_Dax_rewrite} and the result of Part (1), we have 
    \begin{align*}
\Dax\big((\beta_{p_1(1),p_2(1),\gamma(1)}-\beta_{p_1(0),p_2(0),\gamma(0)})(D_0), D_0\big)  
= \epsilon^{(1)} \sum_{j=1}^m \Dax(\epsilon_j\beta_j(D_0),D_0).
    \end{align*}
    By Lemma \ref{lem_self_referential_barbell_Dax} (2), 
    \[
\mathfrak{p}\circ \Dax(\epsilon_j\beta_j(D_0),D_0) = \epsilon_j\,\mathfrak{p}\circ \Dax(\beta_j(D_0),D_0) = \epsilon_j\,\epsilon\, [c_j],
    \]
    where $\epsilon\in\{1,-1\}$ is the constant in Lemma \ref{lem_self_referential_barbell_Dax}. So Part (3) is proved. 
\end{proof}

\subsubsection{Meridian-vertical barbells}
\label{subsubsec_mer_vert_barbell}
\begin{figure}
	\begin{overpic}[width=0.2\textwidth]{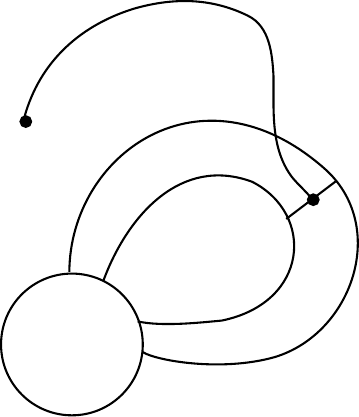}
		\put(3,10){\small$H_0$}
		\put(55,15){\small$h$}
            \put(-3,60){\small$p$}
            \put(62,36){\small$q$}
            \put(70,46){\small$c$}
            \put(55,75){\small$\gamma$}
	\end{overpic}
	\caption{Defining data of a meridian-vertical barbell}
    \label{fig_meridian-vertical-barbell}
\end{figure}

Now we define another class of barbell diffeomorphisms on $M_2$, which will be called the \emph{meridian-vertical barbells}. Let $h\subset H_1$ be a 1-handle of $\Sigma$, let the arc $c$ be a co-core of $h$, and let $q$ be the intersection of $c$ with the core of $h$. Let $p$ be a point in $\inte(H_2)$, let $\gamma$ be a smooth, not necessarily embedded, arc in $\Sigma$ connecting $p,q$ such that $\inte(\gamma)$ is disjoint from $p$ and $c$, and $\gamma$ intersects $c$ at $q$ transversely. See Figure \ref{fig_meridian-vertical-barbell}.

We define a barbell diffeomorphism on $M_2$ associated with $p,q,\gamma$ as follows. Let $m$ be the meridian of $\mathscr{I}_0(h)$ near $\mathscr{I}_0(c)$, let $B$ be a 3--ball bounded by $m$ that intersects $\mathscr{I}_0(h)$ transversely at $\mathscr{I}_0(c)$.  Lift $\gamma$ to an arc $\hat{\gamma}$ that connects $m$ and $\{p\}\times S^2$ such that the interior of $\hat\gamma$ is disjoint from $B$ and $\{p\}\times S^2$. Let $\mathfrak{o}$ be an orientation $m$. Consider the barbell diffeomorphism associated with the triple $(\{p\}\times S^2,m,\hat{\gamma})$, where $m$ is oriented by $\mathfrak{o}$ and $\{p\}\times S^2$ is oriented by the orientation of the $S^2$--factor of $M$. Define $\beta_{p,q,\gamma,\mathfrak{o}}$ to be the extension of the barbell diffeomorphism to $M_2$ by the identity. We call $\beta_{p,q,\gamma,\mathfrak{o}}$ \emph{meridian-vertical barbell diffeomorphism}.

Each meridian-vertical barbell has two associated vertical barbells, which are defined as follows. Extend the co-core $c$ a little bit to $c^\dagger$ such that the endpoints of $c^\dagger$ are in the interior of $H_2$. Let $q',q''$ be the two endpoints of $c^\dagger$. Let $\gamma'$ be the extension of $\gamma$ along $c$ so that the endpoint $q$ is moved to $q'$, and define $\gamma''$ similarly. Consider the vertical barbells diffeomorphisms $\beta_{p,q',\gamma'}$, $\beta_{p,q'',\gamma''}$. By Proposition \ref{prop: connected sum}, we have
\begin{equation}
\label{eqn_meridian_decomp_two_vertical}
   \beta_{p,q',\gamma'}  - \beta_{p,q'',\gamma''} = \pm \beta_{p,q,\gamma,\mathfrak{o}},
\end{equation}
where the sign on the right-hand side depends only on $\mathfrak{o}$, the ordering of $(\gamma',\gamma'')$, and the universal orientation conventions, and does not depend on the choices of $p,\gamma$.

\begin{Proposition}
\label{prop_meridian_vertical_homot}
    Let $p,q$ be as above, suppose $\gamma_1$, $\gamma_2$ that are smooth, not necessarily embedded, arcs in $\Sigma$ connecting $p,q$ such that $\inte(\gamma_j)$ $(j=1,2)$ are disjoint from $p$ and $c$. Suppose $\gamma_1$ and $\gamma_2$ are homotopic relative to the endpoints. Then
    \begin{enumerate}
        \item $\beta_{p,q,\gamma_1,\mathfrak{o}}$ and $\beta_{p,q,\gamma_2,\mathfrak{o}}$ are pseudo-isotopic relative to $\partial M_2$.
        \item $\mathfrak{p}\circ \Dax\big(\beta_{p,q,\gamma_1,\mathfrak{o}}(D_0), \beta_{p,q,\gamma_2,\mathfrak{o}}(D_0)\big)=0.$
    \end{enumerate}
    
\end{Proposition}

\begin{proof}
    Let $q',q''$ and $\gamma_j',\gamma_j''$ $(j=1,2)$ be defined as above. Without loss of generality, we may assume the orientation $\mathfrak{o}$ is chosen so that \eqref{eqn_meridian_decomp_two_vertical} gives
\[
   \beta_{p,q',\gamma_j'}  - \beta_{p,q'',\gamma_j''} = \beta_{p,q,\gamma_j,\mathfrak{o}}
\]
for $j=1,2$. 
By Proposition \ref{prop_vertical_barbell_crossing}, we know that $(\beta_{p,q',\gamma_1'},\beta_{p,q',\gamma_2'})$ and $(\beta_{p,q'',\gamma_1''},\beta_{p,q'',\gamma_2''})$ are pairwise pseudo-isotopic relative to $\partial M_2$, so $\beta_{p,q,\gamma_1,\mathfrak{o}}$ and $\beta_{p,q,\gamma_2,\mathfrak{o}}$ are pseudo-isotopic relative to $\partial M_2$.

We have
\begin{align*}
    &\Dax\big(\beta_{p,q,\gamma_1,\mathfrak{o}}(D_0), \beta_{p,q,\gamma_2,\mathfrak{o}}(D_0)\big) \\
    =& \Dax\big((\beta_{p,q,\gamma_1,\mathfrak{o}}-\beta_{p,q,\gamma_2,\mathfrak{o}})(D_0), D_0\big)\\
    =& \Dax\big( (\beta_{p,q'',\gamma_2''}  - \beta_{p,q',\gamma_2'} + \beta_{p,q',\gamma_1'}- \beta_{p,q'',\gamma_1''})(D_0),D_0\big)\\
    = & \Dax\big( (\beta_{p,q',\gamma_1'}-\beta_{p,q',\gamma_2'})(D_0),( \beta_{p,q'',\gamma_1''}-\beta_{p,q'',\gamma_2''})(D_0)\big)\\
    = &\Dax\big((\beta_{p,q',\gamma_1'}-\beta_{p,q',\gamma_2'})(D_0),D_0\big) - \Dax\big((\beta_{p,q'',\gamma_1''}-\beta_{p,q'',\gamma_2''})(D_0),D_0\big)\\
    = & \Dax\big(\beta_{p,q',\gamma_1'}(D_0),\beta_{p,q',\gamma_2'}(D_0)\big) - \Dax\big(  \beta_{p,q'',\gamma_1''}(D_0),\beta_{p,q'',\gamma_2''}(D_0)\big)
\end{align*}

Let $\gamma(t)$ $(t\in[0,1])$ be a homotopy of from $\gamma_1$ to $\gamma_2$ relative to the endpoints and let $\gamma'(t)$, $\gamma''(t)$ be the associated extensions of $\gamma(t)$. 
After a generic perturbation relative to the endpoints, we may assume that $\{\gamma'(t)\}$, $\{\gamma''(t)\}$ satisfy the conditions of Proposition \ref{prop_vertical_barbell_crossing}.
Since the interiors of $\gamma_1$ and $\gamma_2$ are both disjoint from $c$, we may shrink $c$ and the width of the 1-handle by an isotopy, and further assume that for all $t$ such that $\gamma(t)$ passes through the endpoint $q'$ of $c$, the loop obtained by going from $q'$ to $q$ along $c$, and then from $q$ to $q'$ along $\gamma(t)$ is trivial in $\pi_1(M_2)$.

 We now compute 
\begin{equation}
\label{eqn_two_Dax_compare_merian_vertical_homotopy}
    \mathfrak{p}\circ \Dax\big(\beta_{p,q',\gamma_1'}(D_0),\beta_{p,q',\gamma_2'}(D_0)\big),\,  \mathfrak{p}\circ \Dax\big(\beta_{p,q'',\gamma_1''}(D_0),\beta_{p,q'',\gamma_2''}(D_0)\big) 
\end{equation}
by Proposition \ref{prop_vertical_barbell_crossing}. The Dax invariants are given by sums contributed by the intersections of the interiors of $\gamma(t)$ with $q',q''$, and $p$. By the above assumptions, the intersections of the interiors of $\gamma(t)$ with $q',q''$ give rise to trivial loops and hence the contribution is zero. The intersections of the interior of $\gamma(t)$ with $p$ yield the same contribution to the two Dax invariants in \eqref{eqn_two_Dax_compare_merian_vertical_homotopy}. As a result, we have
\[
 \mathfrak{p}\circ \Dax\big(\beta_{p,q'',\gamma_1''}(D_0),\beta_{p,q'',\gamma_2''}(D_0)\big) -  \mathfrak{p}\circ\Dax\big(  \beta_{p,q',\gamma_1'}(D_0),\beta_{p,q',\gamma_2'}(D_0)\big) = 0,
\]
so the proposition is proved. 
\end{proof}

\subsection{Proof of Proposition \ref{prop_arc_pushing_K_Dax}}

\subsubsection{The universal cover of $\Sigma$}
It will be convenient to identify the fundamental group $\pi_1(\Sigma)$ with the set of fundamental domains in the universal cover of $\Sigma$.  
Recall that $b_0'$ is a fixed point in the interior of $H_2$ and we denote $S_0 = \{b_0'\}\times S^2\subset M$. 
Let $\widetilde{\Sigma}$ be the universal cover of $\Sigma$. The preimages of the 1-skeleton of $\Sigma$ divide $\widetilde{\Sigma}$ into fundamental domains.  Let $\tilde{b}_0'$ be a fixed lifting of $b_0'$, and let $\tilde{b}_0'$ be the base point of $\widetilde{\Sigma}$. Then every fundamental domain $D$ contains exactly one preimage $b_D$ of $\tilde{b}_0'$, and the image of the arc from $\tilde{b}_0'$ to $b_D$ defines an element in $\pi_1(\Sigma,b_0')$, which gives a bijection from the set of fundamental domains to elements in $\pi_1(\Sigma,b_0')$. Recall that $\pi_1(\Sigma,b_0')$ is canonically isomorphic to $\pi$. We use this bijection to identify the set of fundamental domains with $\pi$.  

We say that two fundamental domains are \emph{adjacent} if they have a common edge. 
We say that $g_1,g_2\in \pi$ are \emph{adjacent} if they correspond to adjacent fundamental domains. Note that $g_1,g_2\in \pi$ are adjacent if and only if they differ by the \emph{right} product of a standard generator. Therefore, if $g_1,g_2\in \pi$ are adjacent, then $h g_1,hg_2$ are also adjacent for each $h\in \pi$. 

It will also be convenient to denote the vertical and meridian-vertical barbells using points and arcs on $\widetilde{\Sigma}$. 
Suppose $\beta_{p_1,p_2,\gamma}$ is a vertical barbell diffeomorphism. Let $\tilde\gamma$ be the lifting of $\gamma$ to $\widetilde{\Sigma}$ such that the end points $p_1,p_2$ are lifted to $\tilde{p}_1$ and $\tilde{p}_2$. We define
$\beta_{\tilde{p}_1,\tilde{p}_2,\tilde{\gamma}}:= \beta_{p_1,p_2,\gamma}.$
Similarly, if $\beta_{p,q,\gamma,\mathfrak{o}}$ is a meridian-vertical barbell diffeomorphism and $(\tilde{p}, \tilde{q}, \tilde{\gamma})$ is a lift of $(p,q,\gamma)$, we will also denote $\beta_{p,q,\gamma,\mathfrak{o}}$ as  $\beta_{\tilde{p}, \tilde{q}, \tilde{\gamma},\mathfrak{o}}$. This notation will be convenient later because the lifted endpoints of an arc determine its homotopy class. For example, Proposition \ref{prop_meridian_vertical_homot} can be interpreted as the statement that the effect of a meridian-vertical barbell diffeomorphism on the Dax invariant composed with $\mathfrak{p}$ only depends on the endpoints of $\tilde{\gamma}$ and the orientations. 

Note that when $q$ is a point on the 1-cell of $\Sigma$, a meridian-vertical barbell diffeomorphism $\beta_{b_0',q,\gamma,\mathfrak{o}}$ is the same as the map $\beta_{\hat\gamma}$ in Definition \ref{defn_beta_gamma}, where $\hat\gamma$ is a lifting of $\gamma$ to an embedded arc from $S_0$ and the corresponding 1-cell of $\mathscr{I}_0(\Sigma)$ in $M_2$. Therefore, by Lemma \ref{lem_change_of_lambda0_under_barbell} and \eqref{eqn_meridian_decomp_two_vertical}, we may choose an overall orientation convention so that the following result holds.
\begin{Lemma}
\label{lem_change_of_equi_intersection_by_meridian_vertical}
    Suppose $\beta_{\tilde{b}_0',\tilde{q},\tilde{\gamma},\mathfrak{o}}$ is a meridian-vertical barbell diffeomorphism, where the projection of $(\tilde{q},\tilde{\gamma})$ to $\Sigma$ is $(q,\gamma)$. Let $\gamma',\gamma''$ be as in \eqref{eqn_meridian_decomp_two_vertical}, let $\tilde\gamma',\tilde\gamma''$ be their liftings with starting point $\tilde{b}_0'$, and let $\tilde{q}',\tilde{q}''\in \widetilde{\Sigma}$ be the ending points of the lifted arcs. Suppose $(\gamma',\gamma'')$ are ordered so that 
\begin{equation}
\label{eqn_vertical_meridian_decompose_fixed_sign}
\beta_{\tilde{b}_0',\tilde{q},\tilde{\gamma},\mathfrak{o}} = \beta_{\tilde{b}_0',\tilde{q}',\tilde{\gamma}'} - \beta_{\tilde{b}_0',\tilde{q}'',\tilde{\gamma}''}.
\end{equation}
    Let $D',D''$ be the fundamental domains containing $\tilde{q}',\tilde{q}''$, and let $g',g''\in\pi$ be the corresponding elements in the fundamental group. Then for every $[i]\in \pi_2(M_2;\partial D_0)$, we have
\begin{equation}
\label{eqn_vertical_meridian_on_homotopy_class}
\lambda_0([\beta_{\tilde{b}_0',\tilde{q},\tilde{\gamma},\mathfrak{o}}\circ i] - \lambda_0([i]) = g'-g''-(g')^{-1}+(g'')^{-1}. 
\end{equation}
\end{Lemma}

\begin{proof}
    By Lemma \ref{lem_change_of_lambda0_under_barbell}, we know that 
    \begin{equation}
    \label{eqn_vertical_meridian_on_homotopy_class_sign}
    \lambda_0([\beta_{\tilde{b}_0',\tilde{q},\tilde{\gamma},\mathfrak{o}}\circ i] - \lambda_0([i]) = \pm\big(g'-g''-(g')^{-1}+(g'')^{-1}\big), 
    \end{equation}
    where the sign depends only on $\mathfrak{o}$ and the overall orientation conventions. Note that changing the orientation $\mathfrak{o}$ flips $g_+$ and $g_-$ in Lemma \ref{lem_change_of_lambda0_under_barbell} and also flips $g'$ and $g''$ in \eqref{eqn_vertical_meridian_on_homotopy_class_sign}, so the sign in \eqref{eqn_vertical_meridian_on_homotopy_class_sign} does not depend on $\mathfrak{o}$. As a result, we may choose the orientation conventions so that the sign in \eqref{eqn_vertical_meridian_on_homotopy_class_sign} is positive. 
\end{proof}

\begin{Definition}
    Suppose $\beta_{\tilde{b}_0',\tilde{q},\tilde{\gamma},\mathfrak{o}}$ is a meridian-vertical barbell diffeomorphism. Let $g',g''$ be defined as in Lemma \ref{lem_change_of_equi_intersection_by_meridian_vertical}. We call the ordered pair $(g',g'')$ the \emph{associated adjacent pair} of $\beta_{\tilde{b}_0',\tilde{q},\tilde{\gamma},\mathfrak{o}}$.
\end{Definition}

Proposition \ref{prop_arc_pushing_K_Dax} will be proved as a corollary of the following result.  
\begin{Proposition}
\label{prop_beta_linear_combination_on_Dax}
Suppose 
\[
\beta_{\tilde{b}_0',\tilde{q}_1,\tilde{\gamma}_1,\mathfrak{o}_1},\dots,\beta_{\tilde{b}_0',\tilde{q}_m,\tilde{\gamma}_m,\mathfrak{o}_m}
\]
is a sequence of meridian-vertical barbell diffeomorphisms denoted by defining data on the universal cover of $\Sigma$. Suppose $(g'_j,g''_j)$ is the associated adjacent pair with $\beta_{\tilde{b}_0',\tilde{q}_j,\tilde{\gamma}_j,\mathfrak{o}_j}$. Suppose 
\begin{equation}
\label{eqn_change_of_equi_intersection_by_sum}
\sum_{j=1}^m \big(g'_j-g''_j-(g'_j)^{-1}+(g''_j)^{-1}\big) = 0 \in \mathbb{Z}[\pi].
\end{equation}
Let 
\[
\beta = \sum_{j=1}^m \beta_{\tilde{b}_0',\tilde{q}_j,\tilde{\gamma}_j,\mathfrak{o}_j}\in \pi_0\Diff_\partial(M_2).
\]
Then 
\[
\mathfrak{p}\circ \Dax(D_0,\beta(D_0)) = 0. 
\]
\end{Proposition}

\begin{proof}[Proof of Prop. \ref{prop_arc_pushing_K_Dax} from Prop. \ref{prop_beta_linear_combination_on_Dax}]
    By Corollaries \ref{cor_K_image_equal_K0} and \ref{cor_kerR_gen_by_rho}, we only need to show that if $\alpha\in K$ is generated by elements of the form $\hat \rho_\gamma$, then its image under the arc pushing map $f_\alpha$ satisfies $\mathfrak{p}\circ \Dax(f_\alpha(D_0),D_0)= 0$. 

    Note that the arc-pushing map of an element of the form $\hat \rho_\gamma$ is a meridian-vertical barbell diffeomorphism. We can also always isotope the defining data of a meridian-vertical barbell diffeomorphism so that the starting point of the arc is $b_0'$. So we may assume that $f_\alpha$ is generated by elements of the form $\beta_{\tilde{b}_0',\tilde{q},\tilde{\gamma},\mathfrak{o}}$. Changing the sign of $\mathfrak{o}$ turns the meridian-vertical barbell diffeomorphism to its inverse, so we may assume without loss of generality that 
    \[
f_\alpha = \sum_{j=1}^m \beta_{\tilde{b}_0',\tilde{q}_j,\tilde{\gamma}_j,\mathfrak{o}_j}\in \pi_0\Diff_\partial(M_2)
    \]
    for a sequence of meridian-vertical barbell diffeomorphisms $\beta_{\tilde{b}_0',\tilde{q}_1,\tilde{\gamma}_1,\mathfrak{o}_1}$, $\dots$, $\beta_{\tilde{b}_0',\tilde{q}_m,\tilde{\gamma}_m,\mathfrak{o}_m}$. By Lemma \ref{lem_change_of_equi_intersection_by_meridian_vertical}, the map $f_\alpha$ changes the equivariant intersection number with $D_0$ by the value on the left-hand side of \eqref{eqn_change_of_equi_intersection_by_sum}. Since we assume that $f_\alpha\in K$, the map $f_\alpha$ does not change the equivariant intersection number with $D_0$, so Equation \eqref{eqn_change_of_equi_intersection_by_sum} holds. Therefore, the statement of Proposition \ref{prop_arc_pushing_K_Dax} is implied by the conclusion of Proposition \ref{prop_beta_linear_combination_on_Dax}.
\end{proof}

\subsubsection{Admissibility of sequences of adjacent pairs}

Note that the assumption \eqref{eqn_change_of_equi_intersection_by_sum} in Proposition \ref{prop_beta_linear_combination_on_Dax} only depends on the positions of $\tilde{q}_j$ and the orientations $\mathfrak{o}_j$, but not on the isotopy classes of $\tilde{\gamma}_j$.  The orientation $\mathfrak{o}_j$ and the position of $\tilde{q}_j$ up to isotopy are both determined by the associated adjacent pair. Therefore, the assumption \eqref{eqn_change_of_equi_intersection_by_sum} only depends on the collection of the associated adjacent pairs.  The next lemma shows that the conclusion of Proposition \ref{prop_beta_linear_combination_on_Dax} also only depends on the associated adjacent pairs. 

\begin{Lemma}
\label{lem_meridian_vert_combination_Dax_only_depends_on_fund_domains}
Suppose 
\[
\beta_{\tilde{b}_0',\tilde{q}_1^a,\tilde{\gamma}_1^a,\mathfrak{o}_1^a},\dots,\beta_{\tilde{b}_0',\tilde{q}_m^a,\tilde{\gamma}_m^a,\mathfrak{o}_m^a}
\]
and
\[
\beta_{\tilde{b}_0',\tilde{q}_1^b,\tilde{\gamma}_1^b,\mathfrak{o}_1^b},\dots,\beta_{\tilde{b}_0',\tilde{q}_m^b,\tilde{\gamma}_m^b,\mathfrak{o}_m^b}
\]
are two sequences of meridian-vertical barbell diffeomorphisms. Let $((g_j')^a$, $(g_j'')^a)$ be the associated adjacent pair of $\beta_{\tilde{b}_0',\tilde{q}_j^a,\tilde{\gamma}_j^a,\mathfrak{o}_j^a}$, and let $((g_j')^b$, $(g_j'')^b)$ be the associated adjacent pair of $\beta_{\tilde{b}_0',\tilde{q}_j^b,\tilde{\gamma}_j^b,\mathfrak{o}_j^b}$. Suppose $(g_j')^a = (g_j')^b$, $(g_j'')^a = (g_j'')^b$ for all $j$, and suppose that 
\begin{align*} \sum_{j=1}^m \big((g'_j)^a-(g''_j)^a-((g'_j)^a)^{-1}+((g''_j)^a)^{-1}\big) = 0 \in \mathbb{Z}[\pi].
\end{align*}
Let 
\[
\beta^a = \sum_{j=1}^m \beta_{\tilde{b}_0',\tilde{q}_j^a,\tilde{\gamma}_j^a,\mathfrak{o}_j^a}, \quad 
\beta^b = \sum_{j=1}^m \beta_{\tilde{b}_0',\tilde{q}_j^b,\tilde{\gamma}_j^b,\mathfrak{o}_j^b}.
\]
Then
\[
\mathfrak{p}\circ \Dax(D_0,\beta^a(D_0)) = \mathfrak{p}\circ \Dax(D_0,\beta^b(D_0)).
\]
\end{Lemma}

\begin{proof}
\begin{align*}
&\mathfrak{p}\circ \Dax(D_0,\beta^a(D_0)) - \mathfrak{p}\circ \Dax(D_0,\beta^b(D_0))\\
=&\mathfrak{p}\circ \Dax(\beta^b(D_0),\beta^a(D_0))\\
=&\mathfrak{p}\circ \Dax((\beta^b-\beta^a)(D_0),D_0)\\
=& \mathfrak{p}\circ \Dax\Big(\sum_{j=1}^m (\beta_{\tilde{b}_0',\tilde{q}_j^b,\tilde{\gamma}_j^b,\mathfrak{o}_j^b}-\beta_{\tilde{b}_0',\tilde{q}_j^a,\tilde{\gamma}_j^a,\mathfrak{o}_j^a})(D_0), D_0\Big)\\
=& \sum_{j=1}^m\mathfrak{p}\circ \Dax\Big((\beta_{\tilde{b}_0',\tilde{q}_j^b,\tilde{\gamma}_j^b,\mathfrak{o}_j^b}-\beta_{\tilde{b}_0',\tilde{q}_j^a,\tilde{\gamma}_j^a,\mathfrak{o}_j^a})(D_0), D_0\Big)\\
=& \sum_{j=1}^m\mathfrak{p}\circ \Dax\Big(\beta_{\tilde{b}_0',\tilde{q}_j^b,\tilde{\gamma}_j^b,\mathfrak{o}_j^b}(D_0), \beta_{\tilde{b}_0',\tilde{q}_j^a,\tilde{\gamma}_j^a,\mathfrak{o}_j^a}(D_0)\Big),\\
\end{align*}
where the second-to-last equation follows from Lemma \ref{lem_sum_beta_Dax_rewrite} and Proposition \ref{prop_meridian_vertical_homot} (1), and 
the last line is zero by Proposition \ref{prop_meridian_vertical_homot} (2).
\end{proof}

In view of Lemma \ref{lem_meridian_vert_combination_Dax_only_depends_on_fund_domains},
we introduce the following definition.

\begin{Definition}
\label{defn_admissible_sequence_pixpi}
    Suppose $\{(g_j',g_j'')\}_{1\le j\le m}$ is a sequence of adjacent pairs in $\pi$. Suppose 
\begin{equation}
\label{eqn_change_of_equi_intersection_by_sum_repeat}
\sum_{j=1}^m \big(g'_j-g''_j-(g'_j)^{-1}+(g''_j)^{-1}\big) = 0 \in \mathbb{Z}[\pi].
\end{equation}
We say that the sequence $\{(g_j',g_j'')\}_{1\le j\le m}$ is \emph{admissible}, if there exists a sequence of meridian-vertical barbell diffeomorphisms
$\{\beta_{\tilde{b}_0',\tilde{q}_j,\tilde{\gamma}_j,\mathfrak{o}_j}\}_{1\le j\le m}$
whose associated adjacent pairs are $\{(g_j',g_j'')\}_{1\le j\le m}$, such that
\begin{equation}
\label{eqn_admissible_Dax}
    \mathfrak{p}\circ \Dax(\sum_{j=1}^m \beta_{\tilde{b}_0',\tilde{q}_j,\tilde{\gamma}_j,\mathfrak{o}_j}(D_0),D_0) = 0.
\end{equation}
\end{Definition}

\begin{remark}
\label{rmk_meridian_vert_combination_Dax_only_depends_on_fund_domains}
By Lemma \ref{lem_meridian_vert_combination_Dax_only_depends_on_fund_domains}, Proposition \ref{prop_beta_linear_combination_on_Dax} holds if and only if every sequence $\{(g_j',g_j'')\}_{1\le j\le m}$ satisfying \eqref{eqn_change_of_equi_intersection_by_sum_repeat} holds is admissible. 
\end{remark}

\begin{remark}
\label{rmk_disjoint_union_admissible}
It is clear from the definition that changing the order of the sequence $\{(g_j',g_j'')\}_{1\le j\le m}$ does not affect admissibility. If two sequences are admissible, then their disjoint union is also admissible. Moreover, if the disjoint union of two sequences $\mathcal{G}_1$ and $\mathcal{G}_2$ is admissible and both $\mathcal{G}_1, \mathcal{G}_2$ satisfy \eqref{eqn_change_of_equi_intersection_by_sum_repeat}, then $\mathcal{G}_1$ is admissible if and only if $\mathcal{G}_2$ is admissible. 
\end{remark}

\subsubsection{Verification of admissibility}
Now we prove the admissibility of two special classes of sequences.

\begin{Lemma}
\label{lem_cycle_admissible}
    Let $g_1,\dots,g_{m+1}$ be a sequence of elements in $\pi$ such that $g_1=g_{m+1}$, and $g_j$ is adjacent to $g_{j+1}$ for all $1\le j \le m$. Consider the sequence 
    $\mathcal{G} = \{(g_{j+1},g_j)\}_{1\le j\le m}.$
    Then $\mathcal{G}$ satisfies \eqref{eqn_change_of_equi_intersection_by_sum_repeat} and is admissible. 
\end{Lemma}

\begin{proof}
It is clear that \eqref{eqn_change_of_equi_intersection_by_sum_repeat} holds for $\mathcal{G}$. We need to find a sequence of meridian-vertical barbell diffeomorphisms that satisfies \eqref{eqn_admissible_Dax}. 

Let $b_0''$ be an interior point of $H_2$ that is different from $b_0'$. Let $\tilde{p}\in\widetilde{\Sigma}$ be the lifting of $b_0''$ in the fundamental domain corresponding to $g_1$. Let $\tilde{\gamma}_1$ be an arc from $\tilde{b}_0'$ to $\tilde{p}$ such that its interior does not pass through the preimages of $b_0'$ or $b_0''$. Let $\tilde{\gamma}_2$ be an arc from $\tilde{p}$ to itself such that its interior does not pass through the preimages of $b_0'$ or $b_0''$. Let $\gamma_1,\gamma_2$ be the projections of $\tilde{\gamma}_1$ and $\tilde{\gamma}_2$ to $\Sigma$. We further require that 
\begin{enumerate}
    \item The trajectory of $\tilde{\gamma}_2$ goes through the fundamental domains corresponding to $g_1,\dots,g_{m+1}$ in the given order.
    \item The intersection of ${\gamma}_2$ with the 1-handles of $\Sigma$ are co-core arcs.
    \item The concatenation of ${\gamma}_1$ and ${\gamma}_2$ in $\Sigma$ is a smooth arc and only has transverse double self-intersections.
\end{enumerate}

Let the arc $\tilde{\gamma}_3$ be obtained by first concatenating $\tilde{\gamma}_1$ and $\tilde{\gamma}_2$, and then shrink a little bit near the end point at $\tilde{p}$ so that the end point of the projection $\gamma_3$ does not lie in the interior of the arc. 

In general, suppose $\gamma$ is a smooth oriented arc on $\Sigma$ with only transverse double self-intersections, we define a self-intersection number of $\gamma$ that takes value in $\mathbb{Z}[\mathcal{C}]$ as follows.
For each self-intersection point $z$ of $\gamma$, let $c(z)$ denote the oriented loop that goes from $z$ to $z$ along $\gamma$, with the same orientation as $\gamma$.  Then $c(z)$ as a free loop on $\Sigma$ defines a conjugation class in $\pi$. Define $[c(z)]\in \mathbb{Z}[\mathcal{C}]$ to be the conjugation class of $c(z)$ if $c(z)$ is non-trivial, and zero if $c(z)$ is trivial. Let $\epsilon(z)\in\{1,-1\}$ be the sign of the self-intersection point. By definition, this is the orientation of the basis formed by the two tangent vectors of $\gamma$ at the self-intersection, ordered according to the orientation of the arc. Define 
\begin{equation}
\label{eqn_I_gamma3}
I(\gamma) = \sum_z \epsilon(z)[c(z)] \in\mathbb{Z}[\mathcal{C}],
\end{equation}
where $z$ goes over all self-intersection points of $\gamma$. 

Consider the vertical barbell diffeomorphisms $\beta_{\tilde{b}_0',\tilde{p},\tilde{\gamma}_1}$ and $\beta_{\tilde{b}_0',\tilde{p},\tilde{\gamma}_3}$. By Proposition \ref{prop_vertical_barbell_crossing} (1), we know that $\beta_{\tilde{b}_0',\tilde{p},\tilde{\gamma}_1}$ and $\beta_{\tilde{b}_0',\tilde{p},\tilde{\gamma}_3}$ are pseudo-isotopic relative to $\partial M_2$. 
We compute the value of
\[
\mathfrak{p}\circ\Dax(\beta_{\tilde{b}_0',\tilde{p},\tilde{\gamma}_3}(D_0),  \beta_{\tilde{b}_0',\tilde{p},\tilde{\gamma}_1}(D_0))\in \mathbb{Z}[\mathcal{C}]
\]
using two different methods. The fact that the two methods give the same result will show that $\{(g_j',g_j'')\}_{1\le j\le m}$ is admissible. 

First, let $\tilde{\gamma}(t)$ $(t\in[0,1])$ be a homotopy from $\tilde{\gamma}_1$ to $\tilde{\gamma}_3$ along $\tilde{\gamma}_2$. By definition, the trajectory of the endpoint $\tilde{\gamma}(t)(1)$ is contained in $\tilde{\gamma}_2$. Let ${\gamma}(t)$ be the projection of $\tilde{\gamma}(t)$ to $\Sigma$. 

When $\gamma(t)(1)$ is contained in the interior of $H_2$ and does not lie in the interior of $\gamma(t)$, we have $\beta_{\tilde{b}_0',\tilde{\gamma}(t)(1),\tilde{\gamma}(t)}$ is a well-defined vertical barbell. Define
\[
\beta(t) = \beta_{\tilde{b}_0',\tilde{\gamma}(t)(1),\tilde{\gamma}(t)}.
\]
The isotopy class of $\beta(t)$ changes whenever the endpoint $\gamma(t)(1)$ goes across the interior of $\gamma(t)$ or across a 1-handle of $\Sigma$. By Propositions \ref{prop_vertical_barbell_crossing}, if goes across the interior of the arc, then $\beta(t)$ changes by a diffeomorphism of the form $\epsilon^{(1)}\epsilon(z)\beta(z)$, where $\epsilon^{(1)}\in\{1,-1\}$ is a universal constant and $\beta(z)$ is a self-referential barbell diffeomorphism. The loop $c(z)$ associated with $\beta(z)$ as in Definition \ref{defn_associated_loop_self_referential_barbell} coincides with the loop $c(z)$ in \eqref{eqn_I_gamma3} when viewing $z$ as a self-intersection point of $\gamma_3$.

By \eqref{eqn_meridian_decomp_two_vertical}, if the trajectory of the end point goes across a 1-handle, say, if the endpoint moves from the fundamental domain of $g_j$ to that of $g_{j+1}$, then there is a meridian-vertical barbell diffeomorphism $\beta_{\tilde{b}_0',\tilde{q}_j,\tilde{\gamma}_j,\mathfrak{o}_j}$ such that 
\[
\beta_{\tilde{b}_0',\tilde{q}_j,\tilde{\gamma}_j,\mathfrak{o}_j} = \beta_{\tilde{b}_0',\tilde{\gamma}(t+\delta)(1),\tilde{\gamma}(t+\delta)} -\beta_{\tilde{b}_0',\tilde{\gamma}(t-\delta)(1),\tilde{\gamma}(t-\delta)},
\]
where $t-\delta$ and $t+\delta$ are the values of $t$ before and after $\gamma(t)$ crosses the 1-handle. The adjacent pair associated with $\beta_{\tilde{b}_0',\tilde{q}_j,\tilde{\gamma}_j,\mathfrak{o}_j}$ is $(g_{j+1}, g_j)$.

As a consequence, there exists a sequence of meridian-vertical barbells diffeomorphisms $\{\beta_{\tilde{b}_0',\tilde{q}_j,\tilde{\gamma}_j,\mathfrak{o}_j}\}_{1\le j\le m}$ whose associated adjacent pairs are equal to $\{(g_{j+1}, g_j)\}_{1\le j\le m}$, such that 

\begin{equation}
\label{eqn_diff_beta_gamma3_gamma1_method1}
\beta_{\tilde{b}_0',\tilde{p},\tilde{\gamma}_3}- \beta_{\tilde{b}_0',\tilde{p},\tilde{\gamma}_1} = \sum_{j=1}^m \beta_{\tilde{b}_0',\tilde{q}_j,\tilde{\gamma}_j,\mathfrak{o}_j} + \sum_{z} \epsilon^{(1)}\epsilon(z)\beta(z),
\end{equation}
where the second sum takes over all self-intersection points $z$ of $\gamma_3$ such that at least one sheet is contained in $\gamma_3\setminus \gamma_1$. By Lemma \ref{lem_self_referential_barbell_Dax}, the second sum is pseudo-isotopic to the identity relative to $\partial M_2$, and we have
\[
\mathfrak{p}\circ \Dax\Big(\big( \sum_{z} \epsilon^{(1)}\epsilon(z)\beta(z)\big) (D_0), D_0\Big) = \epsilon^{(1)}\epsilon\cdot (I(\gamma_3) - I(\gamma_1)),
\]
where $\epsilon\in \{1,-1\}$ is the constant in Lemma \ref{lem_self_referential_barbell_Dax}.
Therefore, by \eqref{eqn_diff_beta_gamma3_gamma1_method1} and Lemma \ref{lem_sum_beta_Dax_rewrite}, we have 
\begin{align}
&\mathfrak{p}\circ\Dax(\beta_{\tilde{b}_0',\tilde{p},\tilde{\gamma}_3}(D_0),  \beta_{\tilde{b}_0',\tilde{p},\tilde{\gamma}_1}(D_0)) 
\nonumber \\
=\,&\mathfrak{p}\circ\Dax((\beta_{\tilde{b}_0',\tilde{p},\tilde{\gamma}_3}-\beta_{\tilde{b}_0',\tilde{p},\tilde{\gamma}_1})(D_0),  D_0) 
\nonumber \\
=\,& \mathfrak{p}\circ \Dax\Big(\big( \sum_{z} \epsilon^{(1)}\epsilon(z)\beta(z)\big) (D_0), D_0\Big) 
+   \mathfrak{p}\circ \Dax(\sum_{j=1}^m \beta_{\tilde{b}_0',\tilde{q}_j,\tilde{\gamma}_j,\mathfrak{o}_j}(D_0),D_0)\nonumber\\
=\,& \epsilon^{(1)}\epsilon\cdot (I(\gamma_3) - I(\gamma_1))+    \mathfrak{p}\circ \Dax(\sum_{j=1}^m \beta_{\tilde{b}_0',\tilde{q}_j,\tilde{\gamma}_j,\mathfrak{o}_j}(D_0),D_0).
\label{eqn_Dax_gamma3_gamma1_method1}
\end{align}

Now, we use another method to compute the left-hand side of \eqref{eqn_Dax_gamma3_gamma1_method1}. Let $H$ be a homotopy from $\tilde{\gamma}_1$ to $\tilde{\gamma}_3$ that fixes the starting point at $\tilde{b}_0'$ and such that the trajectory of the ending point does not intersect the preimage of the 1-handles of $\Sigma$. Such a homotopy exists because the endpoints of $\tilde{\gamma}_1$ and $\tilde{\gamma}_3$ are in the same fundamental domain. We also assume that $H$ satisfies the transversality conditions in the setup of Proposition \ref{prop_vertical_barbell_crossing}. Then, there are only finitely many times in the homotopy such that the interior of the arc goes across an endpoint. For each such crossing point $z$, view $z$ as a self-intersection point of the arc, let $c(z)$ be the associated loop defined as above, let $\epsilon^H(z)\in\{1,-1\}$ be the sign defined by the orientation of \eqref{eqn_sign_vertical_barbell_crossing}. Let
\[
I(H) = \sum_z \epsilon^H(z) c(z).
\]
Then by Proposition \ref{prop_vertical_barbell_crossing}, we have
\begin{equation}
\label{eqn_Dax_gamma3_gamma1_method2}
\mathfrak{p}\circ\Dax(\beta_{\tilde{b}_0',\tilde{p},\tilde{\gamma}_3}(D_0),  \beta_{\tilde{b}_0',\tilde{p},\tilde{\gamma}_1}(D_0)) 
 = \epsilon^{(2)} I(H),
\end{equation}
where $\epsilon^{(2)}\in\{1,-1\}$ is the constant in Proposition \ref{prop_vertical_barbell_crossing}.

Now we show that $I(H) = I(\gamma_3) - I(\gamma_1)$. This is because the homotopy $H$ can be decomposed into a sequence of the following operations:
(1) crossings of an endpoint with the interior,
(2) Whitney moves,
(3) twist moves,
(4) regular homotopies.
Whitney moves, twist moves, and regular homotopies do not change the value of the self-intersection number in $\mathbb{Z}[\mathcal{C}]$ of the arc. It is straightforward to track the signs and verify that the change of the self-intersection number from crossings of an endpoint with the interior is given by $I(H)$. Therefore,
\begin{equation}
\label{eqn_I(H)_by_two_arcs}
I(H) = I(\gamma_3) - I(\gamma_1).
\end{equation}
By Proposition \ref{prop_vertical_barbell_crossing}, $\epsilon^{(2)} = \epsilon^{(1)}\epsilon$. It then follows immediately from \eqref{eqn_Dax_gamma3_gamma1_method1}, \eqref{eqn_Dax_gamma3_gamma1_method2}, \eqref{eqn_I(H)_by_two_arcs} that 
\[
  \mathfrak{p}\circ \Dax(\sum_{j=1}^m \beta_{\tilde{b}_0',\tilde{q}_j,\tilde{\gamma}_j,\mathfrak{o}_j}(D_0),D_0)=0,
\]
so $\mathcal{G}$ is admissible.
\end{proof}

Now we establish the admissibility for another class of sequences. Recall that if $g_1,g_2\in \pi$ are adjacent, then  $hg_1,hg_2$ are also adjacent for each $h\in \pi$. 
\begin{Lemma} 
\label{lem_symmetric_ray_admissible}
Let $g_1,\dots,g_{m}$ be a sequence of elements in $\pi$ such that $g_j$ is adjacent to $g_{j+1}$ for all $1\le j \le m-1$. Here, we do not require that $g_m=g_1$. Consider the following two sequences of adjacent pairs 
$
\mathcal{G}_1 = \{(g_1^{-1}g_{j+1},g_1^{-1}g_{j})\}_{1\le j \le m-1}
$
and 
$
\mathcal{G}_2 = \{(g_m^{-1}g_{j},g_{m}^{-1}g_{j+1})\}_{1\le j \le m-1}.
$
Let $\mathcal{G}_1\sqcup \mathcal{G}_2$ be the disjoint union of $\mathcal{G}_1$ and $\mathcal{G}_2$. Then $\mathcal{G}_1\sqcup \mathcal{G}_2$ satisfies \eqref{eqn_change_of_equi_intersection_by_sum_repeat} and is admissible.
\end{Lemma}

\begin{proof}
We first verify that $\mathcal{G}_1\sqcup \mathcal{G}_2$ satisfies \eqref{eqn_change_of_equi_intersection_by_sum_repeat}.
    Let $\sigma:\mathbb{Z}[\mathcal{C}]\to \mathbb{Z}[\mathcal{C}]$ be the involution that takes $[g]$ to $[g^{-1}]$. For a sequence 
$\mathcal{G} = \{(g_j',g_j'')\}_j,$
    let 
    \[
    s(\mathcal{G}) = \sum_{j} ([g_j'] - [g_j''])\in \mathbb{Z}[\mathcal{C}].
    \]
    Then the left-hand side of \eqref{eqn_change_of_equi_intersection_by_sum_repeat} applied to  $\mathcal{G}_1\sqcup \mathcal{G}_2$ is equal to 
    \[
 s(\mathcal{G}_1) + s(\mathcal{G}_2) -\sigma\circ s(\mathcal{G}_1) -\sigma\circ s(\mathcal{G}_2).
    \]
    Note that after cancellation, we have
    \[
    s(\mathcal{G}_1) = [g_1^{-1}g_m] - [g_1^{-1}g_1]= [g_1^{-1}g_m],
    \]
    \[
    s(\mathcal{G}_2) = [g_m^{-1}g_1] - [g_m^{-1}g_m]= [g_m^{-1}g_1],
    \]
    so $ s(\mathcal{G}_1) + s(\mathcal{G}_2)$ is invariant under $\sigma$, and $\mathcal{G}_1\sqcup \mathcal{G}_2$ satisfies \eqref{eqn_change_of_equi_intersection_by_sum_repeat}.

Now we show that $\mathcal{G}_1\sqcup \mathcal{G}_2$ is admissible. We use a similar argument as Proposition \ref{lem_cycle_admissible}. Let $\beta_{\tilde{p},\tilde{q},\tilde{\gamma}}$ be a vertical barbell diffeomorphism, denoted by defining data on the universal cover of $\Sigma$. 
Let $\gamma$ be the projection of $\tilde{\gamma}$. We further require that
\begin{enumerate}
    \item The trajectory of $\tilde{\gamma}$ goes through $g_1,\dots,g_m$ in the given order. 
    \item $\gamma$ only has transverse double self-intersections.
    \item $\gamma$ intersects the 1-handles of $\Sigma$ at co-core arcs. 
\end{enumerate} 
We compute $\beta_{\tilde{p},\tilde{q},\tilde{\gamma}}$ using two methods.

First, let $\tilde{\gamma}(t)$ $(t\in[0,1])$ be a homotopy from the constant path at $\tilde{p}$ to $\tilde{\gamma}$ along $\tilde\gamma$. When the end point $\tilde{\gamma}(t)(1)$ is not in the interior of $\tilde{\gamma}(t)$ or the preimage of a 1-handle of $\Sigma$, the vertical barbell diffeomorphism 
$\beta(t) := \beta_{\tilde{p},\tilde{\gamma}(t)(1),\tilde{\gamma}(t)}$
is well defined. When $t$ is close to $0$, we know that $\beta(t)$ is isotopic to the identity by Corollary \ref{cor_parallel_barbell_compose_to_id}. When changing $t$ from $0$ to $1$, the isotopy class of $\beta(t)$ changes when the endpoint goes across an interior point or across a 1-handle. For each double point $z$ of $ \gamma$, let $c(z)$ be the oriented loop from $z$ to $z$ along $\gamma$, with the same orientation of $\gamma$. Let $\epsilon(z)$ be the sign of the self-intersection of $\gamma$, defined by the orientation of the two tangent vectors of $\gamma$ at $z$ ordered according to the orientation of $\gamma$. Then, by Propositions \ref{prop_vertical_barbell_crossing} and \ref{prop_meridian_vertical_homot}, we have 
\begin{equation}
\label{eqn_beta_pqr_decompse_1}
\beta_{\tilde{p},\tilde{q},\tilde{\gamma}} = \beta(1) = \sum_{j=1}^{m-1} \beta_j  + \sum_{z} \epsilon^{(1)} \epsilon(z) \beta(z),
\end{equation}
where the summation in $z$ goes over all self-intersections of $\gamma$,
the term $\beta(z)$ is a self-referential barbell diffeomorphism with associated loop $c(z)$, the universal constant $\epsilon^{(1)}$ is given by Proposition \ref{prop_vertical_barbell_crossing}, and $\beta_j$ is a meridian-vertical barbell whose associated adjacent pair is $(g_1^{-1}g_{j+1}, g_1^{-1}g_j)$. 
Note that to find the associated adjacent pair of $\beta_j$, we need to take a (left) translation by $g_1^{-1}$ on $\widetilde{\Sigma}$ so that the starting point of the arcs defining the barbells are in the fundamental domain containing the base point $\tilde{b}_0'$. This condition is required by Definition \ref{defn_admissible_sequence_pixpi} and is needed for Lemma \ref{lem_change_of_equi_intersection_by_meridian_vertical}. 

Now we compute $\beta_{\tilde{p},\tilde{q},\tilde{\gamma}}$ by a different method. Let $\overline{(\tilde{\gamma})}$ be the orientation-reversal of $\gamma$, then 
$
\beta_{\tilde{p},\tilde{q},\tilde{\gamma}}+\beta_{\tilde{q},\tilde{p},\overline{(\tilde{\gamma})}}=0,
$
where $0\in \pi_0\Diff_\partial(M_2)$ denotes the isotopy class of the identity map. 
Applying \eqref{eqn_beta_pqr_decompse_1} to $\beta_{\tilde{q},\tilde{p},\overline{(\tilde{\gamma})}}$, we have
\begin{equation}
\label{eqn_beta_pqr_decompse_2}
\beta_{\tilde{q},\tilde{p},\overline{(\tilde{\gamma})}} =  \sum_{j=1}^{m-1} \beta_j'  + \sum_{z} \epsilon^{(1)} (-\epsilon(z)) \beta'(z),
\end{equation}
where $\beta'(z)$ is a self-referential barbell diffeomorphism with associated loop the orientation-reversal of $c(z)$, and $\beta_j'$ is a meridian-vertical barbell diffeomorphism whose associated adjacent pair is $(g_m^{-1}g_j, g_m^{-1}g_{j+1})$. 
There is a negative sign in front of $\epsilon(z)$ in \eqref{eqn_beta_pqr_decompse_2} because when we reverse the orientation of the arc, both tangent vectors at a double point are reversed and their order is also reversed, so the sign of the double point is reversed. 

Therefore, by Lemma \ref{lem_sum_beta_Dax_rewrite} and \eqref{eqn_beta_pqr_decompse_1}, \eqref{eqn_beta_pqr_decompse_2}, we have
\begin{align*}
    0 = & \Dax ((\beta_{\tilde{p},\tilde{q},\tilde{\gamma}}+\beta_{\tilde{q},\tilde{p},\overline{(\tilde{\gamma})}})(D_0), D_0) \\
    = & \Dax\Big(\big(\sum_{j=1}^{m-1}(\beta_j+\beta_j')+ \sum_{z} \epsilon^{(1)} \epsilon(z) (\beta(z)-\beta'(z))\big)(D_0),D_0\Big)\\
    = & \Dax\Big(\big(\sum_{j=1}^{m-1}(\beta_j+\beta_j')\big)(D_0),D_0\Big) + \sum_{z} \epsilon^{(1)}  \epsilon(z)\Dax\big(\beta(z)(D_0), D_0\big) \\
    &\qquad -\sum_{z} \epsilon^{(1)}  \epsilon(z)\Dax\big(\beta'(z)(D_0), D_0\big) 
\end{align*}
By Lemma \ref{lem_self_referential_barbell_Dax}, 
\[
\mathfrak{p}\circ \Dax\big(\beta(z)(D_0), D_0\big) = [c(z)] + [c(z)^{-1}] = \mathfrak{p}\circ \Dax\big(\beta'(z)(D_0), D_0\big).
\]
So we conclude that 
\[
\mathfrak{p}\circ \Dax\Big(\big(\sum_{j=1}^{m-1}(\beta_j+\beta_j')\big)(D_0),D_0\Big) =0.
\]
Therefore, $\mathcal{G}_1\sqcup\mathcal{G}_2$ is admissible. 
\end{proof}

\subsubsection{Finish of the proof}
We use Lemmas \ref{lem_cycle_admissible} and \ref{lem_symmetric_ray_admissible} to show that every sequence of adjacent pairs satisfying \eqref{eqn_change_of_equi_intersection_by_sum_repeat} is admissible.  

\begin{Definition}
   We call the sequence $\mathcal{G}$ in Lemma \ref{lem_cycle_admissible} a \emph{loop} of adjacent pairs.    We call the sequence $\mathcal{G}_1\sqcup \mathcal{G}_2$ in Lemma \ref{lem_symmetric_ray_admissible} a \emph{symmetric ray} of adjacent pairs. 
\end{Definition}

\begin{Definition}
    Let $\mathcal{G}=\{(g_j',g_j'')\}_{1\le j\le m}$ be a sequence of adjacent pairs. Define 
    $-\mathcal{G} = \{(g_j'',g_j')\}_{1\le j\le m}$.
\end{Definition}

\begin{Lemma}
    If $\mathcal{G}$ is admissible, then $-\mathcal{G}$ is admissible.
\end{Lemma}

\begin{proof}
    Reversing the orientation of the meridian of a meridian-vertical barbell diffeomorphism changes the sign of the barbell diffeomorphism and flips the order of its associated adjacent pair. 
\end{proof}

\begin{Lemma}
    \label{lem_decompose_to_rays_and_loops}
    Suppose $\mathcal{G}_0 = \{(g_j',g_j'')\}_{1\le j\le m}$ is a sequence of adjacent pairs satisfying \eqref{eqn_change_of_equi_intersection_by_sum_repeat}. Then there exists a finite collection of sequences $\mathcal{G}_1,\dots,\mathcal{G}_k$, such that 
    \begin{enumerate}
        \item   For each $j$, either $\mathcal{G}_j$ or $-\mathcal{G}_j$ is a symmetric ray of adjacent pairs or its negative,
        \item   $\mathcal{G}_0 \sqcup \mathcal{G}_1\sqcup \dots\sqcup \mathcal{G}_k$
    can be decomposed as the disjoint union of a finite collection of loops of adjacent pairs. 
    \end{enumerate}
\end{Lemma}

\begin{proof}
    For each adjacent pair $(g',g'')$, define $s(g',g'') = [g]-[g'']\in\mathbb{Z}[\mathcal{C}]$. For a sequence $\mathcal{G}$ of adjacent pairs, define $s(\mathcal{G})$ to be the sum of $s(-)$ for all terms in the sequence. Let $\sigma:\mathbb{Z}[\mathcal{C}]\to \mathbb{Z}[\mathcal{C}]$ be the involution that takes $[g]$ to $[g^{-1}]$. Then we have the following observations:
    \begin{enumerate}
        \item $\mathcal{G}$ satisfies \eqref{eqn_change_of_equi_intersection_by_sum_repeat} if and only if $s(\mathcal{G})$ is invariant under $\sigma$.
        \item $\mathcal{G}$ is the disjoint union of a finite collection of loops of adjacent pairs if and only if $s(\mathcal{G})=0$.
        \item The symmetric ray of adjacent pairs $\mathcal{G}$ associated with the adjacent sequence $g_1, \dots, g_m$ satisfies $s(\mathcal{G})=[g_1^{-1}g_m]+ [g_mg_1^{-1}]$.
        \item $s(-\mathcal{G}) = -s(\mathcal{G})$. 
    \end{enumerate}
    By Observation (1) and the assumptions, we know that $s(\mathcal{G}_0)$ is invariant under $\sigma$. By Remark \ref{rmk_involution_fixed_point_C}, we know that the action of $\sigma$ has no fixed point on $\mathcal{C}$. Therefore, $s(\mathcal{G}_0)$ is generated by elements of the form $[g]+[g^{-1}]$. By Observations (3), (4), there exists a finite collection of sequences $\mathcal{G}_1,\dots,\mathcal{G}_k$, each one is either a symmetric ray of adjacent pairs or its negative, such that $\sum_{j=0}^k s(\mathcal{G}_j)=0$. By Observation (2), $\mathcal{G}_0 \sqcup \mathcal{G}_1\sqcup \dots\sqcup \mathcal{G}_k$
    can be decomposed as the disjoint union of a finite collection of loops of adjacent pairs. 
\end{proof}

We can now finish the proof of Proposition \ref{prop_beta_linear_combination_on_Dax}. Recall that this implies Proposition \ref{prop_arc_pushing_K_Dax}.
\begin{proof}[Proof of Proposition \ref{prop_beta_linear_combination_on_Dax}]
 By Lemma \ref{lem_meridian_vert_combination_Dax_only_depends_on_fund_domains} (see also Remark \ref{rmk_meridian_vert_combination_Dax_only_depends_on_fund_domains}), we only need to show that every sequence of adjacent pairs satisfying \eqref{eqn_change_of_equi_intersection_by_sum_repeat} is admissible. Suppose $\mathcal{G}_0$ satisfies \eqref{eqn_change_of_equi_intersection_by_sum_repeat}. Then by Lemma \ref{lem_decompose_to_rays_and_loops} and the earlier results, one can take a disjoint union of $\mathcal{G}_0$ with an admissible sequence to obtain another admissible sequence. By Remark \ref{rmk_disjoint_union_admissible}, $\mathcal{G}_0$ is admissible. 
\end{proof}

\section{The mapping class group of $\Sigma\times S^2$}
\label{sec_Mcg_Sigma_S2}

In this section, we prove Theorem \ref{thm: MCG infinitely generated} using the relative Dax invariant constructed in Definition \ref{def_Dax(i1,i2)}.

Recall that $M=\Sigma\times S^2$ where $\Sigma$ is a closed orientable surface with genus $\ell\ge 1$.
Let $\Mcg(M)$ be the smooth mapping class group of $M$. Let $\Mcg_*(M)$ be the pointed version of the smooth mapping class group; that is, $\Mcg_*(M)$ denotes the group of diffeomorphisms of $M$ fixing the base point $b$, up to smooth isotopies fixing $b$. Let $\operatorname{Aut}(M)$ be the group of homotopy equivalences from $M$ to $M$ up to homotopies, and let $\operatorname{Aut}_*(M)$ be the group of base-point-preserving homotopy equivalences from $M$ to $M$ up to base-point-preserving homotopies. Consider the following commutative diagram
$$
\xymatrix{
  \Mcg_*(M)\ar[r]^{p_0}\ar[d] & \operatorname{Aut}_*(M)\ar[d]\\
  \Mcg(M)\ar[r]^{p_1} & \operatorname{Aut}(M).  
}
$$

We use $\Mcg_{0}(M)$ to denote $\ker p_1$. The next lemma shows that $\ker p_1$ is canonically isomorphic to $\ker p_0$.

\begin{Lemma}\label{lem: Mcg pointed and unpointed}
The map $\iota: \ker p_0\to \ker p_1$ induced by the forgetful map $\Mcg_*(M)\to \Mcg(M)$ is an isomorphism.  
\end{Lemma}
\begin{proof}  
The proof is similar to Lemma \ref{lem_cS_cS'_bijection}. We first show that $\iota$ is surjective. 
Suppose $f:M\to M$ is a diffeomorphism that is homotopic to the identity via a homotopy $h:I\times M\to M$. We may assume that $h$ is smooth. Then $\gamma=h(I\times \{b\})$ is a path from $f(b)$ to $b$. We apply the isotopy extension theorem to $\gamma$ and obtain a diffeomorphism $g:M\to M$ that is isotopic to the identity and sends $f(b)$ to $b$. Then $[g\circ f]\in \ker p_0$ and its image in $\Mcg(M)$ is equal to $[f]$. This shows that $\iota$ is surjective. 

Now we prove the injectivity of $\iota$. Suppose $[f]\in \ker p_0$ and $\iota[f]=1\in \Mcg(M)$. Then there exists an isotopy $h_1:I\to \Diff(M)$ from $f$ to the identity. Let $\operatorname{ev}:\Diff(M)\to M$ be the evaluation map at $b$ and let $\gamma=\operatorname{ev}\circ h_1$. Then we have the commutative diagram
\[
\xymatrix
{\pi_{1}(M,b)\ar[rd]^{\id}\ar[r]^{\operatorname{f_*=\id}}& \pi_{1}(M,b)\ar[d]^{\gamma^{\sharp}}\\
& \pi_{1}(M,b)
}
\]
Hence $\gamma$ represents an element in the center of $\pi_{1}(M)$. 

If the genus $\ell>1$, the group $\pi_{1}(M)$ has a trivial center. So there exists a homotopy $h_2: I\times I\to M$ from $\gamma$ to the identity. We apply the homotopy lifting property
\[
\xymatrix
{(\{0\}\times I)\cup (I\times \{0\})\cup (\{1\}\times I)\ar[d]\ar[rr]^-{f\cup h_1\cup \operatorname{id}}& &\Diff(M)\ar[d]^{\operatorname*{ev}}\\
I\times I\ar[rr]^{h_2}\ar@{.>}[urr]^{\widetilde{h}_2}& & M
}
\]
By restricting $\widetilde{h}_{2}$ to $I\times \{1\}$, we get the pointed isotopy from $f$ to the identity. 

If $\ell=1$, let $\gamma_1: I\to \Sigma$ be the projection of $\gamma$ to $\Sigma$. Note that since $\Sigma$ has genus $1$, it admits an abelian Lie group structure. Consider
\[
\rho_{\gamma_{1}(t)}: M\to M,\quad (x,y)\mapsto (x-\gamma_{1}(t),y).
\]
Let $h'_{1}:I\to \Diff(M)$ be the isotopy defined by 
$h'_{1}(t)=\rho_{\gamma_{1}(t)}\circ h_1(t).$
The loop $\gamma'=\operatorname{ev}\circ h_1'$ projects to a constant path to $\Sigma$ and hence is null-homotopic. So the desired result follows from the same argument as above. 
\end{proof}

Now we prove the following result.
\begin{Proposition}\label{prop_surj_to_Zinfty}
There exist surjective homomorphisms 
\[
\begin{split}
\Phi_0&:\Mcg_0(M)\cong \ker p_0 \twoheadrightarrow \mathbb{Z}^{\infty}\\
\Phi_1&:\Mcg_*(M)\twoheadrightarrow \mathbb{Z}^{\infty}\\
\Phi_2&:\Mcg(M)\twoheadrightarrow \mathbb{Z}^{\infty}.
\end{split}
\]
\end{Proposition}

We start with the construction of $\Phi_0$. Recall that $\mathscr{I}_0:\Sigma\to M$ denotes the standard embedding of $\Sigma$ into $M$. Suppose $i_1,i_2$ are embeddings from $\Sigma \to M$ that are homotopic to $\mathscr{I}_0$ relative to $b_0$, recall that the relative Dax invariant $\Dax(i_1,i_2)\in \mathbb{Z}[\mathcal{C}]$ was defined in Definition \ref{def_Dax(i1,i2)}. 

\begin{Definition}
\label{defn_Psi0_MCG0}
Define $\Psi_0: \ker p_0\to \bZ[\mathcal{C}]^\sigma$ by 
\[
\Psi_0([f]) = \Dax(\mathscr{I}_0,f\circ \mathscr{I}_0). 
\]
\end{Definition}

\begin{Lemma}
The map $\Psi_0$ is a group homomorphism. 
\end{Lemma}

\begin{proof}
Suppose $f,g\in \Mcg_0(M)$. By Lemma \ref{lem_Dax_change_by_MCG0_*}, 
\[
\Dax(\mathscr{I}_0,g\circ \mathscr{I}_0) = \Dax(f\circ \mathscr{I}_0, f\circ g\circ \mathscr{I}_0).
\]
By the additivity of the Dax invariant,
\begin{align*}
\Psi_0([f\circ g]) & = \Dax(\mathscr{I}_0,f\circ g\circ \mathscr{I}_0) \\
& = \Dax(\mathscr{I}_0,g\circ \mathscr{I}_0) + \Dax(g\circ \mathscr{I}_0,f\circ g\circ \mathscr{I}_0) 
\\
& = \Dax(\mathscr{I}_0,g\circ \mathscr{I}_0) + \Dax(\mathscr{I}_0,f\circ \mathscr{I}_0) = \Psi_0([f]) + \Psi_0([g]).
\end{align*}
So $\Psi_0$ is a group homomorphism. 
\end{proof}

By Lemma \ref{lem_self_referential_barbell_Dax}, if $\beta$ is a self-referential barbell diffeomorphism associated with the loop $c$, then $\Psi_0(\beta) =\pm( [c] + [c^{-1}])$, so $\Psi_0$ is surjective. This proves the existence of $\Phi_0$ in Proposition \ref{prop_surj_to_Zinfty}.

To construct $\Phi_{1},\Phi_{2}$, we need more preparations. Recall that $\Sigma_{(0)} = \Sigma\times \{b_1\}\subset M$ is the image of $\mathscr{I}_0$. Also recall that $S_0=\{b'_0\}\times S^2$, where $b_0'$ is a point in the 2-handle of $\Sigma$, and the orientations of $\Sigma$ and $S^2$ are fixed.

\begin{Lemma}\label{lem: diff sends [Sigma] to [Sigma]}
If $f:M\to M$ is a diffeomorphism, then $f_{*}([\Sigma_{(0)}])=\pm[\Sigma_{(0)}]$.    
\end{Lemma}
\begin{proof}
We have 
\[
f_{*}\begin{pmatrix} [S_0]\\
[\Sigma_{(0)}]
\end{pmatrix}=\begin{pmatrix} a,b\\
c, d
\end{pmatrix}\begin{pmatrix} [S_0]\\
[\Sigma_{(0)}]
\end{pmatrix}.
\]
for some $a,b,c,d$ with $ad-bc=\pm 1$. The entry $b$ equals the degree of the composition map
\[
S^2\cong S_0\hookrightarrow M\xrightarrow{f} M \xrightarrow{\text{projection}} \Sigma,
\]
so it must be $0$. Hence $a,d=\pm 1$. Since $[\Sigma_0]\cdot [\Sigma_0]=0$, we have 
$0=f_{*}([\Sigma_0])\cdot f_{*}([\Sigma_0])=2cd.$
So $c=0$.
\end{proof}

Let $f:(M,b)\to (M,b)$ be a diffeomorphism. Let $i:(\Sigma,b_0)\to (M,b)$ be an embedding homotopic to $\mathscr{I}_0$ relative to $b_0$. 
In general, the map $f\circ i$ may not be homotopic to $\mathscr{I}_0$.
Recall that $\pi_1(\Sigma,b_0)$ is canonically isomorphic to $\pi_1(M,b)$. The map $f$ induces an automorphism on $\pi_1(M,b)$. By the Dehn--Nielson theorem, there exists a diffeomorphism $f_\dagger: (\Sigma,b_0)\to (\Sigma,b_0)$, which is unique up to isotopy relative to $b_0$, such that $f_\dagger$ induces the same automorphism as $f$ on $\pi$. 
Define $i^{f}: \Sigma\to M$ to be the composition
\begin{equation}
\label{eqn_def_i^f}
\Sigma\xrightarrow{f_\dagger^{-1}}\Sigma \xrightarrow{i}M \xrightarrow{f}M.
\end{equation}

\begin{Lemma}
The map $i^f$ is homotopic to $\mathscr{I}_0$ relative to $b_0$.
\end{Lemma}

\begin{proof}
    By definition, $i^f$ and $\mathscr{I}_0$ induce the same map on $\pi_1(\Sigma,b_0)\to \pi_1(M,b)$. Let $p$ be the projection from $M$ to $\Sigma$. By Lemma \ref{lem: diff sends [Sigma] to [Sigma]}, the maps $p\circ i^f$ and $p\circ \mathscr{I}_0$ both have degree $0$. Hence $i^f$ is homotopic to $\mathscr{I}_0$ relative to $b_0$ by Lemma \ref{lem_hom_in_F2}.
\end{proof}

\begin{Definition}\label{defi: Dax for MCG}
Define $\Psi_1: \Mcg_*(M)\to \bZ[\mathcal{C}]$ as follows.
Suppose $[f]\in \Mcg_*(M)$ is represented by a diffeomorphism $f:(M,b)\to (M,b)$, we define 
\[
\Psi_1([f])=\Dax(\mathscr{I}_0,\mathscr{I}_0^{f}).
\]
\end{Definition}
\begin{remark}
    $\Psi_1([f])$ is also equal to the relative Dax invariant of the un-parametrized embedded surfaces $\Sigma_{(0)}$ and $f(\Sigma_{(0)})$. 
\end{remark}

By definition, if $f\in\ker p_0\cong \Mcg_0(M)$, then $\mathscr{I}_0^f = f\circ \mathscr{I}_0$. So $\Psi_1$ is an extension of $\Psi_0$.
The map $\Psi_1$ is a map between sets. To understand the behavior of $\Psi_1$ under group operations, we prove the following result, which is an extension of Lemma \ref{lem_Dax_change_by_MCG0_*}.

\begin{Lemma}\label{lem: Dax invariant under general diff}
Let $i_1,i_2:(\Sigma,b_0)\to(M,b)$ be two embeddings that are homotopic to $\mathscr{I}_0$ relative to $b_0$. Let $f:(M,b)\to (M,b)$ be a diffeomorphism. Then we have 
\[
\Dax(i_1^f,i_2^f)=f_* \Dax(i_1,i_2),
\]
where $f_*$ denotes the induced action of $f_*:\pi_1(M,b)\to\pi_1(M,b)$ on $\mathcal{C}$. 
\end{Lemma}

\begin{proof} 
Let $g = (f_\dagger\times\id)^{-1}\circ f $. Then $g:(M,b)\to (M,b)$ induces the identity map on $\pi_1$. By Lemmas \ref{lem_Dax_equiv_diffeo} and \ref{lem_Dax_change_by_MCG0_*}, we have 
\begin{align*}
\Dax(i_1^f,i_2^f) & = \Dax(f\circ i_1(\Sigma), f\circ i_2(\Sigma)) \\
& = \Dax((f_\dagger\times\id)\circ g\circ i_1(\Sigma), (f_\dagger\times\id)\circ g\circ i_2(\Sigma))\\
& = (f_\dagger)_* \Dax(g\circ i_1(\Sigma), g\circ i_2(\Sigma))\\
& = (f_\dagger)_*\Dax(i_1,i_2) = f_* \Dax(i_1,i_2). \qedhere
\end{align*}
\end{proof}

By Lemma \ref{lem: Dax invariant under general diff}, given elements $[f_{1}],[f_{2}]\in \Mcg_{*}(M)$, we have 
\begin{equation}\label{eq: Dax nonadditive}
\begin{split}
\Psi_1([f_{1} f_{2}])=&\Dax(\mathscr{I}_0,\mathscr{I}_0^{f_1 f_2})\\
=&\Dax(\mathscr{I}_0,\mathscr{I}_0^{f_1})+\Dax(\mathscr{I}_0^{f_1},\mathscr{I}_0^{f_1 f_2})\\
=&\Dax(\mathscr{I}_0,\mathscr{I}_0^{f_1})+\Dax(\mathscr{I}_0^{f_1},(\mathscr{I}_0^{f_{2}})^{f_1})\\
=&\Dax(\mathscr{I}_0,\mathscr{I}_0^{f_1})+(f_{1})_*\Dax(\mathscr{I}_0,\mathscr{I}_0^{f_2})\\
=&\Psi_1([f_{1}])+(f_{1})_*\Psi_1([f_{2}]). 
\end{split}
\end{equation}
So $\Psi_1$ is not a group homomorphism. To construct the desired homomorphism $\Phi_1$ in Proposition \ref{prop_surj_to_Zinfty}, consider the action of $\Aut(\pi)\times \mathbb{Z}/2$ on $\mathcal{C}$, where the $\mathbb{Z}/2$ factor acts by taking inverses. We denote the set of orbits by $A$. 
\begin{Lemma}
$A$ is an infinite set.    
\end{Lemma}
\begin{proof} Let $g$ be an element in $\pi$ whose image in the abelization $\operatorname{Ab}(\pi)$ is primitive. Then for every $k\ge 1$, the image of $g^{k}$ in $\operatorname{Ab}(\pi)$ has divisibility $k$. Therefore, the elements $g^{k}$ $(k\geq 1)$ all belong to different orbits. 
\end{proof}
By taking the quotient map from $\mathcal{C}$ to $A$, we get a map 
$\phi_{1}:\mathbb{Z}[\mathcal{C}]\to \mathbb{Z}[A]$. We define $\Phi'_{1}: \Mcg_{*}(M)\to \mathbb{Z}[A]$ by 
$\Phi'_{1}([f]):=\phi_{1}\circ \Psi_1(f).$
\begin{Lemma} $\Phi'_{1}$ is a group homomorphism.
\end{Lemma}
\begin{proof}
Take $[f_1],[f_2]\in \Mcg_{*}(M)$. By \eqref{eq: Dax nonadditive}, we have 
\begin{align*}
\Phi'_{1}([f_{1}][f_{2}]) & =\phi_1 \Big(\Psi_1([f_{1}])+(f_{1})_*\Psi_1([f_{2}]) \Big)\\
& =\phi_1\circ \Psi_1([f_{1}])+\phi_{1} \circ \Psi_1([f_{2}]) =\Phi'_{1}([f_1])+\Phi'_{1}([f_{2}]).
\qedhere
\end{align*}
\end{proof}

\begin{Lemma}
\label{lem_Imm_Phi'_infinite_rk}
The image of $\Phi'_{1}$ is isomorphic to $\bZ^\infty$. 
\end{Lemma}
\begin{proof} 
For each $a\in A$, let $\delta_a\in \mathcal{C}$ be a lifting of $a$ to $\mathcal{C}$. Since $\Psi_0$ is surjective, we have $\delta_{a}+\delta_a^{-1}\in \Imm(\Psi_1)$ for all $a\in A$, so $2a=\phi_1(\delta_{a}+\delta_a^{-1})\in \Imm(\Phi_1')$. Therefore, $\Phi'_{1}$ is of infinite rank.
Since subgroups of $\bZ[A]$ are free abelian groups, we conclude that the image of $\Phi'_{1}$ is isomorphic to $\bZ^\infty$. 
\end{proof}

The existence of $\Phi_1$ in the statement of Proposition \ref{prop_surj_to_Zinfty} follows immediately from Lemma \ref{lem_Imm_Phi'_infinite_rk}.

Next, we consider the exact sequence 
\[
\cdots\to \pi_{1}(M,b)\xrightarrow{ \partial} \Mcg_*(M)\to \Mcg(M)\to 0
\]
coming from the fibration $\Diff_{*}(M)\hookrightarrow \Diff(M)\to M$. The existence of $\Phi_2$ is proved by the following lemma. 

\begin{Lemma}
\label{lem_Phi2'}
The restriction of $\Phi_{1}'$ to the image of $\partial$ is trivial. Therefore, $\Phi_{1}'$ descends to a homomorphism $\Phi_{2}':\Mcg(M)\to \mathbb{Z}[A]$ whose image is isomorphic to $\bZ^\infty$.
\end{Lemma}
\begin{proof}
Given $g\in \pi_1(M,b)$, represented by a loop $\gamma$ in $\Sigma$, we apply the isotopy extension theorem to get a point-pushing diffeomorphism $\widetilde{f}:\Sigma\to \Sigma$. Then $\partial([\gamma])=[f]$, where  $f=\operatorname{id}\times \widetilde{f}: M\to M$. It then follows  from the definitions that $\mathscr{I}_0^f = \mathscr{I}_0$, so $\Psi_{1}([f])=0$, and hence $\Phi_1'([f]) = 0$.  
\end{proof}

This finishes the proof of Proposition \ref{prop_surj_to_Zinfty}. Theorem \ref{thm: MCG infinitely generated} is a consequence of the above constructions and a simple observation as follows.

\begin{proof}[Proof of Theorem \ref{thm: MCG infinitely generated}]
    Consider the homomorphism $\Phi_2':\Mcg(M)\to \bZ[A]$ given by Lemma \ref{lem_Phi2'}. Then the image of a self-referential barbell diffeomorphism is $\pm([c]+[c^{-1}])$, where $c$ is the associated loop. Since self-referential barbell diffeomorphisms are pseudo-isotopic to the identity, the restriction of $\Phi_2'$ to $\Mcg_{PI}(M)$ is of infinite rank. 
\end{proof}

\section{Classification of surfaces with a dual sphere}
\label{sec_classification}
Recall that $\Sigma_{(d)}$ is the image of $\mathscr{I}_d$. 
If $\Sigma_1,\Sigma_2:\Sigma\to \Sigma\times S^2$ are smoothly embedded surfaces that are homotopic to $\Sigma_{(d)}$ (for suitable choices of parameterizations), we have the Dax invariant 
$\Dax(\Sigma_1,\Sigma_2)\in \mathbb{Z}[\mathcal{C}]^\sigma.$
Also recall that $b=(b_0,b_1)\in \Sigma\times S^2$ is a fixed base point.

The main result of this section is the following theorem.

\begin{Theorem}
\label{thm_classification_dual_surface}
Suppose $\Sigma_1,\Sigma_2\subset \Sigma\times S^2$ are smoothly embedded oriented surfaces that intersect $\{b_0\}\times S^2$ transversely and positively at one point, such that both surfaces agree with $\Sigma\times \{b_1\}$ on a neighborhood of $\{b_0\}\times S^2$. Suppose the projections from $\Sigma_{1}$, $\Sigma_2$ to $S^2$ both have degree $d$. Then
\begin{enumerate}
    \item There exists a tubular neighborhood $\nu(\{b_0\}\times S^2)$ of $\{b_0\}\times S^2$, such that the intersection of $\Sigma_1$, $\Sigma_2$, $\Sigma_{(d)}$ with $\nu(\{b_0\}\times S^2)$ are equal to the same 2-disk, and the intersection of the three surfaces with the complement of $\nu(\{b_0\}\times S^2)$ are homotopic relative to boundary in $\Sigma\times S^2\setminus \nu(\{b_0\}\times S^2)$.  
    \item $\Sigma_{1}$, $\Sigma_2$ are isotopic in $M$ if and only if 
$\Dax(\Sigma_{1}, \Sigma_2)= 0 \in \mathbb{Z}[\mathcal{C}]^\sigma.$
\end{enumerate}
\end{Theorem}

\begin{remark}
    In the statement of Theorem \ref{thm_classification_dual_surface} (2), we do not require that the isotopy preserves $\{b_0\}\times S^2$.
\end{remark}

Theorem \ref{thm_classification_dual_surface} implies the following classification result, which is a re-statement of Theorem \ref{intro thm: classification}. 
\begin{Corollary}
\label{cor_classification_dual_surface_deg_d}
Consider the set of embedded oriented surfaces $\Sigma_1\subset \Sigma\times S^2$ that are diffeomorphic to $\Sigma$, intersect $\{b_0\}\times S^2$ transversely and positively at one point, such that the projection from $\Sigma_1$ to $S^2$ has degree $d$. Let $\mathcal{F}_d$ be the equivalence classes of such surfaces up to isotopies in $M$. Then the map 
\begin{align*}
\mathcal{F}_d&\to \mathbb{Z}[\mathcal{C}]^\sigma\\
[\Sigma]&\mapsto \Dax(\Sigma, \Sigma_{(d)})
\end{align*}
is a bijection.
\end{Corollary}

Corollary \ref{cor_classification_dual_surface_deg_d} gives a classification result for embeddings of $\Sigma$ in $\Sigma\times S^2$ that are geometrically dual to $\{b_0\}\times S^2$ up to isotopy.

\begin{proof}[Proof of Corollary \ref{cor_classification_dual_surface_deg_d} from Theorem \ref{thm_classification_dual_surface}]
    Theorem \ref{thm_classification_dual_surface} and the additivity of the Dax invariant imply that the map $\mathcal{F}_d\to \mathbb{Z}[\mathcal{C}]^\sigma$ is injective. The self-referential tube construction shows that the map is surjective. 
\end{proof}

Before proving Theorem \ref{thm_classification_dual_surface}, we establish a technical lemma about the homotopy classes of maps between surfaces. The arguments of the proof of the next lemma will not be needed in the rest of the paper. 

\begin{Lemma}
\label{lem_map_punctured_Sigma}
    Let $\mathring{\Sigma} = \Sigma\setminus \inte(D^2)$, where $D^2$ is a smoothly embedded closed disk in $\Sigma$. Suppose $h: \mathring{\Sigma} \to \mathring{\Sigma}$ is a continuous map that restricts to the identity on the boundary, then $h$ is homotopic to a homeomorphism relative to $\partial \mathring{\Sigma}$.
\end{Lemma}

\begin{proof}
    We use an argument adapted from \cite{culler1981using}.
    Recall that $\ell$ denotes the genus of $\Sigma$. The surface $\mathring{\Sigma}$ has a handle decomposition with one $0$--handle and $2\ell$ 1--handles. Let $c_1,\dots,c_{2\ell}$ be the cocores of the 1--handles. Then the $c_j$'s are disjoint neatly embedded arcs in $\mathring{\Sigma}$, and cutting $\mathring{\Sigma}$ open along $c_1,\dots,c_{2\ell}$ yields a disk. 

    After a homotopy relative to $\partial \mathring{\Sigma}$, we may assume that $h$ is smooth and transverse to $c_1,\dots,c_{2\ell}$, so each $h^{-1}(c_j)$ is a compact $1$--submanifold of $\mathring{\Sigma}$ with exactly two boundary points. We may further assume that $h$ minimizes the total number of connected components of $\{h^{-1}(c_j)\}_{1\le j\le 2\ell}$ among all smooth maps that are homotopic to $h$ and transverse to $c_1,\dots,c_{2\ell}$. 

    \begin{Claim*}
    Each $h^{-1}(c_j)$ contains exactly one connected component. 
    \end{Claim*}    

    \begin{proof}[Proof of the claim]
    Assume the contrary, then for some $j$, the preimage $h^{-1}(c_j)$ contains at least one component that is a circle. Let $\hat c$ denote such a component. Since $h(\hat c)$ is contained in an arc, the restriction $h|_{\hat c}$ is null homotopic. 

    If $\hat c$ is non-trivial in $\mathring{\Sigma}$ (that is, $\hat c$ does not bound a disk on $\mathring{\Sigma}$), let $\mathring{\Sigma}'$ be obtained by cutting $\mathring{\Sigma}$ open along $\hat c$ and gluing the two new boundary components with two disks. Then the map $h$ lifts and extends to a map $h'$ from $\mathring{\Sigma}'$ to $\mathring{\Sigma}$, and $h'$ equals $h$ on the boundary. The genus of the non-closed component of $\mathring{\Sigma'}$ is strictly less than the genus of $\mathring{\Sigma}$.
    Let $\Sigma' = \mathring{\Sigma}' \cup_{\partial D^2} D^2$, then $h'$ extends to a map with degree $1$ from $\Sigma'$ to $\Sigma$, which yields a contradiction because a map from a closed surface with lower genus to a closed surface with higher genus must have degree zero. 

    If $\hat c$ is trivial in $\mathring{\Sigma}$, let $\hat D\subset \mathring{\Sigma}$ be the disk bounded by $\hat c$, let $u: \hat D \to c_j$ be a map such that $u|_{\hat c} = h|_{\hat c}$. Since $h(\hat c)\subset c_j$ and $\pi_1(c_j)$ is trivial, such a map $u$ exists. Let $h_1: \mathring{\Sigma}\to \mathring{\Sigma}$ be defined by 
    \[
    h_1(x) = 
    \begin{cases}
      u(x) \quad \text{ if } x \in \hat D,\\
      h(x) \quad \text{ if } x\notin \hat D.
    \end{cases}
    \]
    Since $\pi_2(\mathring{\Sigma})=0$, we have $u \simeq h|_{\hat D}$ rel $\partial \hat D$, so $h \simeq h_1$ rel $\partial \mathring{\Sigma}$. On the other hand, the map $h_1$ takes a small tubular neighborhood of $\hat D$ to a one-sided collar neighborhood of $c_j$. So after a further perturbation, we can find a map $h_2$ homotopic to $h_1$ rel $\partial \mathring{\Sigma}$ such that $\{h_2^{-1}(c_j)\}_{1\le j\le 2\ell}$ has fewer components, contradicting the assumption on $h$. This finishes the proof of the claim. 
    \end{proof}

    Back to the proof of the lemma. For each $j$, a tubular neighborhood of $h^{-1}(c_j)$ is diffeomorphic to $(-\epsilon,\epsilon)\times h^{-1}(c_j)$, and a tubular neighborhood of $c_j$ is diffeomorphic to $(-\epsilon,\epsilon)\times c_j$. We may choose the tubular neighborhoods so that the map $h$ is given by $\id_{(-\epsilon,\epsilon)}\times h|_{h^{-1}(c_j)}$ with respect to the above diffeomorphisms. Therefore, we may apply a further homotopy to $h$ so that it satisfies all the assumptions above and also
    \[
    h|_{h^{-1}(c_j)}: h^{-1}(c_j)\to c_j
    \]
    is a diffeomorphism for all $j$.

    Let $\Sigma^{(1)}$ be obtained by cutting open $\mathring{\Sigma}$ along $c_1,\dots,c_{2\ell}$. Then $\Sigma^{(1)}$ is a disk. 
    Let $\Sigma^{(2)}$ be obtained by cutting open $\mathring{\Sigma}$ along $\{h^{-1}(c_j)\}_{1\le j\le 2\ell}$. Then $\Sigma^{(2)}$ is a surface whose boundary is a circle and whose Euler number is $1$. Since $\Sigma^{(2)}$ cannot contain any closed components, it must be a disk. The map $h$ then lifts to a map $\tilde{h}: \Sigma^{(2)}\to \Sigma^{(1)}$ that is a homeomorphism on the boundary. Let $\hat{h}: \Sigma^{(2)}\to \Sigma^{(1)}$ be a homeomorphism that equals $\tilde{h}$ on the boundary. Since $\pi_2(D^2)=0$, we have $\hat{h}\simeq \tilde{h}$ relative to the boundary, so $h$ is homotopic to a homeomorphism relative to the boundary. 
\end{proof}

Now we can prove Theorem \ref{thm_classification_dual_surface} (1).

\begin{proof}[Proof of Theorem \ref{thm_classification_dual_surface} (1)]
    Let $i=i_1$ or $i_2$. By the assumptions, there exists an embedding $i:\Sigma \hookrightarrow M = \Sigma\times S^2$ such that $i=\mathscr{I}_0 = \mathscr{I}_d$ near $b_0$ and the image of $i$ is $\Sigma_j$. Since $i$ is geometrically dual to $\{b_0\}\times S^2$, after an isotopy of $\Sigma_j$ relative to a neighborhood of $\{b_0\}\times S^2$, we may assume that there exists a small 2-disk $D\subset \Sigma$ centered at $b_0$ such that the composition map
    \begin{equation}
\label{eqn_classification_dual_surface_(1)_proj_map}
        \Sigma \xhookrightarrow{i} M = \Sigma\times S^2 \xrightarrow{\text{projection}} \Sigma 
    \end{equation}
    is the identity on $D$ and maps $\Sigma\setminus \inte(D)$ to $\Sigma\setminus \inte(D)$. By Lemma \ref{lem_map_punctured_Sigma}, we know that the composition map \eqref{eqn_classification_dual_surface_(1)_proj_map} is homotopic to a homeomorphism relative to $D$. After a reparameterization of the domain of $i$, we may assume that the restriction of \eqref{eqn_classification_dual_surface_(1)_proj_map} to $\Sigma\setminus\inte(D)$ is homotopic to the identity map relative to the boundary via maps to $\Sigma\setminus\inte(D)$. 

    Note that the pointed homotopy class of maps $\langle \Sigma, S^2\rangle$ is classified by the degree because  $S^2$ is the 3-skeleton of $\mathbb{C}P^\infty \simeq K(\mathbb{Z},2)$. So the compositions
    \[
 \Sigma \xhookrightarrow{i} M = \Sigma\times S^2 \xrightarrow{\text{projection}} S^2,
    \]
    \[
 \Sigma \xhookrightarrow{\mathcal{I}_d} M = \Sigma\times S^2 \xrightarrow{\text{projection}} S^2 
    \]
    are homotopic relative to $D$, and the desired result is proved. 
\end{proof}

Theorem \ref{thm_classification_dual_surface} (1) shows that $\Dax(\Sigma_1,\Sigma_2)$ is well-defined. Now we start the proof of Theorem \ref{thm_classification_dual_surface} (2).

Let $\nu(\{b_0\}\times S^2)$ be the tubular neighborhood of $\{b_0\}\times S^2$ given by Theorem \ref{thm_classification_dual_surface} (1). Let 
\[
M_1' = M\setminus \nu(\{b_0\}\times S^2),\, M_2' = M_2\setminus \nu(\{b_0\}\times S^2).
\]
We may assume without loss of generality that $\nu(\{b_0\}\times S^2)$ is sufficiently small, so that both $M_1'$ and $M_2'$ are manifolds with corners. We smooth the corners and view them as smooth manifolds with boundaries. We introduce this new notation because the classification results of embedded disks in \cite{kosanovic2024homotopy,schwartz20214} require a geometrically dual $S^2$ on the boundary. 

Note that $\pi_1(M_1')\cong \pi_1(M_2')$, but their fundamental groups are not isomorphic to $\pi$. 
Let $\pi' = \pi_1(M_2')$, where the base point is taken on $D_0 = \Sigma_{(0)}\cap M_2'$. Then the inclusion map $M_2'\hookrightarrow M_2$ induces a surjection from $\pi'$ to $\pi$. 

Let $i_1,i_2:\Sigma\hookrightarrow \Sigma\times S^2$ be embeddings with images $\Sigma_1,\Sigma_2$ such that $i_1,i_2$ are homotopic to $\mathscr{I}_d$ relative to a neighborhood of $b_0$. The following lemma allows us to take an isotopy supported in $M_1'$ to match the 1-handles of $i_1,i_2$. 

\begin{Lemma}
\label{lem_M2'_match_1_handle}
    Suppose $i=i_1$ or $i_2$. There exists an ambient isotopy $P_t$ $(t\in[0,1])$ of $M$ supported in $M_1'$, such that $P_0=\id$, and $P_1\circ i$ equals $\mathscr{I}_d$ on $H_0\cup H_1$, and $P_1\circ i|_{H_2}$ is a neatly embedded disk in $M_2'$. 
\end{Lemma}

\begin{proof}
Let $\mathring{\Sigma} = i^{-1}(M_1') = \mathscr{I}_d^{-1}(M_1')$. Then $\mathring{\Sigma}$ is the complement of a small open disk in $\Sigma$. 
    By Theorem \ref{thm_classification_dual_surface} (1), we know that $i|_{\mathring{\Sigma}}$ and $\mathscr{I}_d|_{\mathring{\Sigma}}$ are homotopic relative to the boundary in $M_1'$. Since homotopy of $0$-- and $1$--handles in a $4$--manifold implies isotopy, there is an ambient isotopy supported in $M_1'$ that takes $i|_{H_0\cup H_1}$ to $\mathscr{I}_d|_{H_0\cup H_1}$. Taking a further isotopy to move the image of $H_2$ away from the neighborhood of the image of $H_1$, we may assume that the image of $H_2$ is a neatly embedded disk in $M_2'$.
\end{proof}

Now we show that one may take a further isotopy such that $i_1|_{H_2}$ is homotopic to $\mathscr{I}_d|_{H_2}$ in $M_2'$ relative to the boundary.

\begin{Lemma}
\label{lem_M2'_match_2_handle_homotopy}
Let $i:D^2\to M_2'$ be an embedding with the same boundary as $\mathscr{I}_0|_{H_2}$. Suppose the non-equivariant intersection number of $i$ with $D_0$ is $d\in \mathbb{Z}$. Then there exists a diffeomorphism $\varphi\in\Diff_\partial(M_2',\partial M_2')$, such that 
\begin{enumerate}
    \item The extension of $\varphi$ to $M_1'$ by the identity map is isotopic to the identity relative to $\partial M_1'$.
    \item The map $\varphi\circ i$ is homotopic to $\mathscr{I}_d|_{H_2}$ in $M_2'$ relative to the boundary. 
\end{enumerate}
\end{Lemma}

\begin{proof}
The proof is almost verbatim as the proof of Proposition \ref{prop_homotopy_after_barbell} by replacing $M_2$ with $M_2'$ and $\pi$ with $\pi'$. We give a sketch of the argument.

If $i_1,i_2:D^2\to M_2'$ are oriented embedded disks with the same boundary as $D_0$, push the boundary of $i_2$ along a non-zero vector of $T_{b_1}S^2$ and let $\lambda'(i_1,i_2) \in \mathbb{Z}[\pi']$ be the equivariant intersection number of $i_1$ and the pushed $i_2$. The equivariant intersection number only depends on the homotopy classes of $i_1,i_2$ relative to the boundary. 

Similarly, let $\bar{\pi}'$ be the quotient of $\pi'\setminus\{1\}$ by $g\sim g^{-1}$, then every embedding $i:D^2\to M_2'$ has an equivariant self-intersection number $\mu'(i)\in \mathbb{Z}[\bar{\pi}']$ defined by taking the intersection of $i$ with a pushing of $i$ near the boundary along a non-zero vector of $T_{b_1}S^2$. The value of $\mu'(i)$ only depends on the homotopy class of $i$ relative to the boundary. 

Let $\pi_2(M_2';\partial D_0)$ denote the set of maps of disks to $M_2'$ that have the same boundary as $\mathscr{I}_0|_{H_2}$ (which has the same boundary as $\mathscr{I}_d|_{H_2}$), up to homotopy relative to the boundary.
The same argument as Lemma \ref{lem_lambda_0} shows that 
\[
\lambda_0':= \lambda'(-,\mathscr{I}_0|_{H_2}): \pi_2(M_2';\partial D_0)\to  \mathbb{Z}[\pi']
\]
is a bijection.

A direct calculation with explicitly constructed representatives of elements of $\pi_2(M_2';\partial D_0)$ as in Lemma \ref{lem_relation_mu_lambda0} shows that the projection of $\lambda_0'(i)\in \mathbb{Z}[\pi']$ to $\mathbb{Z}[\bar\pi']$ equals $\mu'(i)$. Therefore, if $i$ is embedded, then the projection image of $\lambda_0'([i])$ in $\mathbb{Z}[\bar\pi']$ is zero.

Now consider the meridian-vertical barbell diffeomorphisms as defined in Section \ref{subsubsec_mer_vert_barbell}. Note that the barbells are contained in $M_2'$, and hence they can be viewed as diffeomorphisms of $M_2'$. The meridian-vertical barbell diffeomorphisms are isotopic to the identity in $\Diff_\partial(M_1')$ because the meridian sphere in the definition of the barbell is trivial in $M_1'$.

Similar to Lemma \ref{lem_change_of_lambda0_under_barbell} (see also Lemma \ref{lem_change_of_equi_intersection_by_meridian_vertical}), the action of a meridian-vertical barbell diffeomorphism changes the value of $\lambda_0([i])$ by
\begin{equation}
\label{eqn_g+-_change_homot_pi'}
g_+-g_--g_+^{-1}+g_-^{-1},
\end{equation}
where $g_+,g_-\in \pi'$ represent adjacent fundamental domains in the universal cover of  $\Sigma\setminus\inte(D^2)$. By the same argument as Proposition \ref{prop_homotopy_after_barbell}, elements of the form \eqref{eqn_g+-_change_homot_pi'} generate all expressions of the form $h-h^{-1}$ for $h\in \pi'\setminus\{1\}$. Therefore, after applying finitely many meridian-vertical barbell diffeomorphisms, the value of $\lambda_0'(i)$ can be changed to $d\in \mathbb{Z}[\pi']$, and hence the embedding $i$ becomes homotopic to $\mathscr{I}_d|_{H_2}$ in $M_2'$ relative to the boundary.
\end{proof}

Now we finish the proof of Theorem \ref{thm_classification_dual_surface}.

\begin{proof}[Proof of Theorem \ref{thm_classification_dual_surface} (2)]
The ``only if'' direction follows from the isotopy invariance of the Dax invariant. We prove the ``if'' direction here.

By Lemmas \ref{lem_M2'_match_1_handle} and \ref{lem_M2'_match_2_handle_homotopy}, we may apply ambient isotopies supported in $M_1'$ to $i_1$ and $i_2$, so that after the isotopy:
\begin{enumerate}
    \item $i_1|_{H_0\cup H_1} = i_2|_{H_0\cup H_1} = \mathscr{I}_d|_{H_0\cup H_1}$,
    \item $i_1|_{H_2} = i_2|_{H_2} = \mathscr{I}_d|_{H_2}$ are neatly embedded in $M_2'$ and are homotopic in $M_2'$ relative to the boundary. 
\end{enumerate}
Since the boundary of $\mathscr{I}_d|_{H_2}$ has a geometrically dual sphere in $\partial M_2'$ and $\pi' = \pi_1(M_2')$ has no $2$-torsion, by  \cite[Theorem 1.1]{kosanovic2024homotopy} (and also \cite[Theorem 0.2]{schwartz20214}), for $k=1,2$, the disk $[i_k|_{H_2}]$ is isotopic in $M_2'$ (rel. $\partial M_2'$) to a disk obtained by attaching $[\mathscr{I}_d|_{H_2}]$ with a collection of self-referential tubes associated with an element
\[
\alpha_k =\sum_{j=1}^{m_k} \epsilon_j^{(k)} g_j^{(k)} \in \mathbb{Z}[\pi'\setminus \{1\}],
\]
where $ \epsilon_j^{(k)}\in \{1,-1\}$ and $ g_j^{(k)}\in \pi'\setminus\{1\}$. 
See \cite[Definition 1.7]{schwartz20214} for the definition of self-referential disks. 

The Dax invariant of closed surfaces $\Dax(i_1,i_2)$ is equal to 
\[
\Dax(i_1,i_2)=\sum_{j=1}^{m_1} \epsilon_j^{(1)} \big([g_j^{(1)}]+[(g_j^{(1)})^{-1}]\big) - \sum_{j=1}^{m_2} \epsilon_j^{(2)} \big([g_j^{(2)}]+[(g_j^{(2)})^{-1}]\big) \in \mathbb{Z}[\mathcal{C}],
\]
where $[g_j^{(k)}]$ denotes the conjugation class of the image of $g_j^{(k)}$ in $\pi$. 
By the assumptions, $\Dax(i_1,i_2)=0$.

Note that an isotopy of the self-referential tube that moves the attaching point with $[\mathscr{I}_d|_{H_2}]$ changes the value of $g_j^{(k)}$ by conjugation, and an isotopy of the tube across $\{b_0\}\times S^2$ can change the value of  $g_j^{(k)}$ to every other value with the same image in $\pi = \pi_1(M_2)$. 

Since $\Dax(i_1,i_2)=0 \in \mathbb{Z}[\mathcal{C}]$, after isotopy of the self-referential tubes, we may assume that $\alpha_1=\alpha_2$. Therefore, the relative Dax invariant of $i_1|_{H_2}$ and $i_2|_{H_2}$ as disks in $M_2'$ is zero. By \cite[Theorem 1.1]{kosanovic2024homotopy}, this implies that after the previous isotopies, $i_1|_{H_2}$ and $i_2|_{H_2}$ become isotopic in $M_2'$ relative to the boundary. 
\end{proof}

\section{Non-isotopic symplectic structures on ruled surfaces}
\label{sec_symplectic}
 In this section, we prove the following result, which is a re-statement of Theorem \ref{thm_symplectic_intro}.
 
 \begin{Theorem}\label{thm: symplectic} Let $S^2\to X\xrightarrow{\pi}\Sigma$ be a ruled surface with $g(\Sigma)>0$. For every class $a\in H^2(X;\mathbb{R})$ with $a^2([X])>0$, the space 
 \[
\mathcal{S}_{a}:=\{\text{symplectic form $\omega$ on $X$ with $[\omega]=a$}\}
\]
 has infinitely many components. 
\end{Theorem}
 
 Let $p: L\to \Sigma$ be a complex line bundle. We assume $L$ is either trivial or $\langle c_{1}(L),[\Sigma]\rangle=1$. We use $X$ to denote the fiberwise one-point compactification of $L$. If $L$ is trivial, then $X$ is the product manifold $M=S^2\times \Sigma$ that we have been considering so far; when $\langle c_{1}(L),[\Sigma]\rangle=1$, the manifold $X$ is diffeomorphic to a twisted $S^2$-bundle over $\Sigma$, which will be denoted by $S^2\widetilde{\times}\Sigma$. 
 We use $\Sigma_0\hookrightarrow X$ to denote the $0$-section. For each $b\in \Sigma$, we use $F_{b}$ to denote the fiber over $b$. Fix a base point $b_0\in \Sigma$ and use $F$ to denote $F_{b_0}$. 
 
 The space $X$ has a complex structure $J_0$ such that $\Sigma_0$ and $\{F_{b}\}_{b\in \Sigma}$ are all $J_0$-holomorphic. One may construct a symplectic forms $\omega_{\mu}$ on $X$ for each $\mu>0$, such that $J_0$ is $\omega_\mu$--compatible and that 
 \[
 \operatorname{PD}[\omega_{\mu}]=[\Sigma_0]+\mu [F].
 \]
It is known that such symplectic forms are unique up to isotopy \cite{McDuffRuled}, and every symplectic form on $X$ is diffeomorphic to some $\omega_{\mu}$ up to rescaling by a constant \cite{LalondeMcDuff,LalondeMcDuffJ}. 

 The following proposition is a key ingredient in our proof of Theorem \ref{thm: symplectic}. It is independently known to Hind--Li \cite{HLpreprint}. 

\begin{Proposition}\label{proposition symplectic suface standard}
Let $\Sigma_1, \Sigma_2$ be embedded, connected symplectic surfaces in $(X,\omega_{\mu})$. Assume 
$[\Sigma_1]=[\Sigma_2]=[\Sigma_0]+k[F]$
for some $k\in \mathbb{Z}$. Then $\Sigma_1$ is smoothly isotopic to $\Sigma_2$ under a suitable parametrization. 
\end{Proposition}

\begin{proof}
By the adjunction formula, we have $g(\Sigma_{1})=g(\Sigma_{2})=g(\Sigma)$.
By the symplectic neighborhood theorem, there exist $\omega_{\mu}$-compatible almost complex structures $J_{1}, J_2$ on $X$ such that $\Sigma_{1}$ is $J_{1}$-holomorphic and $\Sigma_2$ is $J_2$-holomorphic. We connect $J_1,J_2$ by  a path $\{J_{t}\}_{1\leq t\leq 2}$ of $\omega_{\mu}$-compatible almost complex structures. Define the moduli space

\begin{equation}
    \cM_I:=\{(u,J_t)|u: S^2\to X, \partial_{J_t}u=0, [u]=[F]\}.
\end{equation}

Given $x\in X$ and any $t\in [1,2]$, there is a unique regular embedded $J_t$-holomorphic sphere passing through $x$  in the class $[F]$ (\cite{LalondeMcDuff}).  
Therefore, we have the following fibration
\begin{equation}\label{eq: moduli space bundle}
    \PSL(2,\mathbb{C})\to \cM_I\xrightarrow{\mathfrak{f}} \Sigma\times I.
\end{equation}
Here $I=[1,2]$, $\mathfrak{f}$ is the quotient by the natural $\PSL(2,\mathbb{C})$--action, which is equivalent to forgetting the parametrization of the curve $S^2$. Consider $\cM_S:=\cM_I\times_{\PSL(2,\mathbb{C})}S^2$.  There is a natural diffeomorphism $\varphi: \cM_S\rightarrow X\times [1,2]$ which respects the projection to $I$. Therefore, the associated bundle of (\ref{eq: moduli space bundle}) gives the smooth fiber bundle
\[
S^2\hookrightarrow X\times [1,2]\xrightarrow{\pi} \Sigma\times [1,2].
\]
For every $(b,t)\in \Sigma\times I$, the fiber $\pi^{-1}(b,t)$ is $J_{t}$-holomorphic.  

We use $\pi_{t}:X\to \Sigma$ $(t\in[1,2])$ to denote the bundle obtained by restricting $\pi$ to $X\times \{t\}$. By positivity of intersections for $J$-holomorphic curves, the surface $\Sigma_{1}$ (resp. $\Sigma_2$) is a smooth section of the bundle $\pi_{1}$ (resp. $\pi_{2}$).

Now we consider the composition  
\begin{equation}\label{eq: projection}
\Sigma\times [1,2]\xrightarrow{\operatorname{projection}} \Sigma\times \{1\}\hookrightarrow \Sigma\times [1,2].    
\end{equation}
If we pull back the bundle $\pi$ using this map, we will get the bundle 
\[
S^2\hookrightarrow X\times [1,2]\xhookrightarrow{\pi_{1}\times \operatorname{id}} \Sigma\times [1,2].
\]
On the other hand, the map (\ref{eq: projection})  is homotopic to $\operatorname{id}_{\Sigma\times \{1,2\}}$. By homotopy invariance of pull-back bundles, the bundles $\pi$ and $\pi_{1}\times \operatorname{id}$ are isomorphic. In other words, there exists a diffeomorphism $h: X\times [1,2]\to X\times [1,2]$ that fits into the commutative diagram 
\[
\xymatrix
{X\times [1,2]\ar[rd]_{\pi_{1}\times \operatorname{id}}\ar[rr]^{h}_{\cong}& & X\times [1,2]\ar[ld]^{\pi}\\
& \Sigma\times [1,2] &}
\]
By replacing $h$ with $h\circ( h|_{X\times \{1\}}\times \operatorname{id})^{-1}$, we may assume $h$ restricts to the identity map on $X\times \{1\}$. Let $f:X\to X$ be defined by restricting $h$ to $X\times \{2\}$. Then $f$ is isotopic to $\operatorname{id}_{X}$ and $\pi_2\circ f = \pi_1$. 

Note that $f(\Sigma_{1})$ and $\Sigma_2$ are both smooth sections of the bundle $X\xrightarrow{\pi_{2}}\Sigma$. Since $f(\Sigma_{1})$ is homologous to $\Sigma_{2}$, and homologous sections of a fixed bundle $S^2\hookrightarrow X\to \Sigma$ are smoothly isotopic, we see that $f(\Sigma_1)$ is smoothly isotopic to $\Sigma_{2}$. Hence $\Sigma_1$ is smoothly isotopic to $\Sigma_2$.
 \end{proof}

Next, we consider a self-referential barbell diffeomorphism $\beta_{p,\gamma}$ on $X$.
When $X= S^2\times \Sigma$, this has been introduced in Section \ref{sec_dax_invariant}. The same construction works when $X=S^2\widetilde{\times }\Sigma$. In both cases, we require that the associated loop $c(\gamma)$ is non-contractible.  We use $f:X\to X$ to denote  $\beta_{p,\gamma}$ and use $\Sigma_{m}$ to denote the embedded surface $f^m(\Sigma_0)$.

\begin{Lemma}\label{lem: sigma_m nonisotopic} For every positive integer $m$, the surface $\Sigma_m$ is not smoothly isotopic to $\Sigma_0$.
\end{Lemma}
\begin{proof}
When $X=\Sigma\times S^2$, by Lemma \ref{lem_self_referential_barbell_Dax}, we have  
\[
\Dax(\Sigma_0,\Sigma_m)=\pm m([c(\gamma)]+[c(\gamma)^{-1}])\neq 0.
\]
So $\Sigma_m$ is not isotopic to $\Sigma_0$. 

When $X=\Sigma\widetilde{\times }S^2$, we take a two-fold cover $\widetilde{\Sigma}\to \Sigma$ such that the preimage of $c(\gamma)$ has two components, denoted by $c_{1}(\gamma)$ and $c_{2}(\gamma)$. Let $\widetilde{\Sigma}_{m}$ be the preimage of $\Sigma$ under the covering map $\widetilde{\Sigma}\times S^2\to X$, then again by Lemma \ref{lem_self_referential_barbell_Dax}, we have  
\[
\Dax(\widetilde{\Sigma}_0,\widetilde{\Sigma}_m)=\pm m([c_{1}(\gamma)]+[c_{1}(\gamma)^{-1}]+[c_{2}(\gamma)]+[c_{2}(\gamma)^{-1}])\neq 0.
\]
Hence $\widetilde{\Sigma}_m$ is not isotopic to $\widetilde{\Sigma}_0$. So $\Sigma_m$ is not isotopic to $\Sigma_0$.
\end{proof}

\begin{proof}[Proof of Theorem \ref{thm: symplectic}]
After rescaling, we may assume $a=[\omega_{\mu}]$. For $m\geq 0$, consider the symplectic form $\omega_{m,\mu}:=(f^{m})^*(\omega_{\mu})\in \mathcal{S}_{a}.$
It suffices to show that for all $n< m$, the symplectic forms $\omega_{n,\mu}$ and $\omega_{m,\mu}$ are not in the same component of $\mathcal{S}_{a}$. By pulling back both symplectic forms via the diffeomorphism $f^{-n}$,  we may reduce to the case $n=0$. Now suppose $\omega_{m,\mu}$ and $\omega_{\mu}$ are in the same component of $\mathcal{S}_{a}$. Then, by Moser's argument, there exists a symplectomorphism $\varphi$ with respect to $\omega_\mu$, such that the diffeomorphism $f^{m}$ is smoothly isotopic to $\varphi$. Hence $\Sigma_{m}=f^{m}(\Sigma_0)$ is smoothly isotopic to $\varphi(\Sigma_0)$. Since  $\varphi(\Sigma_0)$ is a symplectic surface in $(X,\omega_{\mu})$, it is in turn smoothly isotopic to $\Sigma_0$ by Proposition \ref{proposition symplectic suface standard}. This implies that $\Sigma_m$ is smoothly isotopic to $\Sigma_0$, which is a contradiction to Lemma \ref{lem: sigma_m nonisotopic}.
\end{proof}

\bibliographystyle{amsalpha}
\bibliography{references}

@article{Igusa2021,
  title={Second obstruction to pseudoisotopy in dimension 3},
  author={Igusa, Kiyoshi},
  journal={arXiv preprint, arXiv:2112.08293},
  year={2021}
}

@article{Singh2021,
  title={Pseudo-isotopies and diffeomorphisms of 4-manifolds},
  author={Singh, Oliver},
  journal={arXiv preprint, arXiv:2111.15658},
  year={2021}
}

@article{FGHK2024,
  title={Grasper families of spheres in {$S^2\times D^2$} and barbell diffeomorphisms of {$S^1\times S^2\times I$}},
  author={Fern\'{a}ndez, Eduardo and  Gay, David T. and  Hartman, Daniel and  Kosanovi\'{c}, Danica},
  journal={arXiv preprint, arXiv:2412.07467},
  year={2024}
}

@article{kosanovic2024space,
  title={A space level light bulb theorem in all dimensions},
  author={Kosanovic, Danica and Teichner, Peter},
  journal={Comment. Math. Helv},
  volume={99},
  number={4},
  pages={799--860},
  year={2024}
}

@article{culler1981using,
  title={Using surfaces to solve equations in free groups},
  author={Culler, Marc},
  journal={Topology},
  volume={20},
  number={2},
  pages={133--145},
  year={1981},
  publisher={Pergamon}
}

@article{schwartz20214,
  title={A 4-dimensional light bulb theorem for disks},
  author={Schwartz, Hannah},
  journal={arXiv preprint arXiv:2109.13397},
  year={2021}
}

@inproceedings{dax1972etude,
  title={Etude homotopique des espaces de plongements},
  author={Dax, Jean-Pierre},
  booktitle={Annales scientifiques de l'{\'E}cole Normale Sup{\'e}rieure},
  volume={5},
  number={2},
  pages={303--377},
  year={1972}
}

@article{lin2025mapping,
  title={Mapping class groups of 4-manifolds with one-handles},
  author={Lin, Jianfeng and Xie, Yi and Zhang, Boyu},
  journal={arXiv preprint arXiv:2501.11821},
  year={2025}
}

@article{gabai2021self,
  title={Self-referential discs and the light bulb lemma},
  author={Gabai, David},
  journal={Commentarii Mathematici Helvetici},
  volume={96},
  number={3},
  year={2021}
}

@article{BG2023,
  title={On the automorphism groups of hyperbolic manifolds},
  author={ Budney, Ryan and  Gabai, David
},
  journal={arXiv preprint, arXiv:2303.05010},
  year={2023}
}

@article{Watanabe2020,
  title={Theta-graph and diffeomorphisms of some 4-manifolds},
  author={  Watanabe, Tadayuki
},
  journal={arXiv preprint,  arXiv:2005.09545},
  year={2020}
}

@article{BG2019,
  title={Knotted 3-balls in ${S}^4$},
  author={ Budney, Ryan and  Gabai, David
},
  journal={arXiv preprint, arXiv:1912.09029},
  year={2019}
}

@article {AB,
    AUTHOR = {Atiyah, M. F. and Bott, R.},
     TITLE = {The {Y}ang-{M}ills equations over {R}iemann surfaces},
   JOURNAL = {Philos. Trans. Roy. Soc. London Ser. A},
  FJOURNAL = {Philosophical Transactions of the Royal Society of London.
              Series A. Mathematical and Physical Sciences},
    VOLUME = {308},
      YEAR = {1983},
    NUMBER = {1505},
     PAGES = {523--615},
      ISSN = {0080-4614},
   MRCLASS = {14F05 (14D22 32G13 32L05 53C05 58E05 58E20 81E13)},
  MRNUMBER = {702806},
MRREVIEWER = {Martin A. Guest},
       DOI = {10.1098/rsta.1983.0017},
       URL = {https://doi.org/10.1098/rsta.1983.0017},
}

@article {Gabai20,
    AUTHOR = {Gabai, David},
     TITLE = {The 4-dimensional light bulb theorem},
   JOURNAL = {J. Amer. Math. Soc.},
  FJOURNAL = {Journal of the American Mathematical Society},
    VOLUME = {33},
      YEAR = {2020},
    NUMBER = {3},
     PAGES = {609--652},
      ISSN = {0894-0347},
   MRCLASS = {57K40 (57N35)},
  MRNUMBER = {4127900},
MRREVIEWER = {Sergey M. Finashin},
       DOI = {10.1090/jams/920},
       URL = {https://doi.org/10.1090/jams/920},
}

@article {KT24,
    AUTHOR = {Kosanovi\'{c}, Danica and Teichner, Peter},
     TITLE = {A new approach to light bulb tricks: disks in 4-manifolds},
   JOURNAL = {Duke Math. J.},
  FJOURNAL = {Duke Mathematical Journal},
    VOLUME = {173},
      YEAR = {2024},
    NUMBER = {4},
     PAGES = {673--721},
      ISSN = {0012-7094},
   MRCLASS = {57K45 (57R40 58D10)},
  MRNUMBER = {4734552},
MRREVIEWER = {Richard Stong},
       DOI = {10.1215/00127094-2023-0036},
       URL = {https://doi.org/10.1215/00127094-2023-0036},
}

@article{kosanovic2024homotopy,
  title={On homotopy groups of spaces of embeddings of an arc or a circle: the {D}ax invariant},
  author={Kosanovi{\'c}, Danica},
  journal={Transactions of the American Mathematical Society},
  volume={377},
  number={02},
  pages={775--805},
  year={2024}
}

@book{MSIntro,
	author = {McDuff, D. and Salamon, D.},
	edition = {Third},
	publisher = {Oxford University Press, Oxford},
	series = {Oxford Graduate Texts in Mathematics},
	title = {Introduction to symplectic topology},
	year = {2017}}

@article{Gromov85,
	author = {Gromov, M.},
	journal = {Inventiones mathematicae},
	number = {2},
	pages = {307--347},
	publisher = {Springer-Verlag Berlin/Heidelberg},
	title = {Pseudo holomorphic curves in symplectic manifolds},
	volume = {82},
	year = {1985}}

@incollection {LalondeMcDuffJ,
    AUTHOR = {Lalonde, Fran\c{c}ois and McDuff, Dusa},
     TITLE = {{$J$}-curves and the classification of rational and ruled
              symplectic {$4$}-manifolds},
 BOOKTITLE = {Contact and symplectic geometry ({C}ambridge, 1994)},
    SERIES = {Publ. Newton Inst.},
    VOLUME = {8},
     PAGES = {3--42},
 PUBLISHER = {Cambridge Univ. Press, Cambridge},
      YEAR = {1996},
   MRCLASS = {57R15 (57N13 57R57 58D10)},
  MRNUMBER = {1432456},
MRREVIEWER = {Margaret F. Symington},
}

@article {LalondeMcDuff,
    AUTHOR = {Lalonde, Fran\c{c}ois and McDuff, Dusa},
     TITLE = {The classification of ruled symplectic {$4$}-manifolds},
   JOURNAL = {Math. Res. Lett.},
  FJOURNAL = {Mathematical Research Letters},
    VOLUME = {3},
      YEAR = {1996},
    NUMBER = {6},
     PAGES = {769--778},
      ISSN = {1073-2780},
   MRCLASS = {57R15 (57R57)},
  MRNUMBER = {1426534},
MRREVIEWER = {Ricardo F. Vila Freyer},
       DOI = {10.4310/MRL.1996.v3.n6.a5},
       URL = {https://doi.org/10.4310/MRL.1996.v3.n6.a5},
}

@article{wang2024note,
  title={A note on a nontrivial connected component of the space of symplectic structures},
  author={Wang, Luya},
  journal={arXiv preprint arXiv:2410.16216},
  year={2024}
}

@article{HLpreprint,
  title={Private communications},
  author={Hind, Richard and Li, Jun},
}

@article {SalamonSurvey,
    AUTHOR = {Salamon, Dietmar},
     TITLE = {Uniqueness of symplectic structures},
   JOURNAL = {Acta Math. Vietnam.},
  FJOURNAL = {Acta Mathematica Vietnamica},
    VOLUME = {38},
      YEAR = {2013},
    NUMBER = {1},
     PAGES = {123--144},
      ISSN = {0251-4184,2315-4144},
   MRCLASS = {53D35},
  MRNUMBER = {3098204},
MRREVIEWER = {Philippe\ Rukimbira},
       DOI = {10.1007/s40306-012-0004-x},
       URL = {https://doi.org/10.1007/s40306-012-0004-x},
}

@book {TaubesGrSW,
    AUTHOR = {Taubes, Clifford Henry},
     TITLE = {Seiberg {W}itten and {G}romov invariants for symplectic
              {$4$}-manifolds},
    SERIES = {First International Press Lecture Series},
    VOLUME = {2},
      NOTE = {Edited by Richard Wentworth},
 PUBLISHER = {International Press, Somerville, MA},
      YEAR = {2000},
     PAGES = {vi+401},
      ISBN = {1-57146-061-6},
   MRCLASS = {53D45 (32Q65 53D35 57R17 57R57)},
  MRNUMBER = {1798809},
MRREVIEWER = {Ignasi Mundet-Riera},
}

@article {Abreu,
    AUTHOR = {Abreu, Miguel and McDuff, Dusa},
     TITLE = {Topology of symplectomorphism groups of rational ruled
              surfaces},
   JOURNAL = {J. Amer. Math. Soc.},
  FJOURNAL = {Journal of the American Mathematical Society},
    VOLUME = {13},
      YEAR = {2000},
    NUMBER = {4},
     PAGES = {971--1009},
      ISSN = {0894-0347},
   MRCLASS = {57R15 (53D35 57S05)},
  MRNUMBER = {1775741},
MRREVIEWER = {Margaret F. Symington},
       DOI = {10.1090/S0894-0347-00-00344-1},
       URL = {https://doi.org/10.1090/S0894-0347-00-00344-1},
}

@article {McDuffRuled,
    AUTHOR = {McDuff, Dusa},
     TITLE = {The structure of rational and ruled symplectic
              {$4$}-manifolds},
   JOURNAL = {J. Amer. Math. Soc.},
  FJOURNAL = {Journal of the American Mathematical Society},
    VOLUME = {3},
      YEAR = {1990},
    NUMBER = {3},
     PAGES = {679--712},
      ISSN = {0894-0347},
   MRCLASS = {58F05 (53C15 57R50 58C10)},
  MRNUMBER = {1049697},
MRREVIEWER = {Jean-Claude Sikorav},
       DOI = {10.2307/1990934},
       URL = {https://doi.org/10.2307/1990934},
}

@incollection {TaubesGr,
    AUTHOR = {Taubes, Clifford H.},
     TITLE = {{$\rm Gr=SW$}: counting curves and connections [1761081]},
 BOOKTITLE = {Seiberg {W}itten and {G}romov invariants for symplectic
              4-manifolds},
    SERIES = {First Int. Press Lect. Ser.},
    VOLUME = {2},
     PAGES = {275--401},
 PUBLISHER = {Int. Press, Somerville, MA},
      YEAR = {2000},
   MRCLASS = {53D45 (32Q65 53D35 57R17 57R57)},
  MRNUMBER = {1798140},
}

@article {TaubesSW,
    AUTHOR = {Taubes, Clifford Henry},
     TITLE = {The {S}eiberg-{W}itten and {G}romov invariants},
   JOURNAL = {Math. Res. Lett.},
  FJOURNAL = {Mathematical Research Letters},
    VOLUME = {2},
      YEAR = {1995},
    NUMBER = {2},
     PAGES = {221--238},
      ISSN = {1073-2780},
   MRCLASS = {57R57 (57R15 58D15 58D29)},
  MRNUMBER = {1324704},
MRREVIEWER = {Jim A. Bryan},
       DOI = {10.4310/MRL.1995.v2.n2.a10},
       URL = {https://doi.org/10.4310/MRL.1995.v2.n2.a10},
}

@article {Ruan94,
    AUTHOR = {Ruan, Yongbin},
     TITLE = {Symplectic topology on algebraic {$3$}-folds},
   JOURNAL = {J. Differential Geom.},
  FJOURNAL = {Journal of Differential Geometry},
    VOLUME = {39},
      YEAR = {1994},
    NUMBER = {1},
     PAGES = {215--227},
      ISSN = {0022-040X},
   MRCLASS = {14J15 (14J29 14J30 57R15 57R57 58F05)},
  MRNUMBER = {1258920},
MRREVIEWER = {I. Dolgachev},
       URL = {http://projecteuclid.org/euclid.jdg/1214454682},
}

@book {CieliebakHprinciple,
    AUTHOR = {Cieliebak, K. and Eliashberg, Y. and Mishachev, N.},
     TITLE = {Introduction to the {$h$}-principle},
    SERIES = {Graduate Studies in Mathematics},
    VOLUME = {239},
      NOTE = {Second edition},
 PUBLISHER = {American Mathematical Society, Providence, RI},
      YEAR = {2024},
     PAGES = {xvii+363},
}

@article {McEx,
    AUTHOR = {McDuff, Dusa},
     TITLE = {Examples of symplectic structures},
   JOURNAL = {Invent. Math.},
  FJOURNAL = {Inventiones Mathematicae},
    VOLUME = {89},
      YEAR = {1987},
    NUMBER = {1},
     PAGES = {13--36},
      ISSN = {0020-9910,1432-1297},
   MRCLASS = {58F05 (53C57 57R50 58C10)},
  MRNUMBER = {892186},
MRREVIEWER = {Jean-Claude\ Sikorav},
       DOI = {10.1007/BF01404672},
       URL = {https://doi.org/10.1007/BF01404672},
}

\end{document}